\documentclass[10pt]{amsart}
\usepackage{amsmath}
\usepackage{amssymb}
\usepackage{amsfonts}
\usepackage{pxfonts}
\usepackage{amsthm}
\usepackage{amsxtra}
\usepackage{rotating}
\usepackage{float}
\usepackage[all,web,arc,poly]{xy}
\UseCrayolaColors
\usepackage{epsfig}
\usepackage{array}
\usepackage{pifont}
\usepackage{multirow}
\usepackage{graphics}

\theoremstyle{plain}
\newtheorem{theorem}{Theorem}[section]
\newtheorem{corollary}{Corollary}[section]
\newtheorem{lemma}{Lemma}[section]
\newtheorem{proposition}{Proposition}[section]

\theoremstyle{definition}

\theoremstyle{remark}
\newtheorem{remark}{Remark}[section]

\newcommand{\C}{\mathbb C}
\newcommand{\R}{\mathbb R}
\newcommand{\Z}{\mathbb Z}
\newcommand{\N}{\mathbb N}
\newcommand{\T}{\mathbb T}

\newcommand{\DySt}{\displaystyle}
\newcommand{\ScSt}{\scriptstyle}

\newcommand{\half}{
        {\lower0.00ex\hbox{\raise.6ex\hbox{\the\scriptfont0 1}
                           \kern-.5em\slash\kern-.1em\lower.45ex
                                     \hbox{\the\scriptfont0 2}}}}
\newcommand{\quarter}{
        {\lower0.00ex\hbox{\raise.6ex\hbox{\the\scriptfont0 1}
                           \kern-.5em\slash\kern-.1em\lower.45ex
                                     \hbox{\the\scriptfont0 4}}}}
\newcommand{\tquarter}{
        {\lower0.00ex\hbox{\raise.6ex\hbox{\the\scriptfont0 3}
                           \kern-.5em\slash\kern-.1em\lower.45ex
                                     \hbox{\the\scriptfont0 4}}}}
\newcommand{\eighth}{
        {\lower0.00ex\hbox{\raise.6ex\hbox{\the\scriptfont0 1}
                           \kern-.5em\slash\kern-.1em\lower.45ex
                                     \hbox{\the\scriptfont0 8}}}}
\newcommand{\othird}{
        {\lower0.00ex\hbox{\raise.6ex\hbox{\the\scriptfont0 1}
                           \kern-.5em\slash\kern-.1em\lower.45ex
                                     \hbox{\the\scriptfont0 3}}}}

\begin{document}

\title[]
{Bi-orthogonal systems on the unit circle, Regular Semi-Classical Weights and Integrable Systems - II.}

\author{N.S.~Witte}
\address{Department of Mathematics and Statistics,
University of Melbourne,Victoria 3010, Australia}
\email{\tt n.witte@ms.unimelb.edu.au}

\begin{abstract}
We derive the Christoffel-Geronimus-Uvarov transformations of a system of bi-orthogonal 
polynomials and associated functions on the unit circle, that is to say the 
modification of the system corresponding to a rational modification of the 
weight function. 
In the specialisation of the weight function to the regular semi-classical
case with an arbitrary number of regular singularities $ \{z_1, \ldots, z_M \} $ 
the bi-orthogonal system is known to be isomonodromy preserving with respect to
deformations of the singular points. If the zeros and poles of the 
Christoffel-Geronimus-Uvarov factors coincide with the singularities then we have the 
Schlesinger transformations of this isomonodromic system. Compatibility of the 
Schlesinger transformations with the other structures of the system - the recurrence 
relations, the spectral derivatives and deformation derivatives is explicitly 
deduced. Various forms of Hirota-Miwa equations are derived for the $ \tau $-functions 
or equivalently Toeplitz determinants of the system.
\end{abstract}

\maketitle

\section{Motivations}
\setcounter{equation}{0}
The unitary group $U(N)$ with Haar (uniform) measure
has eigenvalue probability density function (see e.g.~\cite[Chapter 2]{rmt_Fo})
\begin{equation}
   \frac{1}{(2 \pi )^N N!} \prod_{1 \le j < k \le N} | z_k - z_j |^2,
 \quad z_l \coloneqq e^{i \theta_l} \in \T, \quad \theta_l \in (-\pi,\pi] ,
\label{UN_Haar}
\end{equation}
where $ \T = \{z \in \C: |z|=1 \} $.
One of the motivations of our study is to characterise averages over $U \in U(N)$ of class functions $ w(U) $ 
(i.e. symmetric functions of the eigenvalues of $ U $ only) which have
the factorization property $ \prod_{l=1}^N w(z_l) $ for
$ \{z_1,\dots, z_N \} \in {\rm Spec}(U) $. Such functions $ w(z) $ can be interpreted as 
weights. Introducing the Fourier components $\{w_l\}_{l\in \Z}$ of the weight
$ w(z) = \sum_{l=-\infty}^\infty w_l z^l $, due to
the well known Heine identity \cite{ops_Sz}
\begin{equation}
  \Big \langle \prod_{l=1}^N w(z_l) \Big \rangle_{U(N)} =
  \det[ w_{i-j} ]_{i,j=0,\dots,N-1},
\label{Haar_avge}
\end{equation}
we are equivalently studying Toeplitz determinants.

Such averages over the unitary group are ubiquitous in many applications to
mathematical physics,  
in particular the gap probabilities and characteristic polynomial averages in the 
circular ensembles of random 
matrix theory \cite{rmt_Fo},\cite{AvM_2002},\cite{FW_2004},
the spin-spin correlations of the planar Ising model \cite{McCW_1973},\cite{JM_1980},
the density matrix of a system of impenetrable bosons on the ring 
\cite{FFGW_2002a} and probability distributions for various classes of non-intersecting 
lattice path problems \cite{Fo_2003}.
All of these applications can be subsumed within a generic class of weights, termed 
the regular semi-classical class and its degenerations.
The simplest non-trivial example of a regular semi-classical weight is
\begin{equation}
   w(z) = t^{-\mu}z^{-\omega-\mu}(1+z)^{2\omega_1}(1+tz)^{2\mu} \;
   \begin{cases}
     1, \quad \theta \in (-\pi,\pi-\phi) \cr
     1-\xi, \quad \theta \in (\pi-\phi,\pi]
   \end{cases} ,
\label{M=3_scw}
\end{equation}
where $ \mu $, $ \omega=\omega_1+i\omega_2 $ are complex parameters and $ \xi $,
$ t=e^{i\phi} $ are complex variables. 
It has been previously shown \cite{FW_2004} that the average (\ref{Haar_avge}) 
with weight (\ref{M=3_scw}) can, as a function of $ t $, be characterised as a 
$ \tau$-function for the Painlev\'e VI system. The examples mentioned above then 
are characterised by solutions to the sixth Painlev\'e equation with the parameters 
given in Table \ref{M=3_table}.
One of our primary motivations for our study is to deduce recurrence relations that
can be used to characterise (also in an algorithmic sense) this class of unitary
group averages.

\begin{table}[h]
\begin{tabular}{cccc|c}\label{M=3_table}
$ \mu $	& $ \omega_1 $	& $ \omega_2 $	& $ \xi $	& Interpretation \\
\hline
$ 0 $ & $ 0 $ & $ 0 $ & $ \xi $ & 
\begin{minipage}[c][2cm][c]{8cm}
{Generating function for the probability of exactly $ k $ eigenvalues of a random unitary matrix with 
phases in the sector $ (\pi-\phi,\pi] $, Dyson CUE}\end{minipage} \\
\hline
$ k $ & $ 0 $ & $ 0 $ & $ 0 $ & 
\begin{minipage}[c][1.5cm][c]{8cm}
{$k$-th moment of the characteristic polynomial for random unitary matrices}\end{minipage} \\
\hline
$ 1/4 $	& $ 1/4 $ & $ -i/2 $ & $ 0 $ & 
\begin{minipage}[c][1.5cm][c]{8cm}
{Diagonal spin-spin correlations of square lattice Ising model}
\end{minipage} \\
\hline
$ 1/2 $	& $ 0 $	& $ 0 $	& $ 2 $ & 
\begin{minipage}[c][1.5cm][c]{8cm}
{Density matrix of impenetrable Bose gas confined on a circle}\end{minipage}  \\
\hline
\end{tabular}
\vskip0.5cm
\caption{Examples of the regular semi-classical weight (\ref{M=3_scw}) occurring 
in random matrix theory and condensed matter physics}
\end{table}

However one notices that members of this class of weights are not necessarily 
positive or even real valued, $ \overline{w(z)} \neq w(z) $ for $ z \in \T $, where 
the bar denotes the complex conjugate, and consequently the Toeplitz matrices are 
non-Hermitian $ \bar{w}_n \neq w_{-n} $. Systems of orthogonal polynomials 
constructed with such weights defined on the unit circle are naturally split
into a bi-orthogonal pair of polynomials (see (\ref{orthog:a}),(\ref{orthog:b})).
Following the pioneering works on orthogonal polynomial systems (OPS) on the unit 
circle by Szeg{\H o} and Geronimus (see the monographs \cite{ops_Sz}, \cite{Ge_1961} and the 
contemporary state of affairs in Simon \cite{Simon_I_2005},\cite{Simon_II_2005}) the 
study of bi-orthogonal polynomial systems (BOPS) was initiated by Baxter \cite{Ba_1961}
motivated in part by the desire to analyse non-hermitian Toeplitz matrices and applications
in probability theory.
In an independent development a study of systems of orthogonal Laurent polynomials
was begun by Jones and Thron \cite{JT_1982} which has its origins in the theory of a particular type of
analytic continued fraction, the $T$-fractions, in multipoint Pad\'e approximation and the strong forms
of the Stieltjes and Hamburger moment problems. In subsequent work \cite{JTN_1984},
\cite{HvR_1986} it was realised that these orthogonal Laurent polynomial systems
(OLPS) are equivalent to a pair of bi-orthogonal polynomials, and furthermore that
these were found to be equivalent to Baxter's BOPS in \cite{Pa_1985}. In fact one can construct
a system of bi-orthogonal Laurent polynomials through the defining recurrence relations,
with the consequence that our conclusions apply equally to the 
bi-orthogonal Laurent polynomial systems or BOLPS.

A systematic study of BOPS on the unit circle and the 
regular semi-classical class was initiated in \cite{FW_2004a} where a number
of basic structures were laid out - the well-known recurrence relation, 
the spectral derivatives of the system and the deformation derivatives with
respect to the free singular points of the weight. 
Subsequently a similar program in the context of BOLPS was undertaken in \cite{BG_2007}.
The present study is a 
continuation of the program \cite{FW_2004a} and here the central theme is the rational 
modification of the weight and the consequences for the system and its structures.

Why do we consider rational modification? The answer to this question is two-fold -
there are direct applications of such formulae and because a given regular
semi-classical weight should be seen as a member of a family of weights whose 
exponents differ by integer increments or decrements and such systems are
linked by Schlesinger transformations. 
The historical literature has focused on the rational modification of weights
defined on the real line or the interval and the formulae expressing the
new orthogonal polynomials in terms of the original ones were derived -
see the work of Uvarov \cite{Uv_1959}, \cite{Uv_1969} following the pioneering
study of Christoffel \cite{Ch_1858}. There has emerged a modern period of interest in 
this subject, see for example \cite{Zh_1997}, \cite{SVZ_1997}, and \cite{BM_2004}.

Formulae for the polynomials orthogonal with respect to weights changed by
polynomial multiplication and division written in terms of the original system
were extended to those defined on the unit circle in the
works \cite{GM_1987}, \cite{GM_1991}, \cite{IR_1992} and \cite{GM_1993}.  
The extension of such formulae to BOPS on the unit circle was undertaken,
however in an incomplete way, in \cite{Ru_1994}.
The most general rational modification of a weight $ w(z) $ on the unit circle 
takes the form
\begin{align}
   w\left[\begin{array}{cc} \ScSt K; & \ScSt K^* \\ \ScSt L; & \ScSt L^* \end{array};z \right] 
  & = w\left[\begin{array}{ccccccc}
             \ScSt \alpha_1 & \ldots & \ScSt \alpha_K &; & \ScSt \alpha^*_1 & \ldots & \ScSt \alpha^*_{K^*} \\
             \ScSt \beta_1  & \ldots & \ScSt \beta_L  &; & \ScSt \beta^*_1  & \ldots & \ScSt \beta^*_{L^*} 
            \end{array};z \right] ,
  \label{C-Uxfm:a} \\
  & = \frac{\prod^{K}_{i=1}(z-\alpha_i)\prod^{K^*}_{j=1}(1-\alpha^*_jz^{-1})}
          {\prod^{L}_{k=1}(z-\beta_k)\prod^{L^*}_{l=1}(1-\beta^*_lz^{-1})} w(z) ,
      \qquad \alpha_i,\alpha^*_j,\beta_k,\beta^*_l \in \C .
  \label{C-Uxfm:b}
\end{align}
Whilst this contains redundancies it manifestly exhibits interesting symmetries 
and is rendered real and positive under simple conditions on the parameters.
In addition it is this case, or more precisely the polynomial form, that arises
in studies on the symmetrised model of last passage percolation and non-intersecting
lattice paths, e.g. see Formulae (2.3), (3.9) and (6.12) of \cite{FR_2007} which
were first derived in \cite{BR_2001}.  
The general rational modification (\ref{C-Uxfm:b}) is also fundamental in the 
calculation of the averages of ratios of characteristic polynomials (which are 
an example of a class function) with respect to the Haar measure
(here the original weight is trivial, $ w(z)=1 $). In this situation we are, in effect, 
considering the average (\ref{Haar_avge}) with respect to a rational weight function. 
A contemporary interest in such averages arises from analytic number theory where 
the average (\ref{Haar_avge}) of (\ref{C-Uxfm:b}) with $ w(z)=1 $ is conjectured
to reproduce the integral of ratios of the Riemann zeta function along the critical
line \cite{CFS_2005}, \cite{CFZ_2007}, \cite{CFZ_2005}, in an asymptotic sense as
$ N \to \infty $. This conjecture has been enlarged to apply to similar averages 
of certain families of $L$-functions, which now involve the calculation of averages of
ratios of characteristic polynomials over the Haar measure for the compact classical
groups $ Sp(N) $, $ O^+(2N) $ and $ O^+(2N+1) $. 

\begin{sideways}
\begin{minipage}[c][12cm][c]{20cm}{
Our penultimate result gives the bi-orthogonal polynomials $ \Phi_n, \Phi^*_n $ 
and their associated functions $ \Xi_n, \Xi^*_n $
for the modified weight defined by (\ref{C-Uxfm:b}) as a 
$ (K+K^*+L+L^*+1)\times (K+K^*+L+L^*+1) $ determinant with a block structure 
whose elements are the bi-orthogonal polynomials $ \phi_n, \phi^*_n $ and their associated
functions $ \xi_n, \xi^*_n $ 
(see (\ref{orthog:a}), (\ref{orthog:b}), (\ref{eps:a}), (\ref{eps:b}), (\ref{eps:c}) for their definitions)
corresponding to the original weight. For example the bi-orthogonal polynomial $ \Phi_n $ is given by
\begin{multline}
  \prod^{K}_{i=1}(z-\alpha_i)\prod^{K^*}_{j=1}(1-\alpha^*_j z^{-1})\,\Phi_n(z) \\ 
  = C\cdot
  \det \begin{pmatrix}
	[z^{L-c}\phi_{n-L^*-L+c}(z)]_{c=0,..,L-1}	&
	[\phi_{n-L^*+c}(z)]_{c=0,..,L^*-1} &
	[z^{-c}\phi_{n+c}(z)]_{c=0,..,K^*} &
	[z^{-K^*}\phi_{n+K^*+c}(z)]_{c=1,..,K} \cr
	[\alpha_r^{L-c}\phi_{n-L^*-L+c}(\alpha_r)]_{{r=1,..,K}\atop{c=0,..,L-1}}	&
	[\phi_{n-L^*+c}(\alpha_r)]_{{r=1,..,K}\atop{c=0,..,L^*-1}} &
	[\alpha_r^{-c}\phi_{n+c}(\alpha_r)]_{{r=1,..,K}\atop{c=0,..,K^*}} &
	[\alpha_r^{-K^*}\phi_{n+K^*+c}(\alpha_r)]_{{r=1,..,K}\atop{c=1,..,K}} \cr
	[(\alpha^*_r)^{L-c}\phi_{n-L^*-L+c}(\alpha^*_r)]_{{r=1,..,K^*}\atop{c=0,..,L-1}}	&
	[\phi_{n-L^*+c}(\alpha^*_r)]_{{r=1,..,K^*}\atop{c=0,..,L^*-1}} &
	[(\alpha^*_r)^{-c}\phi_{n+c}(\alpha^*_r)]_{{r=1,..,K^*}\atop{c=0,..,K^*}} &
	[(\alpha^*_r)^{-K^*}\phi_{n+K^*+c}(\alpha^*_r)]_{{r=1,..,K^*}\atop{c=1,..,K}} \cr
	[\beta_r^{L-c}\xi_{n-L^*-L+c}(\beta_r)]_{{r=1,..,L}\atop{c=0,..,L-1}}	&
	[\xi_{n-L^*+c}(\beta_r)]_{{r=1,..,L}\atop{c=0,..,L^*-1}} &
	[\beta_r^{-c}\xi_{n+c}(\beta_r)]_{{r=1,..,L}\atop{c=0,..,K^*}} &
	[\beta_r^{-K^*}\xi_{n+K^*+c}(\beta_r)]_{{r=1,..,L}\atop{c=1,..,K}} \cr
	[(\beta^*_r)^{L-c}\xi_{n-L^*-L+c}(\beta^*_r)]_{{r=1,..,L^*}\atop{c=0,..,L-1}}	&
	[\xi_{n-L^*+c}(\beta^*_r)]_{{r=1,..,L^*}\atop{c=0,..,L^*-1}} &
	[(\beta^*_r)^{-c}\xi_{n+c}(\beta^*_r)]_{{r=1,..,L^*}\atop{c=0,..,K^*}} &
	[(\beta^*_r)^{-K^*}\xi_{n+K^*+c}(\beta^*_r)]_{{r=1,..,L^*}\atop{c=1,..,K}} \cr
       \end{pmatrix} .
\label{KK*LL*shift}
\end{multline}
The constant $ C $ is independent of $ z $. The other components of the bi-orthogonal
system are given by the above formula with the following substitutions: 
$ \prod^{K}_{i=1}(z-\alpha_i)\prod^{K^*}_{j=1}(1-\alpha^*_j z^{-1})\,\Phi^*_n(z) $ is given by
(\ref{KK*LL*shift}) with $ C \mapsto C^* $, $ \phi_n \mapsto \phi^*_n $ in the first row, and 
second and third block rows and $ \xi_n \mapsto \xi^*_n $ in the last
two block rows; 
$ \prod^{L}_{i=1}(z-\beta_i)\prod^{L^*}_{j=1}(1-\beta^*_j z^{-1})\,\Xi_n(z) $ is given by
(\ref{KK*LL*shift}) with $ \phi_n(z) \mapsto \xi_n(z) $ in the first row;
and
$ \prod^{L}_{i=1}(z-\beta_i)\prod^{L^*}_{j=1}(1-\beta^*_j z^{-1})\,\Xi^*_n(z) $ is given by
(\ref{KK*LL*shift}) with $ C \mapsto C^* $, $ \phi_n \mapsto \xi^*_n $ in the first row,
$ \phi_n \mapsto \phi^*_n $ in the second and third block rows and 
$ \xi_n \mapsto \xi^*_n $ in the last two block rows.
To find the leading coefficient $ \textsc{K}_n $ of the modified polynomials one can use 
the ratio formula
\begin{equation*}
   \textsc{K}^2_n = \frac{2[z^n]\Phi_n(z)}{[z^n] \Xi_n(z)} = \frac{2\Phi^*_n(0)}{\lim_{z\to \infty} \Xi^*_n(z)} .
\end{equation*}}
\end{minipage}
\end{sideways}
\vfill\eject

There is a recent parallel body of work \cite{AV_2003}, \cite{Be_2004},
\cite{HO_2006}, \cite{HO_2007} which treats a BOPS
defined with a measure $ d\mu(x,y) $ on a two-dimensional domain $ (x,y) \subset \R^2 $ 
and the rational modifications of this measure. This is more general in one sense 
than the problem we consider here and there are significant simplifications occurring 
in our system such as the existence of closed form expressions for the Christoffel-Darboux 
sums which are not present in theirs. As a consequence the resulting expressions have
rather different forms. In the main these works only consider the 
normalisation constants of the modified system and are completely formal in their 
approach. In addition there is no attempt to develop a theory of the 
semi-classical weight case and its connection to isomonodromy preserving systems.  

For the class of regular semi-classical weights on the unit circle, taking for the 
purposes of illustration the simple example 
\begin{equation}
  w(z) = \begin{cases} \prod_{j=0}^{M} (z-z_j)^{\rho_j} \\
                       \prod_{j=0}^{M} (1-z_jz^{-1})^{\rho^*_j}
         \end{cases} ,
\label{}
\end{equation}
it was first shown in \cite{FW_2004a} that this bi-orthogonal system
is monodromy preserving in the spectral variable $ z $ with respect to deformations 
of the singularities $ z_j $. 
This fact is the unit circle analogue of the result for orthogonal polynomial 
systems on the line with regular semi-classical weights, first shown by Magnus in 
1995 \cite{Ma_1995a}. Either of these systems constitute a realisation of a Garnier
system and one can construct a dictionary to translate one into the other,
and we refer the reader to \cite{Wi_2007b} for details.
In the regular semi-classical context the Christoffel-Geronimus-Uvarov transformations, under 
the specialisation of positioning the zeros or poles of the rational factors at 
the singular points of the weight, are the Schlesinger transformations of this
Garnier system. The fundamental Schlesinger transformations shift the exponents
$ \rho_j \mapsto \rho_j\pm 1 $, $ \rho^*_j \mapsto \rho^*_j\pm 1 $, 
so that we are in effect considering a family of
regular semi-classical weights with exponents displaced by integral values.

The theory of Schlesinger transformations for Garnier systems was
initiated by Schlesinger \cite{Sch_1912} and Garnier \cite{Ga_1912} where only the regular
singular case was treated, and the theory for irregular singular points was completed 
in the contributions by the Kyoto school 
\cite{JM_1980a},\cite{JM_1980b},\cite{JM_1980c},\cite{JM_1981a},\cite{JM_1981b}.
Schlesinger transformations are a special type of B\"acklund transformation, which
were systematically studied in the simplest non-trivial case of a Garnier system, i.e. the sixth
Painlev\'e system by Okamoto \cite{Ok_1987a}. 
All the B\"acklund transformations can be composed of fundamental transformations,
generated by the elementary reflections and automorphisms of the extended affine
Weyl groups $ D^{(1)}_4 $ for $ M=3 $ and $ B^{(1)}_{M+1} $ for $ M\geq 4 $ and
the permutation group $ S_{M+1} $ \cite{Ki_1990},\cite{IKSY_1991},\cite{Ts_2003a},\cite{Su_2005}.
 
In the final section
we present an alternative construction of the Schlesinger transformation theory to
the studies mentioned above, starting from the BOPS with a
regular semi-classical weight, which constitutes the most general classical solution to the
Garnier system. Using all of the structures of such a system from the approximation
theory point of view we can recover all of the known results, and some novel
results to be reported in \cite{Wi_2007b}, in a very efficient and transparent manner. 
Furthermore the present study can also be seen as the completion of the tasks begun in \cite{FW_2004a},
in that this work didn't cover the aspects of the integrable system relating to 
the Schlesinger transformations or the systems of bi-linear difference equations
governing the $\tau$-functions.
Our approach is to employ the simplest of classical arguments from approximation 
theory although we fully recognise that all of the results presented here could
be found in an elegant manner using Riemann-Hilbert methods. A natural formulation
of the bi-orthogonal system as a Riemann-Hilbert problem was given in 
\cite{FW_2004a} and it remains an interesting task to employ these techniques on the questions
posed here.

\section{General Structures of Bi-orthogonality}
\setcounter{equation}{0}

Let $ \T $ denote the unit circle with $ \zeta=e^{i\theta} $, $ \theta\in (-\pi,\pi] $. 
Consider an absolutely continuous weight $ w(z) $ and assume that its Fourier decomposition 
has meaning in a domain containing $ \T $,
\begin{equation}
  w(z) = \sum_{k=-\infty}^{\infty} w_{k}z^k, \quad
  w_{k} = \int_{\mathbb T} \frac{d\zeta}{2\pi i\zeta} w(\zeta)\zeta^{-k} .
\label{Fcoeff}
\end{equation}

We define the bi-orthogonal polynomials $ \{\phi_n(z),\bar{\phi}_n(z)\}^{\infty}_{n=0} $
with respect to the weight $ w(z) $ with support contained
within the unit circle by the orthogonality relation
\begin{equation}
  \int_{\T} \frac{d\zeta}{2\pi i\zeta} w(\zeta)\phi_m(\zeta)\bar{\phi}_n(\bar{\zeta})
   = \delta_{m,n} .
\label{onorm}
\end{equation}
Alternatively one can express this definition in terms of orthogonality with 
respect to the monomial basis
\begin{align}
  \int_{\T} \frac{d\zeta}{2\pi i\zeta} w(\zeta)\phi_n(\zeta)\bar{\zeta}^m
  & = \begin{cases}  0 & m<n \\ 1/\kappa_n & m=n \end{cases} ,
\label{orthog:a} \\
  \int_{\T} \frac{d\zeta}{2\pi i\zeta} w(\zeta)\zeta^m\bar{\phi}_n(\bar{\zeta})
  & = \begin{cases}  0 & m<n \\ 1/\kappa_n & m=n \end{cases} .
\label{orthog:b}
\end{align}
Notwithstanding the notation,
$ \bar{\phi}_n $ is not in general equal to the complex conjugate of $ \phi_n $.
We set
\begin{align}
  \phi_n(z)       &= \kappa_nz^n + \lambda_nz^{n-1} + \mu_nz^{n-2} +\ldots +\phi_n(0) ,
  \nonumber\\
  \bar{\phi}_n(z) &= \kappa_nz^n + \bar{\lambda}_nz^{n-1} + \bar{\mu}_nz^{n-2} +\ldots +\bar{\phi}_n(0) , 
  \nonumber
\end{align}
where again $ \bar{\lambda}_n $, $ \bar{\mu}_n $, $ \bar{\phi}_n(0) $ are not in general 
equal to the corresponding complex conjugate. We have chosen the leading coefficients
of both polynomials to be the same without any loss of generality as this preserves
the relations arising from the orthogonal polynomial theory with, of course, the 
distinction made above. The initial values of the polynomials are 
$ \phi_0(z) = \bar{\phi}_0(z) = \kappa_0 $ and the normalisation implies $ w_0\kappa_0^2 = 1 $.

Denote the $ U(n) $ average by $ I_n[w] $,
\begin{equation}
  I_n[w] = \Big \langle \prod_{l=1}^n w(z_l) \Big \rangle_{U(n)}
         = \det[ w_{j-k} ]_{j,k=0,\dots,n-1}. 
\label{Toep}
\end{equation}
By a result of Baxter \cite{Ba_1961} the existence of the bi-orthogonal system is 
guaranteed if and only if $ I_n[w] \neq 0 $ for $ n \in\Z_{\geq 0} $.
We define a sequence of $r$-coefficients by
\begin{equation}
  r_n = \frac{\phi_n(0)}{\kappa_n}, \quad 
  \bar{r}_n = \frac{\bar{\phi}_n(0)}{\kappa_n} ,
  \quad r_0 = \bar{r}_0 = 1,
\end{equation}
which differs slightly from the standard definition of the reflection or Verblunsky coefficients
$ \alpha_n $ in that $ \alpha_n = -\bar{r}_{n+1} $.
It is a well known result in the theory of Toeplitz determinants that
\begin{equation}
   \frac{I_{n+1}[w] I_{n-1}[w]}{(I_{n}[w])^2}
   = 1 - r_{n}\bar{r}_n .
\label{I0}
\end{equation}
The coefficients are related by many coupled equations, two of the simplest being
\begin{equation}
   \kappa_n^2 = \kappa_{n-1}^2 + \phi_n(0)\bar{\phi}_n(0), \quad
   \frac{\lambda_{n}}{\kappa_{n}}-\frac{\lambda_{n-1}}{\kappa_{n-1}} = r_{n}\bar{r}_{n-1}.
   \label{l}
\end{equation}

Introduce the reciprocal polynomial $ \phi^*_n(z) $ of the $ n$th degree polynomial
$ \bar{\phi}_n(z) $ by
\begin{equation}
  \phi^*_n(z) \coloneqq z^n\bar{\phi}_n(1/z) .
\label{recip}
\end{equation}
Fundamental to our study \cite{FW_2004a} is the matrix
\begin{equation}
   Y_n(z) \coloneqq 
   \begin{pmatrix}
          \phi_n(z)   &  \xi_n(z)/w(z) \cr
          \phi^*_n(z) & -\xi^*_n(z)/w(z) \cr
   \end{pmatrix} ,\quad n\geq 0,
\label{Ydefn}
\end{equation}
where the associated functions are defined, for $ n > 0 $,
\begin{align}
   \xi_n(z)
   &\coloneqq \int_{\T}\frac{d\zeta}{2\pi i\zeta}\frac{\zeta+z}{\zeta-z}w(\zeta)
                   \phi_n(\zeta) ,
   \label{eps:a} \\
   \xi^*_n(z)
   &\coloneqq -z^n\int_{\T}\frac{d\zeta}{2\pi i\zeta}\frac{\zeta+z}{\zeta-z}w(\zeta)
                   \bar{\phi}_n(\bar{\zeta}) ,
   \label{eps:b} \\
   &\coloneqq \frac{1}{\kappa_n}-\int_{\T}\frac{d\zeta}{2\pi i\zeta}\frac{\zeta+z}{\zeta-z}w(\zeta)
                   \phi^*_n(\zeta) .
   \label{eps:c}
\end{align}
A central object in the theory is the Carath\'eodory function, or generating
function of the Toeplitz elements,
\begin{equation}
   F(z) \coloneqq \int_{\T}\frac{d\zeta}{2\pi i\zeta}\frac{\zeta+z}{\zeta-z}w(\zeta) .
\label{Cfun}
\end{equation}
The initial values of the associated functions are defined as
$ \xi_0(z) = \kappa_0[w_0+F(z)] $ and $ \xi^*_0(z) = \kappa_0[w_0-F(z)] $. 

A consequence of the definitions (\ref{orthog:a}), (\ref{orthog:b}) are the following determinantal
representations for the polynomials,
\begin{equation}
  \phi_{n}(z) = \frac{\kappa_n}{I_{n}}
         \det \begin{pmatrix}
                     w_{0}      & \ldots & w_{-j}       & \ldots        & w_{-n} \cr
                     \vdots     & \vdots & \vdots & \vdots      & \vdots \cr
                     w_{n-1}    & \ldots & w_{n-j-1}    & \ldots        & w_{-1} \cr
                     1  & \ldots & z^j  & \ldots        & z^n  \cr
              \end{pmatrix} ,
  \label{DetRep:a}
\end{equation}
\begin{equation}
  \phi^*_{n}(z) = \frac{\kappa_n}{I_{n}}
         \det \begin{pmatrix}
                     w_{0}      & \ldots & w_{-n+1}     & z^n  \cr
                     \vdots     & \vdots & \vdots       & \vdots \cr
                     w_{n-j}    & \ldots & w_{-j+1}     & z^j  \cr
                     \vdots     & \vdots & \vdots       & \vdots \cr
                     w_{n}      & \ldots & w_{1}        & 1  \cr
              \end{pmatrix} .
  \label{DetRep:b}
\end{equation}
The associated functions have representations analogous to (\ref{IntRep:a},\ref{IntRep:b})
\begin{equation}
  \xi_{n}(z)
   = \frac{\kappa_n}{I_{n}}
         \det \begin{pmatrix}
                     w_{0}      & \ldots & w_{-j}       & \ldots        & w_{-n} \cr
                     \vdots     & \vdots & \vdots & \vdots      & \vdots \cr
                     w_{n-1}    & \ldots & w_{n-j-1}    & \ldots        & w_{-1} \cr
                     g_0(z)     & \ldots & g_j(z)       & \ldots        & g_n(z)  \cr
              \end{pmatrix} ,
  \label{DetRep:c}
\end{equation}
\begin{equation}
  \xi^*_{n}(z)
   = -\frac{\kappa_n}{I_{n}}
         \det \begin{pmatrix}
                     w_{0}      & \ldots & w_{-n+1}     & g_n(z)  \cr
                     \vdots     & \vdots & \vdots       & \vdots \cr
                     w_{n-j}    & \ldots & w_{-j+1}     & g_j(z)  \cr
                     \vdots     & \vdots & \vdots       & \vdots \cr
                     w_{n}      & \ldots & w_{1}        & g_0(z)  \cr
              \end{pmatrix} ,
  \label{DetRep:d}
\end{equation}
where
\begin{equation}
   g_{j}(z) \coloneqq 2z\int_{\T}\frac{d\zeta}{2\pi i\zeta}\frac{\zeta^j}{\zeta-z}w(\zeta) .
\end{equation}

A further consequence of the definitions is the difference system \cite{Ge_1961}
\begin{equation}
   Y_{n+1} \coloneqq K_n Y_{n}
   = \frac{1}{\kappa_n}
       \begin{pmatrix}
              \kappa_{n+1} z   & \phi_{n+1}(0) \cr
              \bar{\phi}_{n+1}(0) z & \kappa_{n+1} \cr
       \end{pmatrix} Y_{n} ,
\label{Yrecur:a}
\end{equation}
or its inverse
\begin{equation}
   Y_{n}
   = \frac{1}{\kappa_n}
       \begin{pmatrix}
              \frac{\DySt \kappa_{n+1}}{\DySt z}   & -\frac{\DySt \phi_{n+1}(0)}{\DySt z} \cr
              -\bar{\phi}_{n+1}(0) & \kappa_{n+1} \cr
       \end{pmatrix} Y_{n+1} .
\label{Yrecur:b}
\end{equation}
We also note that
\begin{equation}
   \det Y_{n}(z)
   = -\frac{2z^n}{w(z)} .
\label{Ydet}
\end{equation}

If one writes the coupled first order difference system (\ref{Yrecur:a}) as one decoupled
second order difference equation for the monic polynomial $ \varphi_n := \phi_n/\kappa_n $
one finds
\begin{equation}
  \varphi_{n+1}-\frac{r_{n+1}}{r_{n}}\varphi_{n} 
   = z\left[ \varphi_{n}-\frac{r_{n+1}}{r_{n}}(1-r_{n}\bar{r}_{n})\varphi_{n-1} \right] .
\end{equation}
These are the recurrence relations that serve to define one of the bi-orthogonal
Laurent polynomials \cite{JT_1982},\cite{JTN_1984},\cite{HvR_1986}
(the partner polynomial is found under the substitution $ r_n \mapsto \bar{r}_{n} $,
$ \varphi_n \mapsto \bar{\varphi}_n $).

\begin{theorem}[\cite{Ge_1961}]
The Casoratians of the solutions
$ \phi_n, \phi^*_n, \xi_n, \xi^*_{n} $ are
\begin{align}
  \phi_{n+1}(z)\xi_n(z) - \xi_{n+1}(z)\phi_n(z)
  & = 2\frac{\phi_{n+1}(0)}{\kappa_n}z^n ,
  \label{Cas:a} \\
  \phi^*_{n+1}(z)\xi^*_n(z) - \xi^*_{n+1}(z)\phi^*_n(z)
  & = 2\frac{\bar{\phi}_{n+1}(0)}{\kappa_n}z^{n+1} ,
  \label{Cas:b} \\
  \phi_{n}(z)\xi^*_n(z) + \xi_{n}(z)\phi^*_n(z)
  & = 2z^n ,
  \label{Cas:c}
\end{align}
valid for $ n \geq 0 $.
\end{theorem}

The bi-orthogonal system satisfies the spectral differential system 
\cite{FW_2004a} (with the abbreviation $ d/dz \coloneqq {}' $)
\begin{multline}
   \frac{d}{dz}Y_{n} \coloneqq A_n Y_{n}
   \\
   = \frac{1}{W(z)}
       \begin{pmatrix}
              -\left[ \Omega_n(z)+V(z)
                     -\dfrac{\kappa_{n+1}}{\kappa_n}z\Theta_n(z)
               \right]
            & \dfrac{\phi_{n+1}(0)}{\kappa_n}\Theta_n(z)
            \cr
              -\dfrac{\bar{\phi}_{n+1}(0)}{\kappa_n}z\Theta^*_n(z)
            &  \Omega^*_n(z)-V(z)
                     -\dfrac{\kappa_{n+1}}{\kappa_n}\Theta^*_n(z)
            \cr
       \end{pmatrix} Y_{n} ,
\label{YzDer}
\end{multline}
for $ n \geq 0 $ under the sufficient conditions of the existence of
the moments
\begin{equation*} 
  \int_{\T} \frac{d\zeta}{2\pi i\zeta} w(\zeta) \zeta^k
  \frac{[\log w(\zeta)]'-[\log w(z)]'}{\zeta-z}, \quad k \in \Z.
\end{equation*} 
The particular parameterisation of $ A_n $ into the spectral coefficients
$ W(z), V(z) $ and $ \Omega_n(z),\Omega^*_n(z), \Theta_n(z), \Theta^*_n(z) $
will serve a purpose when we specialise to the regular semi-classical 
weights. We see that
\begin{equation}
   {\rm Tr}A_{n} = \frac{n}{z}-\frac{w'}{w} .
\label{traceA}
\end{equation}
The compatibility relation between the recurrence relations (\ref{Yrecur:a}) and 
the spectral derivative (\ref{YzDer}) is
\begin{equation}
   K'_{n} = A_{n+1}K_{n}-K_{n}A_{n} ,
\label{Ycompat}
\end{equation}
and this is equivalent to recurrences in $ n $ for the spectral coefficients,
which can be found in \cite{FW_2004a} but are not required here.

We will require the leading order terms in expansions of
$ \phi_n(z), \phi^*_n(z), \xi_n(z), \xi^*_{n}(z) $ both inside and
outside the unit circle.
\begin{corollary}
The bi-orthogonal polynomials $ \phi_n(z), \phi^*_n(z) $ have the following
expansions for $ |z| < 1 $
\begin{multline}
   \phi_n(z) =
      \phi_n(0)
       + \dfrac{1}{\kappa_{n-1}}
         (\kappa_n\phi_{n-1}(0)+\phi_n(0)\bar{\lambda}_{n-1})z \\
       + \left[\dfrac{\kappa_n}{\kappa_{n-1}\kappa_{n-2}}
               (\kappa_{n-1}\phi_{n-2}(0)+\phi_{n-1}(0)\bar{\lambda}_{n-2})
               +\dfrac{\phi_{n}(0)\bar{\mu}_{n-1}}{\kappa_{n-1}} \right]z^2
       + {\rm O}(z^3) ,
\label{phiexp:a}
\end{multline}
\begin{equation}
   \phi^*_n(z) =
      \kappa_n + \bar{\lambda}_n z + \bar{\mu}_n z^2 + {\rm O}(z^{3}) ,
\label{phiexp:b}
\end{equation}
whilst the associated functions have the expansions
\begin{multline}
   \dfrac{\kappa_n}{2}\xi_n(z) =
      z^n - \dfrac{\bar{\lambda}_{n+1}}{\kappa_{n+1}}z^{n+1} \\
   + \left[\dfrac{\bar{\lambda}_{n+1}\bar{\lambda}_{n+2}}{\kappa_{n+1}\kappa_{n+2}}
           -\dfrac{\bar{\mu}_{n+2}}{\kappa_{n+2}}\right]z^{n+2} \\
   + \left[\dfrac{\bar{\lambda}_{n+1}\bar{\mu}_{n+3}}{\kappa_{n+1}\kappa_{n+3}}
           +\dfrac{\bar{\lambda}_{n+3}\bar{\mu}_{n+2}}{\kappa_{n+2}\kappa_{n+3}}
           -\dfrac{\bar{n}_{n+3}}{\kappa_{n+3}}
           -\dfrac{\bar{\lambda}_{n+1}\bar{\lambda}_{n+2}\bar{\lambda}_{n+3}}{\kappa_{n+1}\kappa_{n+2}\kappa_{n+3}}
     \right]z^{n+3} 
       + {\rm O}(z^{n+4}) ,
\label{epsexp:a}
\end{multline}
\begin{multline}
   \dfrac{\kappa_n}{2}\xi^*_n(z) =
         \dfrac{\bar{\phi}_{n+1}(0)}{\kappa_{n+1}}z^{n+1}
       + \left( \dfrac{\bar{\phi}_{n+2}(0)}{\kappa_{n+2}}
               -\dfrac{\bar{\phi}_{n+1}(0)\bar{\lambda}_{n+2}}{\kappa_{n+1}\kappa_{n+2}}
         \right)z^{n+2} \\
       + \left[ \dfrac{\bar{\phi}_{n+3}(0)}{\kappa_{n+3}}
               -\dfrac{\bar{\phi}_{n+1}(0)\bar{\mu}_{n+3}}{\kappa_{n+1}\kappa_{n+3}}
               +\dfrac{\bar{\lambda}_{n+3}(0)}{\kappa_{n+3}}\left(
                   \dfrac{\bar{\phi}_{n+1}(0)\bar{\lambda}_{n+2}}{\kappa_{n+1}\kappa_{n+2}}
                  -\dfrac{\bar{\phi}_{n+2}(0)}{\kappa_{n+2}}\right)
         \right]z^{n+3}
       + {\rm O}(z^{n+4}) .
\label{epsexp:b}
\end{multline}
The large argument expansions $ |z| > 1 $ are 
\begin{equation}
   \phi_n(z) =
      \kappa_n z^n + \lambda_n z^{n-1} + \mu_n z^{n-2} + {\rm O}(z^{n-3}) ,
\label{phiexp:c}
\end{equation}
\begin{multline}
   \phi^*_n(z) =
      \bar{\phi}_n(0) z^n
       + \dfrac{1}{\kappa_{n-1}}
         (\kappa_n\bar{\phi}_{n-1}(0)+\bar{\phi}_n(0)\lambda_{n-1})z^{n-1} \\
       + \left[\dfrac{\kappa_n}{\kappa_{n-1}\kappa_{n-2}}
               (\kappa_{n-1}\bar{\phi}_{n-2}(0)+\bar{\phi}_{n-1}(0)\lambda_{n-2})
               +\dfrac{\bar{\phi}_{n}(0)\mu_{n-1}}{\kappa_{n-1}} \right]z^{n-2} 
       + {\rm O}(z^{n-3}) ,
\label{phiexp:d}
\end{multline}
whilst the associated functions have the expansions
\begin{equation}
   \dfrac{\kappa_n}{2} \xi_n(z) =
         \dfrac{\phi_{n+1}(0)}{\kappa_{n+1}}z^{-1}
       + \left(\dfrac{\kappa^2_n}{\kappa^2_{n+1}}
               \dfrac{\phi_{n+2}(0)}{\kappa_{n+2}}
               -\dfrac{\phi_{n+1}(0)}{\kappa_{n+1}}
                \dfrac{\lambda_{n+1}}{\kappa_{n+1}}
         \right)z^{-2}
       + {\rm O}(z^{-3}) ,
\label{epsexp:c}
\end{equation}
\begin{equation}
   \dfrac{\kappa_n}{2} \xi^*_n(z) =
       1 - \dfrac{\lambda_{n+1}}{\kappa_{n+1}}z^{-1} + \left(
         \dfrac{\lambda_{n+2}\lambda_{n+1}}{\kappa_{n+2}\kappa_{n+1}}
        -\dfrac{\mu_{n+2}}{\kappa_{n+2}}
         \right)z^{-2} + {\rm O}(z^{-3}) .
\label{epsexp:d}
\end{equation}
\end{corollary}
From the Riemann-Hilbert perspective \cite{FW_2004a} it is apparent that our
bi-orthogonal system demands the existence of two domains, one containing $ z=0 $
and the other containing $ z=\infty $, where the above expansions are valid in
order to specify the system.

\section{Christoffel-Geronimus-Uvarov Transformations}
\setcounter{equation}{0}
Starting with the determinantal representations of the polynomials and associated functions,
given in (\ref{DetRep:a}), (\ref{DetRep:b}), (\ref{DetRep:c}) and 
(\ref{DetRep:d}),
it is easy to recast these as multiple-integrals using the Heine formula. We observe
that the following integrals occur with the weight modified by a rational factor
\begin{align}
  I_{n}[w(\zeta)(\zeta-z)] 
  & = (-)^nI_{n}[w(\zeta)]\frac{\phi_{n}(z)}{\kappa_n}, \quad n\geq 0,
  \label{IntRep:a} \\
  I_{n}[w(\zeta)(1-z\zeta^{-1})] 
  & = I_{n}[w(\zeta)]\frac{\phi^*_{n}(z)}{\kappa_n}, \quad n\geq 0,
  \label{IntRep:b} \\
  I_{n}[w(\zeta)\frac{1}{(1-z\zeta^{-1})}] 
  & = I_{n}[w(\zeta)]\frac{\kappa_{n-1}\xi_{n-1}(z)}{2z^{n-1}}, \quad n>0, |z|\neq 1,
  \label{IntRep:c} \\
  I_{n}[w(\zeta)\frac{1}{(\zeta-z)}] 
  & = (-)^{n}I_{n}[w(\zeta)]\frac{\kappa_{n-1}\xi^*_{n-1}(z)}{2z^{n}}, \quad n>0, |z|\neq 1.
  \label{IntRep:d}
\end{align}
These are the analogues of the Christoffel formula or the simplest cases of the Uvarov 
formulae for orthogonal polynomials \cite{Uv_1959}, \cite{Uv_1969}, 
and are specialisations of the Ismail and Ruedemann formulae for OPS on the unit circle
\cite{IR_1992} and those of Ruedemann \cite{Ru_1994} for BOPS on the unit circle. 
We seek to generalise the above expressions and derive formulae for the BOPS corresponding to an
arbitrarily rationally modified weight (\ref{C-Uxfm:b}) in terms of the original BOPS. 

However Ismail and Ruedemann's study is incomplete for our purposes in a number of ways -
only the modification of a real, positive weight on the unit circle by a quotient of a real, 
non-negative polynomial and a real, positive polynomial is considered. This doesn't cover 
the case of complex weights, or equivalently non-hermitian Toeplitz determinants, and 
consequently the bi-orthogonal systems associated with this type of weight. In addition this type of
transformation is a composite of fundamental transformations and it is not possible to 
dis-entangle these. In a subsequent work \cite{Ru_1994} the extension to the bi-orthogonal 
system was made, but still only treating the case of a composite transformation.
In both of these studies formulae were given only 
for the modified polynomials $ \phi_n, \phi^*_n $ whereas we require expressions for the complete 
system, that is to say for $ \xi_n, \xi^*_n $ as well. Another technical caveat
applying to the results of \cite{IR_1992} is that the degree of the polynomial
$ n $ cannot be less than the degree of the denominator polynomial $ L+L^* $.
The studies by Godoy and Marcell{\'a}n \cite{GM_1987},\cite{GM_1991},\cite{GM_1993}
are limited in similar ways, and have the additional complication that expressions
for the modified orthogonal polynomial system are given in terms of determinants 
with elements containing reproducing polynomials. 
Formulae for the polynomials and recurrence coefficients of the bi-orthogonal Laurent system
under a modification of the weight by multiplication with a single factor, and division
by a simple factor plus addition of a mass point were derived by Zhedanov \cite{Zh_1998}.   
While we do not consider the problem of additional mass points at all here we will extend
his results to arbitrary rational modification. 

We will carry out our task in seven steps, each one exhibiting our methods of deduction in 
their simplest setting, rather than give the proof for the general case. 
In the notation of the most general modification
(\ref{C-Uxfm:a}), (\ref{C-Uxfm:b}) we will derive formulae for the cases -
$ K\neq 0 $, $ n>L\neq 0 $, $ n<L\neq 0 $, $ K,L\neq 0 $, $ K^*\neq 0 $, $ L^*\neq 0 $ and
$ K^*,L^*\neq 0 $. The most general case has been given in (\ref{KK*LL*shift}).
In each step we will denote the modified bi-orthogonal polynomials and associated functions by
the same symbols $ \Phi_n(z) $, $ \Phi^*_n(z) $, $  \Xi_n(z) $, $  \Xi^*_n(z) $
and it should be appreciated by the reader that the formulae only apply in their specific context. 
We also denote the modified leading and $r$-coefficients by $ \textsc{K}_n, \textsc{R}_n, \bar{\textsc{R}}_n$;
$ C, C^* $ are arbitrary  constants of proportionality independent of the spectral variable $ z $;
and lastly we denote the space of polynomials in $ z $ with degree at most $ n $ by $ \Pi_{n}(z) $.
The same normalisation condition applies for the modified system, i.e. (\ref{onorm}) holds.
We will assume hereafter that the following {\it generic conditions} apply -
\begin{itemize}
 \item
 that the zeros and poles are not pairwise coincident
 $ \alpha_j\neq \alpha_k, \beta_j\neq \beta_k, \alpha^*_j\neq \alpha^*_k, \beta^*_j\neq \beta^*_k, j\neq k $,
 in the proofs for convenience although confluent formulae can be derived from our results,
 \item
 assume the existence of the original BOPS and of the modified BOPS, 
 i.e. that $ r_n\bar{r}_n \neq 0,1 $ and $ \textsc{R}_n\bar{\textsc{R}}_n \neq 0,1 $ for all $ n \in \N $.
\end{itemize}

\begin{proposition}\label{Kmod}
For the polynomial modification of the weight
\begin{equation}
   w\left[\begin{array}{c} \ScSt K \\ \ScSt \cdot \end{array};z \right] = \prod^{K}_{j=1}(z-\alpha_j) w(z) ,
\end{equation}
the corresponding bi-orthogonal polynomials are given by
\begin{equation}
  \prod^{K}_{j=1}(z-\alpha_j)\,\Phi_n(z) =
  (-1)^{K}b_{n,n+K}
   \frac{\det \begin{pmatrix}
                     \phi_{n}(z)	& \ldots & \phi_{n+K}(z)	\cr
                     \phi_{n}(\alpha_1)	& \ldots & \phi_{n+K}(\alpha_1)	\cr
                     \vdots		& \vdots & \vdots		\cr
                     \phi_{n}(\alpha_K)	& \ldots & \phi_{n+K}(\alpha_K)	\cr
              \end{pmatrix}}
        {\det \begin{pmatrix}
                     \phi_{n}(\alpha_1)	& \ldots & \phi_{n+K-1}(\alpha_1)	\cr
                     \vdots		& \vdots & \vdots		\cr
                     \phi_{n}(\alpha_K)	& \ldots & \phi_{n+K-1}(\alpha_K)	\cr
              \end{pmatrix}}, \quad n\geq 0,
\label{Kmod:a}
\end{equation}
\begin{equation}
  \prod^{K}_{j=1}(z-\alpha_j)\,\Phi^*_n(z) =
  (-1)^{K}\hat{b}_{n,n+K}
   \frac{\det \begin{pmatrix}
                     \phi^*_{n}(z)	& \ldots & \phi^*_{n+K}(z)	\cr
                     \phi^*_{n}(\alpha_1)	& \ldots & \phi^*_{n+K}(\alpha_1)	\cr
                     \vdots		& \vdots & \vdots		\cr
                     \phi^*_{n}(\alpha_K)	& \ldots & \phi^*_{n+K}(\alpha_K)	\cr
              \end{pmatrix}}
        {\det \begin{pmatrix}
                     \phi^*_{n}(\alpha_1)	& \ldots & \phi^*_{n+K-1}(\alpha_1)	\cr
                     \vdots		& \vdots & \vdots		\cr
                     \phi^*_{n}(\alpha_K)	& \ldots & \phi^*_{n+K-1}(\alpha_K)	\cr
              \end{pmatrix}}, \quad n\geq 0,
\label{Kmod:b}
\end{equation}
and the associated functions
\begin{equation}
   \Xi_n(z) =
  (-1)^{K}b_{n,n+K}
   \frac{\det \begin{pmatrix}
                     \xi_{n}(z)	& \ldots & \xi_{n+K}(z)	\cr
                     \phi_{n}(\alpha_1)	& \ldots & \phi_{n+K}(\alpha_1)	\cr
                     \vdots		& \vdots & \vdots		\cr
                     \phi_{n}(\alpha_K)	& \ldots & \phi_{n+K}(\alpha_K)	\cr
              \end{pmatrix}}
        {\det \begin{pmatrix}
                     \phi_{n}(\alpha_1)	& \ldots & \phi_{n+K-1}(\alpha_1)	\cr
                     \vdots		& \vdots & \vdots		\cr
                     \phi_{n}(\alpha_K)	& \ldots & \phi_{n+K-1}(\alpha_K)	\cr
              \end{pmatrix}}, \quad n\geq 0,
\label{Kmod:c}
\end{equation}
\begin{equation}
   \Xi^*_n(z) =
  (-1)^{K}\hat{b}_{n,n+K}
   \frac{\det \begin{pmatrix}
                     \xi^*_{n}(z)	& \ldots & \xi^*_{n+K}(z)	\cr
                     \phi^*_{n}(\alpha_1)	& \ldots & \phi^*_{n+K}(\alpha_1)	\cr
                     \vdots		& \vdots & \vdots		\cr
                     \phi^*_{n}(\alpha_K)	& \ldots & \phi^*_{n+K}(\alpha_K)	\cr
              \end{pmatrix}}
        {\det \begin{pmatrix}
                     \phi^*_{n}(\alpha_1)	& \ldots & \phi^*_{n+K-1}(\alpha_1)	\cr
                     \vdots		& \vdots & \vdots		\cr
                     \phi^*_{n}(\alpha_K)	& \ldots & \phi^*_{n+K-1}(\alpha_K)	\cr
              \end{pmatrix}}, \quad n\geq 0 .
\label{Kmod:d}
\end{equation}
Here $  b_{n,n+K} = \textsc{K}_n/\kappa_{n+K} $, 
$ \hat{b}_{n,n+K} = \bar{\Phi}_n(0)/\bar{\phi}_{n+K}(0) $ 
and the leading and $r$-coefficients are given by
\begin{equation}
   \prod^{K}_{i=1}\alpha_{i}\textsc{K}_n^2 
   = (-)^K
   \frac{\det \begin{pmatrix}
                     \kappa_{n}	& \ldots & \kappa_{n+K}	\cr
                     \phi^*_{n}(\alpha_1)	& \ldots & \phi^*_{n+K}(\alpha_1)	\cr
                     \vdots		& \vdots & \vdots		\cr
                     \phi^*_{n}(\alpha_K)	& \ldots & \phi^*_{n+K}(\alpha_K)	\cr
              \end{pmatrix}}
        {\det \begin{pmatrix}
                     1/\kappa_{n}	& \ldots & 1/\kappa_{n+K}	\cr
                     \phi^*_{n}(\alpha_1)	& \ldots & \phi^*_{n+K}(\alpha_1)	\cr
                     \vdots		& \vdots & \vdots		\cr
                     \phi^*_{n}(\alpha_K)	& \ldots & \phi^*_{n+K}(\alpha_K)	\cr
              \end{pmatrix}} ,
\label{Kmod:e}
\end{equation} 
and
\begin{equation}
  \textsc{R}_n = \frac{1}{\kappa_{n+K}\prod^{K}_{i=1}\alpha_{i}}
   \frac{\det \begin{pmatrix}
                     \phi_{n}(0)	& \ldots & \phi_{n+K}(0)	\cr
                     \phi_{n}(\alpha_1)	& \ldots & \phi_{n+K}(\alpha_1)	\cr
                     \vdots		& \vdots & \vdots		\cr
                     \phi_{n}(\alpha_K)	& \ldots & \phi_{n+K}(\alpha_K)	\cr
              \end{pmatrix}}
        {\det \begin{pmatrix}
                     \phi_{n}(\alpha_1)	& \ldots & \phi_{n+K-1}(\alpha_1)	\cr
                     \vdots		& \vdots & \vdots		\cr
                     \phi_{n}(\alpha_K)	& \ldots & \phi_{n+K-1}(\alpha_K)	\cr
              \end{pmatrix}} ,
\label{Kmod:f}
\end{equation}
and
\begin{equation}
  \bar{\textsc{R}}_n = \phi_{n+K}(0)\prod^{K}_{i=1}\alpha_{i}
   \frac{\det \begin{pmatrix}
                     \phi^*_{n}(\alpha_1)	& \ldots & \phi^*_{n+K-1}(\alpha_1)	\cr
                     \vdots		& \vdots & \vdots		\cr
                     \phi^*_{n}(\alpha_K)	& \ldots & \phi^*_{n+K-1}(\alpha_K)	\cr
              \end{pmatrix}}
        {\det \begin{pmatrix}
                     \kappa_{n} 	& \ldots & \kappa_{n+K}	\cr
                     \phi^*_{n}(\alpha_1)	& \ldots & \phi^*_{n+K}(\alpha_1)	\cr
                     \vdots		& \vdots & \vdots		\cr
                     \phi^*_{n}(\alpha_K)	& \ldots & \phi^*_{n+K}(\alpha_K)	\cr
              \end{pmatrix}} .
\label{Kmod:g}
\end{equation}
\end{proposition}
\begin{proof}
Taking (\ref{Kmod:a}) first we note that $ \prod^K_{j=1}(z-\alpha_j)\,\Phi_n(z) \in \Pi_{n+K}(z) $ 
so that 
\begin{equation}
  \prod^K_{j=1}(z-\alpha_j)\,\Phi_n(z) = \sum^{n+K}_{j=0}b_{n,j}\phi_j(z) ,     
\end{equation}
with 
\begin{equation}
  b_{n,j} = \int\frac{d\zeta}{2\pi i\zeta}w(\zeta)\prod^K_{j=1}(\zeta-\alpha_j)\,\Phi_n(\zeta)\bar{\phi}_j(\bar{\zeta}) . 
\end{equation}
However by orthogonality (\ref{orthog:b}) $ b_{n,j} $ vanishes for $ j<n $ 
and in fact
\begin{equation}
  \prod^K_{j=1}(z-\alpha_j)\,\Phi_n(z) = \sum^{n+K}_{j=n}b_{n,j}\phi_j(z) .     
\end{equation}
To determine the non-vanishing coefficients we observe that there are $ K $ 
vanishing conditions
\begin{equation}
  \sum^{n+K}_{j=n}b_{n,j}\phi_j(\alpha_l) = 0, \quad l=1,\ldots, K ,     
\end{equation}
which enable the $ \{b_{n,j}\}^{n+K-1}_{j=n} $ to be solved in terms of $ b_{n,n+K} $.
This establishes (\ref{Kmod:a}). Result (\ref{Kmod:c}) follows from (\ref{Kmod:a}) 
upon using the definition (\ref{eps:a}).

To derive (\ref{Kmod:b}) we start with
\begin{equation}
  \prod^K_{j=1}(z-\alpha_j)\,\Phi^*_n(z) = \sum^{n+K}_{j=0}\hat{b}_{n,j}\phi^*_j(z) .     
\end{equation}
For $ m=0,\ldots,n-1 $ we see that
\begin{align}
   0 & = \int\frac{d\zeta}{2\pi i\zeta}w\left[\begin{array}{c} \ScSt K \\ \ScSt \cdot \end{array};\zeta \right]
                    \zeta^m\bar{\Phi}_n(\bar{\zeta}) ,
     \\
     & = \sum^{n+K}_{j=0}\hat{b}_{n,j}
         \int\frac{d\zeta}{2\pi i\zeta}w(\zeta)\zeta^{m-n+j}\bar{\phi}_j(\bar{\zeta}) .
\end{align}
However the integral given above vanishes for $ j\geq n-m $, and as
\begin{equation}
  \int\frac{d\zeta}{2\pi i\zeta}w(\zeta)\bar{\zeta}\bar{\phi}_j(\bar{\zeta}), \ldots,
  \int\frac{d\zeta}{2\pi i\zeta}w(\zeta)\bar{\zeta}^{n}\bar{\phi}_j(\bar{\zeta}) ,
\end{equation}
are non-zero we conclude that $ \hat{b}_{n,j}=0 $ for $ j=0,\ldots,n-1 $ and 
thus
\begin{equation}
  \prod^K_{j=1}(z-\alpha_j)\,\Phi^*_n(z) = \sum^{n+K}_{j=n}\hat{b}_{n,j}\phi^*_j(z) .     
\end{equation}
Again we have $ K $ relations 
\begin{equation}
  \sum^{n+K}_{j=n}\hat{b}_{n,j}\phi^*_j(\alpha_l) = 0, \quad l=1,\ldots, K ,     
\end{equation}
allowing us to solve for $ \{\hat{b}_{n,j}\}^{n+K-1}_{j=n} $ in terms of
$ \hat{b}_{n,n+K} $. This yields (\ref{Kmod:b}). From this result we can derive
(\ref{Kmod:d}) using the first definition (\ref{eps:b}) and (\ref{recip}).
This gives us the formula
\begin{equation}
   \Xi^*_n(z) = -(-1)^{K}\hat{b}_{n,n+K}z^n
   \frac{\det \begin{pmatrix}
                     \ldots & 
   \int\frac{\DySt d\zeta}{\DySt2\pi i\zeta}w(\zeta)\frac{\DySt\zeta+z}{\DySt\zeta-z}\zeta^{j}\bar{\phi}_{n+j}(\bar{\zeta})
                                                	& \ldots	\cr
                     \ldots & \phi^*_{n+j}(\alpha_1)	& \ldots	\cr
                            & \vdots									& \ldots	\cr
                     \ldots & \phi^*_{n+j}(\alpha_K)	& \ldots	\cr
              \end{pmatrix}}
        {\det (\phi^*_{n+j}(\alpha_k))} ,
\end{equation} 
but we rewrite the integral appearing here as
\begin{multline}
 \int\frac{d\zeta}{2\pi i\zeta}w(\zeta)\frac{\zeta+z}{\zeta-z}\zeta^{j}\bar{\phi}_{n+j}(\bar{\zeta})
 = \int\frac{d\zeta}{2\pi i\zeta}w(\zeta)(\zeta+z)\frac{\zeta^j-z^j}{\zeta-z}\bar{\phi}_{n+j}(\bar{\zeta})
  \\
  + z^j\int\frac{d\zeta}{2\pi i\zeta}w(\zeta)\frac{\zeta+z}{\zeta-z}\bar{\phi}_{n+j}(\bar{\zeta}) .     
\end{multline}
The first term vanishes because $ (\zeta+z)\frac{\DySt\zeta^j-z^j}{\DySt\zeta-z}\in \Pi_{j}(\zeta) $
and the second term evaluates to $ -z^{-n}\xi^*_{n+j}(z) $.

The relation for $ b_{n,n+K} $ preceding (\ref{Kmod:e}) follows from examining
the leading coefficient of (\ref{Kmod:a}) and comparison with (\ref{phiexp:c}).
Similarly the relation for $ \hat{b}_{n,n+K} $ follows from the leading coefficient
of (\ref{Kmod:b}) and using (\ref{phiexp:d}). The formula (\ref{Kmod:e}) is found
from the leading coefficient of the $ z\to\infty $ expansion in (\ref{epsexp:d})
and the corresponding coefficient of the $ z\to 0 $ expansion in (\ref{phiexp:b}).
The relation (\ref{Kmod:f}) is derived from the ratio of the leading and trailing
coefficients of (\ref{Kmod:a}) and the relation (\ref{Kmod:g}) from the 
leading coefficient of the $ z\to\infty $ expansion in (\ref{epsexp:d}) and 
the result for $ \hat{b}_{n,n+K} $. 
\end{proof}

\begin{proposition}\label{Lmod}
For the reciprocal polynomial modification of the weight
\begin{equation}
   w\left[\begin{array}{c} \ScSt \cdot \\ \ScSt L \end{array};z \right]
    = \frac{1}{\prod^{L}_{j=1}(z-\beta_j)} w(z) ,
\end{equation}
the corresponding bi-orthogonal polynomials are given by
\begin{equation}
  \Phi_n(z) = \frac{(-1)^L}{\kappa^L_n}
   \frac{\det \begin{pmatrix}
                     z^L\phi_{n-L}(z)	& \ldots & \phi_{n}(z)	\cr
                     \beta_1^{-n+L}\xi_{n-L}(\beta_1)	& \ldots & \beta_1^{-n}\xi_{n}(\beta_1)	\cr
                     \vdots		& \vdots & \vdots		\cr
                     \beta_L^{-n+L}\xi_{n-L}(\beta_L)	& \ldots & \beta_L^{-n}\xi_{n}(\beta_L)	\cr
              \end{pmatrix}}
        {\det \begin{pmatrix}
                     1/\kappa_{n-L}	& \ldots & 1/\kappa_{n}	\cr
                     \beta_1^{-n+L-1}\xi_{n-L}(\beta_1)	& \ldots & \beta_1^{-n-1}\xi_{n}(\beta_1)	\cr
                     \vdots		& \vdots & \vdots		\cr
                     \beta_L^{-n+L-1}\xi_{n-L}(\beta_L)	& \ldots & \beta_L^{-n-1}\xi_{n}(\beta_L)	\cr
              \end{pmatrix}}, \quad n\geq L,
\label{Lmod:a}
\end{equation}
\begin{equation}
  \Phi^*_n(z) = (-1)^{L}b'_{n,n}
   \frac{\det \begin{pmatrix}
                     z^L\phi^*_{n-L}(z)	& \ldots & \phi^*_{n}(z)	\cr
                     \beta_1^{-n+L}\xi^*_{n-L}(\beta_1)	& \ldots & \beta_1^{-n}\xi^*_{n}(\beta_1)	\cr
                     \vdots		& \vdots & \vdots		\cr
                     \beta_L^{-n+L}\xi^*_{n-L}(\beta_L)	& \ldots & \beta_L^{-n}\xi^*_{n}(\beta_L)	\cr
              \end{pmatrix}}
        {\det \begin{pmatrix}
                     \beta_1^{-n+L}\xi^*_{n-L}(\beta_1)	& \ldots & \beta_1^{-n+1}\xi^*_{n-1}(\beta_1)	\cr
                     \vdots		& \vdots & \vdots		\cr
                     \beta_L^{-n+L}\xi^*_{n-L}(\beta_L)	& \ldots & \beta_L^{-n+1}\xi^*_{n-1}(\beta_L)	\cr
              \end{pmatrix}}, \quad n\geq L,
\label{Lmod:b}
\end{equation}
and the associated functions
\begin{equation}
  \prod^{L}_{j=1}(z-\beta_j)\,\Xi_n(z) = C
   \frac{\det \begin{pmatrix}
                     z^L\xi_{n-L}(z)	& \ldots & \xi_{n}(z)	\cr
                     \beta_1^{-n+L}\xi_{n-L}(\beta_1)	& \ldots & \beta_1^{-n}\xi_{n}(\beta_1)	\cr
                     \vdots		& \vdots & \vdots		\cr
                     \beta_L^{-n+L}\xi_{n-L}(\beta_L)	& \ldots & \beta_L^{-n}\xi_{n}(\beta_L)	\cr
              \end{pmatrix}}
        {\det \begin{pmatrix}
                     \beta_1^{-n+L}\xi_{n-L}(\beta_1)	& \ldots & \beta_1^{-n+1}\xi_{n-1}(\beta_1)	\cr
                     \vdots		& \vdots & \vdots		\cr
                     \beta_L^{-n+L}\xi_{n-L}(\beta_L)	& \ldots & \beta_L^{-n+1}\xi_{n-1}(\beta_L)	\cr
              \end{pmatrix}}, \quad n\geq L,
\label{Lmod:c}
\end{equation}
\begin{equation}
  \prod^{L}_{j=1}(z-\beta_j)\,\Xi^*_n(z) = (-1)^{L}b'_{n,n}z^n
   \frac{\det \begin{pmatrix}
                     z^{-n+L}\xi^*_{n-L}(z)	& \ldots & z^{-n}\xi^*_{n}(z)	\cr
                     \beta_1^{-n+L}\xi^*_{n-L}(\beta_1)	& \ldots & \beta_1^{-n}\xi^*_{n}(\beta_1)	\cr
                     \vdots		& \vdots & \vdots		\cr
                     \beta_L^{-n+L}\xi^*_{n-L}(\beta_L)	& \ldots & \beta_L^{-n}\xi^*_{n}(\beta_L)	\cr
              \end{pmatrix}}
        {\det \begin{pmatrix}
                     \beta_1^{-n+L}\xi^*_{n-L}(\beta_1)	& \ldots & \beta_1^{-n+1}\xi^*_{n-1}(\beta_1)	\cr
                     \vdots		& \vdots & \vdots		\cr
                     \beta_L^{-n+L}\xi^*_{n-L}(\beta_L)	& \ldots & \beta_L^{-n+1}\xi^*_{n-1}(\beta_L)	\cr
              \end{pmatrix}}, \quad n\geq L.
\label{Lmod:d}
\end{equation}
Here $ b'_{n,n} = \textsc{K}_n/\kappa_n $ and the leading and $r$-coefficients 
are given by
\begin{equation}
  \textsc{K}_n^2 = (-)^L
   \frac{\det \begin{pmatrix}
                     \kappa_{n-L}	& \ldots & \kappa_{n}	\cr
                     \beta_1^{-n+L}\xi_{n-L}(\beta_1)	& \ldots & \beta_1^{-n}\xi_{n}(\beta_1)	\cr
                     \vdots		& \vdots & \vdots		\cr
                     \beta_L^{-n+L}\xi_{n-L}(\beta_L)	& \ldots & \beta_L^{-n}\xi_{n}(\beta_L)	\cr
              \end{pmatrix}}
        {\det \begin{pmatrix}
                     1/\kappa_{n-L}	& \ldots & 1/\kappa_{n}	\cr
                     \beta_1^{-n+L-1}\xi_{n-L}(\beta_1)	& \ldots & \beta_1^{-n-1}\xi_{n}(\beta_1)	\cr
                     \vdots		& \vdots & \vdots		\cr
                     \beta_L^{-n+L-1}\xi_{n-L}(\beta_L)	& \ldots & \beta_L^{-n-1}\xi_{n}(\beta_L)	\cr
              \end{pmatrix}} ,
\end{equation} 
and
\begin{equation}
  \textsc{R}_n = (-)^L\phi_n(0)
   \frac{\det \begin{pmatrix}
                     \beta_1^{-n+L}\xi_{n-L}(\beta_1)	& \ldots & \beta_1^{-n+1}\xi_{n-1}(\beta_1)	\cr
                     \vdots		& \vdots & \vdots		\cr
                     \beta_L^{-n+L}\xi_{n-L}(\beta_L)	& \ldots & \beta_L^{-n+1}\xi_{n-1}(\beta_L)	\cr
              \end{pmatrix}}
        {\det \begin{pmatrix}
                     \kappa_{n-L}	& \ldots & \kappa_{n}	\cr
                     \beta_1^{-n+L}\xi_{n-L}(\beta_1)	& \ldots & \beta_1^{-n}\xi_{n}(\beta_1)	\cr
                     \vdots		& \vdots & \vdots		\cr
                     \beta_L^{-n+L}\xi_{n-L}(\beta_L)	& \ldots & \beta_L^{-n}\xi_{n}(\beta_L)	\cr
              \end{pmatrix}} ,
\end{equation}
and
\begin{equation}
  \bar{\textsc{R}}_n = \frac{(-)^L}{\kappa_n}
   \frac{\det \begin{pmatrix}
                     \bar{\phi}_{n-L}(0)	& \ldots & \bar{\phi}_{n}(0)	\cr
                     \beta_1^{-n+L}\xi_{n-L}(\beta_1)	& \ldots & \beta_1^{-n}\xi_{n}(\beta_1)	\cr
                     \vdots		& \vdots & \vdots		\cr
                     \beta_L^{-n+L}\xi_{n-L}(\beta_L)	& \ldots & \beta_L^{-n}\xi_{n}(\beta_L)	\cr
              \end{pmatrix}}
        {\det \begin{pmatrix}
                     \beta_1^{-n+L}\xi_{n-L}(\beta_1)	& \ldots & \beta_1^{-n+1}\xi_{n-1}(\beta_1)	\cr
                     \vdots		& \vdots & \vdots		\cr
                     \beta_L^{-n+L}\xi_{n-L}(\beta_L)	& \ldots & \beta_L^{-n+1}\xi_{n-1}(\beta_L)	\cr
              \end{pmatrix}} .
\end{equation}
\end{proposition}
\begin{proof}
We begin with the derivation of (\ref{Lmod:b}). Writing
\begin{equation}
  \bar{\Phi}_n(z) = \sum^{n}_{j=0}b'_{n,j}\bar{\phi}_j(z) ,
\end{equation}
we have
\begin{equation}
  b'_{n,j} =  \int\frac{d\zeta}{2\pi i\zeta}w\left[\begin{array}{c} \ScSt \cdot \\ \ScSt L \end{array};\zeta \right]
              \prod^{L}_{l=1}(\zeta-\beta_l)\phi_j(\zeta)\bar{\Phi}_{n}(\bar{\zeta}) .
\end{equation}
We observe that $ \prod^{L}_{l=1}(\zeta-\beta_l)\phi_j(\zeta) \in \Pi_{j+L}(\zeta) $ 
and therefore $ b'_{n,j} = 0 $ for $ 0\leq j \leq n-1-L $. Now consider the
integral 
\begin{equation}
  0 =  \int\frac{d\zeta}{2\pi i\zeta}w\left[\begin{array}{c} \ScSt \cdot \\ \ScSt L \end{array};\zeta \right]
       \prod^{L}_{l=1\neq m}(\zeta-\beta_l)\bar{\Phi}_{n}(\bar{\zeta}) ,
\end{equation} 
for $ m=1,\ldots, L $ and subject to $ n\geq L $. This can be rewritten
\begin{align}
  0 & = \sum^{n}_{j=n-L}b'_{n,j} 
    \int\frac{d\zeta}{2\pi i\zeta}\frac{w(\zeta)}{\zeta-\beta_m}\bar{\phi}_{j}(\bar{\zeta}) ,
    \\
    & = -\sum^{n}_{j=n-L}b'_{n,j}\beta^{-j}_m\xi^*_j(\beta_m) ,
\end{align}
which yields $ L $ conditions. These enable us to solve for $ b'_{n,n-L}, \ldots, b'_{n,n-1} $
in terms of $ b'_{n,n} $ and we have the result (\ref{Lmod:b}). From this formula
we proceed to the derivation of (\ref{Lmod:d}). Using the definition (\ref{eps:b})
on this last result we require an evaluation of the integral
\begin{equation}
  \int\frac{d\zeta}{2\pi i\zeta}\frac{\zeta+z}{\zeta-z}w(\zeta)
       \frac{1}{\prod^{L}_{l=1}(\zeta-\beta_l)}\bar{\phi}_{n-j}(\bar{\zeta}) .
\end{equation}
We employ the partial fraction decomposition (assuming $ z \neq \beta_l $)
\begin{multline}
  \frac{\zeta+z}{\zeta-z}\frac{1}{\prod^{L}_{l=1}(\zeta-\beta_l)} 
 = \frac{\zeta+z}{\zeta-z}\frac{1}{\prod^{L}_{l=1}(z-\beta_l)} 
   - \frac{1}{\prod^{L}_{l=1}(z-\beta_l)} 
   \\
   + \sum^{L}_{m=1}\frac{1}{\zeta-\beta_m}\frac{\beta_m+z}{(\beta_m-z)\prod^{L}_{l=1\neq m}(\beta_m-\beta_l)} ,
\end{multline}
and evaluate the above integral as
\begin{equation} 
   -\frac{1}{\prod^{L}_{l=1}(z-\beta_l)}z^{-n+j}\xi^*_{n-j}(z) 
   + \sum^{L}_{m=1}\frac{\beta_m+z}{(\beta_m-z)\prod^{L}_{l=1\neq m}(\beta_m-\beta_l)}
                   \frac{1}{2}\beta^{-n+j-1}_m\xi^*_{n-j}(\beta_m) .
\end{equation}
When this result is inserted into the first row of the determinant formula 
(\ref{Lmod:b}) we see that the last term of this is a column independent
linear combination of the last $ L $ rows of the determinant. Thus we can eliminate
this term and arrive at (\ref{Lmod:d}).

To establish (\ref{Lmod:a}) we take an approach differing from the preceding
arguments. We postulate the determinant
\begin{equation}
   \det \begin{pmatrix}
                     \ldots & z^{L-j}\phi_{n-L+j}(z)	& \ldots 	\cr
                     \ldots & \beta_1^{-n+L-j}\xi_{n-L+j}(\beta_1)	& \ldots 	\cr
                            & \vdots		& \cr
                     \ldots & \beta_L^{-n+L-j}\xi_{n-L+j}(\beta_L)	& \ldots	\cr
        \end{pmatrix} ,
\label{star}
\end{equation}
and examine its orthogonality properties with respect to the weight
$ w\left[\begin{array}{c} \ScSt \cdot \\ \ScSt L \end{array};\zeta \right] $.
Let $ 0\leq m < n $ and $ 0\leq k\coloneqq L-j\leq L $ (we know $ n\geq k $), and
consider the integral
\begin{equation}
  \int\frac{d\zeta}{2\pi i\zeta}w\left[\begin{array}{c} \ScSt \cdot \\ \ScSt L \end{array};\zeta \right]
       \zeta^k\phi_{n-k}(\zeta)\bar{\zeta}^m
 =  \int\frac{d\zeta}{2\pi i\zeta}w(\zeta) \phi_{n-k}(\zeta)\frac{\zeta^{k-m}}{\prod^{L}_{l=1}(\zeta-\beta_l)} .
\end{equation}
In the case $ m>k $ we have the partial fraction decomposition
\begin{equation}
  \frac{1}{\zeta^N\prod^{L}_{l=1}(\zeta-\beta_l)} 
 = \frac{e_0}{\zeta^N}+ \cdots + e_N 
   + \sum^{L}_{p=1}\frac{1}{\zeta-\beta_p}\frac{\beta_p^{-N}}{\prod^{L}_{l=1\neq p}(\beta_p-\beta_l)} ,
   \quad N\geq 0 .
\end{equation}
This can be inserted into the preceding integral and, by orthogonality and (\ref{eps:a}),
we see the first $ N+1 $ terms vanish leaving
\begin{equation} 
   \sum^{L}_{p=1}\frac{\beta^{-m+k-1}_p}{\prod^{L}_{l=1\neq p}(\beta_p-\beta_l)}
                   \frac{1}{2}\xi_{n-k}(\beta_p), \quad n>m\geq k .
\end{equation}
In the remainder of cases, i.e. when $ 0\leq m\leq k $, we use the decomposition
\begin{equation}
  \frac{\zeta^N}{\prod^{L}_{l=1}(\zeta-\beta_l)} 
 = \delta_{N,L}
   + \sum^{L}_{p=1}\frac{1}{\zeta-\beta_p}\frac{\beta_p^{N}}{\prod^{L}_{l=1\neq p}(\beta_p-\beta_l)} ,
   \quad N\geq 0 .
\label{pf2}
\end{equation}
Utilising this form in the integral we are left with
\begin{equation} 
   \sum^{L}_{p=1}\frac{\beta^{k-m-1}_p}{\prod^{L}_{l=1\neq p}(\beta_p-\beta_l)}
                   \frac{1}{2}\xi_{n-k}(\beta_p), \quad n>m\geq k .
\end{equation}
In either case we conclude that integration of (\ref{star}) with the modified weight 
yields the determinant
\begin{equation}
   \det \begin{pmatrix}
                     \ldots & \sum^L_{p=1}f_p \beta_p^{-n+k}\xi_{n-k}(\beta_p)	& \ldots 	\cr
                     \ldots & \beta_1^{-n+k}\xi_{n-k}(\beta_1)	& \ldots 	\cr
                            & \vdots		& \cr
                     \ldots & \beta_L^{-n+k}\xi_{n-k}(\beta_L)	& \ldots	\cr
        \end{pmatrix} ,
\end{equation}
where the $ f_p $ are independent of the column index. Therefore the sum in the 
first row is a linear combination of the lower $ L $ rows and the determinant 
vanishes. Thus $ \Phi_n(z) $ is proportional to (\ref{star}) by (\ref{orthog:a}).
All that remains is to settle the normalisation.
To show that (\ref{Lmod:c}) is true we utilise the foregoing result. Here we are
lead to a determinant which has the following integral as entries of its first row, 
\begin{equation}
  \int\frac{d\zeta}{2\pi i\zeta}w\left[\begin{array}{c} \ScSt \cdot \\ \ScSt L \end{array};\zeta \right]
       \frac{\zeta+z}{\zeta-z}\zeta^k\phi_{n-k}(\zeta) .
\end{equation}
Then we can apply (\ref{pf2}) to $ \zeta^k/\prod^{L}_{l=1}(\zeta-\beta_l) $ and 
so obtain the evaluation
\begin{equation} 
   \frac{z^k}{\prod^{L}_{l=1}(z-\beta_l)}\xi_{n-k}(z)
  - \sum^{L}_{p=1}\frac{\beta^{n-1}_p}{\prod^{L}_{l=1\neq p}(\beta_p-\beta_l)}
                  \frac{z+\beta_p}{z-\beta_p}\frac{1}{2}\beta^{-n+k}_p\xi_{n-k}(\beta_p) .
\end{equation}
By an identical argument to that employed with (\ref{Lmod:a}) we conclude that only 
the first term survives in the determinant and (\ref{Lmod:c}) follows.
\end{proof}

We need to extend the preceding proposition to cover the case $ n < L $, which 
is the subject of the following proposition.
\begin{proposition}
For the reciprocal polynomial modification of the weight 
\begin{equation}
   w\left[\begin{array}{c} \ScSt \cdot \\ \ScSt L \end{array};z \right] 
   = w\left[\begin{array}{ccc} \ScSt \cdot & \ldots & \ScSt \cdot \\ \ScSt \beta_1 & \ldots & \ScSt \beta_{L}
            \end{array};z \right]
   = \frac{w(z)}{\prod^{L}_{j=1}(z-\beta_j)} ,
\end{equation}
with $ n < L $ the corresponding bi-orthogonal polynomials are given by
\begin{equation}
  \Phi_n(z) = C\cdot
   \det \begin{pmatrix}
         0 & \ldots & 0 & z^n\phi_{0}(z)	& \ldots & \phi_{n}(z)	\cr
         \beta_1^{L-n-1} & \ldots & 1 & \xi_{0}(\beta_1)	& \ldots & \beta_1^{-n}\xi_{n}(\beta_1)	\cr
         \vdots & & \vdots & \vdots		& \vdots & \vdots		\cr
         \beta_L^{L-n-1} & \ldots & 1 & \xi_{0}(\beta_L)	& \ldots & \beta_L^{-n}\xi_{n}(\beta_L)	\cr
              \end{pmatrix} ,
\label{Lmod:e}
\end{equation}
\begin{equation}
  \Phi^*_n(z) = C^*\cdot
   \det \begin{pmatrix}
         0 & \ldots & 0 & z^n\phi^*_{0}(z)	& \ldots & \phi^*_{n}(z)	\cr
         \beta_1^{L-n} & \ldots & \beta_1 & \xi^*_{0}(\beta_1)	& \ldots & \beta_1^{-n}\xi^*_{n}(\beta_1)	\cr
         \vdots & & \vdots & \vdots		& \vdots & \vdots		\cr
         \beta_L^{L-n} & \ldots & \beta_L & \xi^*_{0}(\beta_L)	& \ldots & \beta_L^{-n}\xi^*_{n}(\beta_L)	\cr
              \end{pmatrix} ,
\label{Lmod:f}
\end{equation}
and the associated functions
\begin{equation}
  \prod^{L}_{j=1}(z-\beta_j)\,\Xi_n(z) = C\cdot
   \det \begin{pmatrix}
         z^{L-1} & \ldots & z^n & z^n\xi_{0}(z)	& \ldots & \xi_{n}(z)	\cr
         \beta_1^{L-n-1} & \ldots & 1 & \xi_{0}(\beta_1)	& \ldots & \beta_1^{-n}\xi_{n}(\beta_1)	\cr
         \vdots & & \vdots & \vdots		& \vdots & \vdots		\cr
         \beta_L^{L-n-1} & \ldots & 1 & \xi_{0}(\beta_L)	& \ldots & \beta_L^{-n}\xi_{n}(\beta_L)	\cr
              \end{pmatrix} ,
\label{Lmod:g}
\end{equation}
\begin{equation}
  \prod^{L}_{j=1}(z-\beta_j)\,\Xi^*_n(z) = C^*\cdot
   \det \begin{pmatrix}
         z^{L} & \ldots & z^{n+1} & z^n\xi^*_{0}(z)	& \ldots & \xi^*_{n}(z)	\cr
         \beta_1^{L-n} & \ldots & \beta_1 & \xi^*_{0}(\beta_1)	& \ldots & \beta_1^{-n}\xi^*_{n}(\beta_1)	\cr
         \vdots & & \vdots & \vdots		& \vdots & \vdots		\cr
         \beta_L^{L-n} & \ldots & \beta_L & \xi^*_{0}(\beta_L)	& \ldots & \beta_L^{-n}\xi^*_{n}(\beta_L)	\cr
              \end{pmatrix} .
\label{Lmod:h}
\end{equation}
\end{proposition}
\begin{proof}
Orthogonality of the proposed formulae (\ref{Lmod:e}) and (\ref{Lmod:f}) with respect
to the modified weight according to the criteria (\ref{orthog:a}) and (\ref{orthog:b})
is easily proved using the methods of the preceding proposition and the 
interpolation identity
\begin{equation}
  \sum^L_{p=1} \frac{\beta^r_p}{\prod^L_{j\neq p}(\beta_p-\beta_j)} = 0, \quad 
  r=0,\ldots, L-2 .
\end{equation}
The remaining formulae (\ref{Lmod:g}) and (\ref{Lmod:h}) then follow from these two 
again using the methods of the preceding proposition.
\end{proof}

\begin{remark}
There is a simple observation that renders the above formulae (\ref{Lmod:e}-\ref{Lmod:h})
as transparently obvious - the inverse recurrence relation (\ref{Yrecur:b}) allows 
for a backwards recurrence and therefore gives meaning to the bi-orthogonal system
with negative indices $ n < 0 $. Recall that 
\begin{equation}
   Y_{0}
   = \begin{pmatrix}
              \kappa_0   &  \frac{\DySt \kappa_0(w_0+F)}{\DySt w} \cr
              \kappa_0   & -\frac{\DySt \kappa_0(w_0-F)}{\DySt w} \cr
     \end{pmatrix} ,
\label{}
\end{equation}
so that using (\ref{Yrecur:b}) we find $ \phi_{-1}(z) = \phi^*_{-1}(z) = 0 $
and $ \frac{1}{2}\kappa_{-1}\xi_{-1} = 1/z $ and 
$ \frac{1}{2}\kappa_{-1}\xi^*_{-1} = 1 $, assuming $ \kappa_{-1} \neq 0 $. 
Proceeding in this way we next find that $ \phi_{-2}(z) = \phi^*_{-2}(z) = 0 $
and in fact all the polynomials with negative argument can be set to zero.
The solutions of the backward recurrences for the associated functions can be
shown to be polynomials in $ 1/z $, namely with $ N \in \mathbb{N} $
\begin{align}
   \frac{1}{2}\kappa_{-N}\xi_{-N}(z) & = 
   \frac{\phi_{-N+1}(0)}{\kappa_{-N}}\frac{1}{z} + \cdots + \frac{1}{z^N} ,
   \\
   \frac{1}{2}\kappa_{-N}\xi^*_{-N}(z) & = 
   1 + \cdots + \frac{\bar{\phi}_{-N+1}(0)}{\kappa_{-N}}\frac{1}{z^{N-1}} ,
\end{align}
subject to the conditions $ \kappa_{-N}, \phi_{-N}(0), \bar{\phi}_{-N}(0) $
can be defined in some way and the $ \kappa_{-N} $ are non-zero. We note that the above formulae
are precisely the terminating forms of the general expansions 
(\ref{epsexp:a},\ref{epsexp:b},\ref{epsexp:c},\ref{epsexp:d}).
In defining the bi-orthogonal
system with negative indices in this way we have to break the relationship
between the polynomials and their associated functions as given in relations
(\ref{eps:a},\ref{eps:b},\ref{eps:c}).

Then using elementary column operations on (\ref{Lmod:e}-\ref{Lmod:h}) we can
recast these as
\begin{equation}
  \Phi_n(z) = C\cdot
   \det \begin{pmatrix}
         z^L\phi_{-(L-n)}(z) & \ldots & z^n\phi_{0}(z)	& \ldots & \phi_{n}(z)	\cr
         \beta_1^{L-n}\xi_{-(L-n)}(\beta_1) & \ldots & \xi_{0}(\beta_1)	& \ldots & \beta_1^{-n}\xi_{n}(\beta_1)	\cr
         \vdots & & \vdots		& \vdots & \vdots		\cr
         \beta_L^{L-n}\xi_{-(L-n)}(\beta_L) & \ldots & \xi_{0}(\beta_L) & \ldots & \beta_L^{-n}\xi_{n}(\beta_L)	\cr
              \end{pmatrix} ,
\label{}
\end{equation}
\begin{equation}
  \Phi^*_n(z) = C^*\cdot
   \det \begin{pmatrix}
         z^L\phi^*_{-(L-n)}(z) & \ldots & z^n\phi^*_{0}(z)	& \ldots & \phi^*_{n}(z)	\cr
         \beta_1^{L-n}\xi^*_{-(L-n)}(\beta_1) & \ldots & \xi^*_{0}(\beta_1)	& \ldots & \beta_1^{-n}\xi^*_{n}(\beta_1)	\cr
         \vdots & & \vdots		& \vdots & \vdots		\cr
         \beta_L^{L-n}\xi^*_{-(L-n)}(\beta_L) & \ldots & \xi^*_{0}(\beta_L)	& \ldots & \beta_L^{-n}\xi^*_{n}(\beta_L)	\cr
              \end{pmatrix} ,
\label{}
\end{equation}
and the associated functions
\begin{multline}
  \prod^{L}_{j=1}(z-\beta_j)\,\Xi_n(z) \\
  = C\cdot
   \det \begin{pmatrix}
         z^{L}\xi_{-(L-n)}(z) & \ldots & z^{n+1}\xi_{-1}(z) & z^n\xi_{0}(z)	& \ldots & \xi_{n}(z)	\cr
         \beta_1^{L-n}\xi_{-(L-n)}(\beta_1) & \ldots & \beta_1\xi_{-1}(\beta_1) & \xi_{0}(\beta_1) & \ldots & \beta_1^{-n}\xi_{n}(\beta_1)	\cr
         \vdots & & \vdots & \vdots		& \vdots & \vdots		\cr
         \beta_L^{L-n}\xi_{-(L-n)}(\beta_L) & \ldots & \beta_L\xi_{-1}(\beta_L) & \xi_{0}(\beta_L) & \ldots & \beta_L^{-n}\xi_{n}(\beta_L)	\cr
              \end{pmatrix} ,
\label{}
\end{multline}
\begin{multline}
  \prod^{L}_{j=1}(z-\beta_j)\,\Xi^*_n(z) \\
  = C^*\cdot
   \det \begin{pmatrix}
         z^{L}\xi^*_{-(L-n)}(z) & \ldots & z^{n+1}\xi^*_{-1}(z) & z^n\xi^*_{0}(z) & \ldots & \xi^*_{n}(z)	\cr
         \beta_1^{L-n}\xi^*_{-(L-n)}(\beta_1) & \ldots & \beta_1\xi^*_{-1}(\beta_1) & \xi^*_{0}(\beta_1) & \ldots & \beta_1^{-n}\xi^*_{n}(\beta_1)	\cr
         \vdots & & \vdots & \vdots		& \vdots & \vdots		\cr
         \beta_L^{L-n}\xi^*_{-(L-n)}(\beta_L) & \ldots & \beta_L\xi^*_{-1}(\beta_L) & \xi^*_{0}(\beta_L) & \ldots & \beta_L^{-n}\xi^*_{n}(\beta_L)	\cr
              \end{pmatrix} .
\label{}
\end{multline}
\end{remark}

Propositions \ref{Kmod} and \ref{Lmod} can be combined into one result covering 
the case of rational modification of the weight.
\begin{corollary}\label{Rmod}
For a general rational modification, subject to $ n \geq L $,
\begin{equation}
   w\left[\begin{array}{c} \ScSt K \\ \ScSt L \end{array};z \right]
    = \frac{\prod^{K}_{j=1}(z-\alpha_j)}{\prod^{L}_{j=1}(z-\beta_j)} w(z) ,
\end{equation}
we have the corresponding bi-orthogonal system
\begin{equation}
  \prod^{K}_{j=1}(z-\alpha_j)\,\Phi_n(z) =
  C\cdot
  \det \begin{pmatrix}
         z^L\phi_{n-L}(z)	& \ldots & \phi_{n}(z)	& \ldots & \phi_{n+K}(z)	\cr
         \alpha_1^L\phi_{n-L}(\alpha_1)	& \ldots & \phi_{n}(\alpha_1)	& \ldots & \phi_{n+K}(\alpha_1)	\cr
         \vdots	&  & \vdots		&  & \vdots		\cr
         \alpha_K^L\phi_{n-L}(\alpha_K)	& \ldots & \phi_{n}(\alpha_K)	& \ldots & \phi_{n+K}(\alpha_K)	\cr
         \beta_1^{-n+L}\xi_{n-L}(\beta_1)	& \ldots & \beta_1^{-n}\xi_{n}(\beta_1)	& \ldots & \beta_1^{-n}\xi_{n+K}(\beta_1) \cr
         \vdots		&  & \vdots	&  & \vdots		\cr
         \beta_L^{-n+L}\xi_{n-L}(\beta_L)	& \ldots & \beta_L^{-n}\xi_{n}(\beta_L)	& \ldots & \beta_L^{-n}\xi_{n+K}(\beta_L) \cr
       \end{pmatrix} ,
\label{Rmod:a}
\end{equation}
\begin{equation}
  \prod^{K}_{j=1}(z-\alpha_j)\,\Phi^*_n(z) =
  C^*\cdot
  \det \begin{pmatrix}
         z^L\phi^*_{n-L}(z)	& \ldots & \phi^*_{n}(z)	& \ldots & \phi^*_{n+K}(z)		\cr
         \alpha_1^L\phi^*_{n-L}(\alpha_1)	& \ldots & \phi^*_{n}(\alpha_1)	& \ldots & \phi^*_{n+K}(\alpha_1)	\cr
         \vdots			&	 & \vdots 		& 	 & \vdots			\cr
         \alpha_K^L\phi^*_{n-L}(\alpha_K)	& \ldots & \phi^*_{n}(\alpha_K)	& \ldots & \phi^*_{n+K}(\alpha_K)	\cr
         \beta_1^{-n+L}\xi^*_{n-L}(\beta_1)& \ldots & \beta_1^{-n}\xi^*_{n}(\beta_1)	& \ldots & \beta_1^{-n}\xi^*_{n+K}(\beta_1)\cr
         \vdots			& 	 & \vdots 		& 	 & \vdots		\cr
         \beta_L^{-n+L}\xi^*_{n-L}(\beta_L)& \ldots & \beta_L^{-n}\xi^*_{n}(\beta_L)	& \ldots & \beta_L^{-n}\xi^*_{n+K}(\beta_L)\cr
       \end{pmatrix} ,
\label{Rmod:b}
\end{equation}
\begin{equation}
  \prod^{L}_{j=1}(z-\beta_j)\,\Xi_n(z) =
  C\cdot
  \det \begin{pmatrix}
         z^L\xi_{n-L}(z)		& \ldots & \xi_{n}(z)	& \ldots & \xi_{n+K}(z)	\cr
         \alpha_1^L\phi_{n-L}(\alpha_1)	& \ldots & \phi_{n}(\alpha_1)	& \ldots & \phi_{n+K}(\alpha_1)	\cr
         \vdots	&  & \vdots		&  & \vdots		\cr
         \alpha_K^L\phi_{n-L}(\alpha_K)	& \ldots & \phi_{n}(\alpha_K)	& \ldots & \phi_{n+K}(\alpha_K)	\cr
         \beta_1^{-n+L}\xi_{n-L}(\beta_1)	& \ldots & \beta_1^{-n}\xi_{n}(\beta_1)	& \ldots & \beta_1^{-n}\xi_{n+K}(\beta_1) \cr
         \vdots		&  & \vdots	&  & \vdots		\cr
         \beta_L^{-n+L}\xi_{n-L}(\beta_L)	& \ldots & \beta_L^{-n}\xi_{n}(\beta_L)	& \ldots & \beta_L^{-n}\xi_{n+K}(\beta_L) \cr
       \end{pmatrix} ,
\label{Rmod:c}
\end{equation}
\begin{equation}
  \prod^{L}_{j=1}(z-\beta_j)\,\Xi^*_n(z) =
  C^*\cdot
  \det \begin{pmatrix}
         z^{L}\xi^*_{n-L}(z)	& \ldots & \xi^*_{n}(z)	& \ldots & \xi^*_{n+K}(z)	\cr
         \alpha_1^{L}\phi^*_{n-L}(\alpha_1)	 & \ldots & \phi^*_{n}(\alpha_1)	& \ldots & \phi^*_{n+K}(\alpha_1)	\cr
         \vdots	&  & \vdots 		&  & \vdots		\cr
         \alpha_K^{L}\phi^*_{n-L}(\alpha_K)	 & \ldots & \phi^*_{n}(\alpha_K)	& \ldots & \phi^*_{n+K}(\alpha_K)	\cr
         \beta_1^{-n+L}\xi^*_{n-L}(\beta_1) & \ldots & \beta_1^{-n}\xi^*_{n}(\beta_1) & \ldots & \beta_1^{-n}\xi^*_{n+K}(\beta_1)	\cr
         \vdots	&  & \vdots 		&  & \vdots		\cr
         \beta_L^{-n+L}\xi^*_{n-L}(\beta_L) & \ldots & \beta_L^{-n}\xi^*_{n}(\beta_L)	& \ldots & \beta_L^{-n}\xi^*_{n+K}(\beta_L)	\cr
       \end{pmatrix} .
\label{Rmod:d}
\end{equation}
\end{corollary}
\begin{proof}
This follows by using the methods of Propositions \ref{Kmod} and \ref{Lmod}.
\end{proof}

We now present a "conjugated" version of Proposition \ref{Kmod}.
\begin{proposition}\label{Kmod*}
For the polynomial modification of the weight
\begin{equation}
   w\left[\begin{array}{ccc} \cdot & ; & \ScSt K^* \\ \cdot & ; & \ScSt \cdot \end{array};z \right]
   = \prod^{K^*}_{j=1}(1-\alpha^*_jz^{-1}) w(z) ,
\end{equation}
the corresponding bi-orthogonal polynomials are given by
\begin{equation}
  \prod^{K^*}_{j=1}(1-\alpha^*_j z^{-1})\,\Phi_n(z) = C\cdot
   \det \begin{pmatrix}
          \phi_{n}(z)	& \ldots & z^{-K^*}\phi_{n+K^*}(z)	\cr
          \phi_{n}(\alpha^*_1)	& \ldots & (\alpha^*_1)^{-K^*}\phi_{n+K^*}(\alpha^*_1)	\cr
          \vdots		& \vdots & \vdots		\cr
          \phi_{n}(\alpha^*_{K^*})	& \ldots & (\alpha^*_{K^*})^{-K^*}\phi_{n+K^*}(\alpha^*_{K^*})	\cr
        \end{pmatrix} ,
\label{Kmod*:a}
\end{equation}
\begin{equation}
  \prod^{K^*}_{j=1}(1-\alpha^*_j z^{-1})\,\Phi^*_n(z) = C^*\cdot
   \det \begin{pmatrix}
          \phi^*_{n}(z)	& \ldots & z^{-K^*}\phi^*_{n+K^*}(z)	\cr
          \phi^*_{n}(\alpha^*_1)	& \ldots & (\alpha^*_1)^{-K^*}\phi^*_{n+K^*}(\alpha^*_1)	\cr
          \vdots		& \vdots & \vdots		\cr
          \phi^*_{n}(\alpha^*_{K^*})	& \ldots & (\alpha^*_{K^*})^{-K^*}\phi^*_{n+K^*}(\alpha^*_{K^*})	\cr
        \end{pmatrix} ,
\label{Kmod*:b}
\end{equation}
and the associated functions
\begin{equation}
   \Xi_n(z) = C\cdot
   \det \begin{pmatrix}
          \xi_{n}(z)	& \ldots & z^{-K^*}\xi_{n+K^*}(z)	\cr
          \phi_{n}(\alpha^*_1)	& \ldots & (\alpha^*_1)^{-K^*}\phi_{n+K^*}(\alpha^*_1)	\cr
          \vdots		& \vdots & \vdots		\cr
          \phi_{n}(\alpha^*_{K^*})	& \ldots & (\alpha^*_{K^*})^{-K^*}\phi_{n+K^*}(\alpha^*_{K^*})	\cr
        \end{pmatrix} ,
\label{Kmod*:c}
\end{equation}
\begin{equation}
   \Xi^*_n(z) = C^*\cdot
   \det \begin{pmatrix}
          \xi^*_{n}(z)	& \ldots & z^{-K^*}\xi^*_{n+K^*}(z)	\cr
          \phi^*_{n}(\alpha^*_1)	& \ldots & (\alpha^*_1)^{-K^*}\phi^*_{n+K^*}(\alpha^*_1)	\cr
          \vdots		& \vdots & \vdots		\cr
          \phi^*_{n}(\alpha^*_{K^*})	& \ldots & (\alpha^*_{K^*})^{-K^*}\phi^*_{n+K^*}(\alpha^*_{K^*})	\cr
        \end{pmatrix} .
\label{Kmod*:d}
\end{equation}
\end{proposition}
\begin{proof}
This result follows by employing the methods of Proposition \ref{Kmod}.
\end{proof}

There is also a "conjugated" analogue of Proposition \ref{Lmod} and is given in
the following statement.
\begin{proposition}\label{Lmod*}
For the reciprocal polynomial modification of the weight
\begin{equation}
   w\left[\begin{array}{ccc} \cdot & ; & \ScSt \cdot \\ \cdot & ; & \ScSt L^* \end{array};z \right]
    = \frac{1}{\prod^{L^*}_{j=1}(1-\beta^*_jz^{-1})} w(z) ,
\end{equation}
the corresponding bi-orthogonal polynomials are given by
\begin{equation}
  \Phi_n(z) = C\cdot
   \det \begin{pmatrix}
          \phi_{n-L^*}(z)	& \ldots & \phi_{n}(z)	\cr
          (\beta^*_1)^{-n}\xi_{n-L^*}(\beta^*_1)	& \ldots & (\beta^*_1)^{-n}\xi_{n}(\beta^*_1)	\cr
          \vdots		& \vdots & \vdots		\cr
          (\beta^*_{L^*})^{-n}\xi_{n-L^*}(\beta^*_{L^*})	& \ldots & (\beta^*_{L^*})^{-n}\xi_{n}(\beta^*_{L^*})	\cr
        \end{pmatrix}, \quad n\geq L^*,
\label{Lmod*:a}
\end{equation}
\begin{equation}
  \Phi^*_n(z) = C^*\cdot
   \det \begin{pmatrix}
          \phi^*_{n-L^*}(z)	& \ldots & \phi^*_{n}(z)	\cr
          (\beta^*_1)^{-n}\xi^*_{n-L^*}(\beta^*_1)	& \ldots & (\beta^*_1)^{-n}\xi^*_{n}(\beta^*_1)	\cr
          \vdots		& \vdots & \vdots		\cr
          (\beta^*_{L^*})^{-n}\xi^*_{n-L^*}(\beta^*_{L^*})	& \ldots & (\beta^*_{L^*})^{-n}\xi^*_{n}(\beta^*_{L^*})	\cr
        \end{pmatrix}, \quad n\geq L^*,
\label{Lmod*:b}
\end{equation}
and the associated functions
\begin{equation}
  \prod^{L^*}_{j=1}(1-\beta^*_j z^{-1})\,\Xi_n(z) = C\cdot
   \det \begin{pmatrix}
          \xi_{n-L^*}(z)	& \ldots & \xi_{n}(z)	\cr
          (\beta^*_1)^{-n}\xi_{n-L^*}(\beta^*_1)	& \ldots & (\beta^*_1)^{-n}\xi_{n}(\beta^*_1)	\cr
          \vdots		& \vdots & \vdots		\cr
          (\beta^*_{L^*})^{-n}\xi_{n-L^*}(\beta_{L^*})	& \ldots & (\beta^*_{L^*})^{-n}\xi_{n}(\beta^*_{L^*})	\cr
        \end{pmatrix}, \quad n\geq L^*,
\label{Lmod*:c}
\end{equation}
\begin{equation}
  \prod^{L^*}_{j=1}(1-\beta^*_j z^{-1})\,\Xi^*_n(z) = C^*\cdot
   \det \begin{pmatrix}
          \xi^*_{n-L^*}(z)	& \ldots & \xi^*_{n}(z)	\cr
          (\beta^*_1)^{-n}\xi^*_{n-L^*}(\beta^*_1)	& \ldots & (\beta^*_1)^{-n}\xi^*_{n}(\beta^*_1)	\cr
          \vdots		& \vdots & \vdots		\cr
          (\beta^*_{L^*})^{-n}\xi^*_{n-L^*}(\beta^*_{L^*})	& \ldots & (\beta^*_{L^*})^{-n}\xi^*_{n}(\beta^*_{L^*})	\cr
        \end{pmatrix}, \quad n\geq L^*.
\label{Lmod*:d}
\end{equation}
\end{proposition}
\begin{proof}
This result follows by employing the methods of Proposition \ref{Lmod}.
\end{proof}

Analogous to Corollary \ref{Rmod}, Propositions \ref{Kmod*} and \ref{Lmod*} can 
be combined into one result covering the case of a "conjugated" rational modification 
of the weight.
\begin{corollary}\label{Rmod*}
For a general rational modification, subject to $ n \geq L^* $,
\begin{equation}
   w\left[\begin{array}{ccc} \cdot & ; & \ScSt K^* \\ \cdot & ; & \ScSt L^* \end{array};z \right]
    = \frac{\prod^{K^*}_{j=1}(1-\alpha^*_jz^{-1})}{\prod^{L^*}_{j=1}(1-\beta^*_jz^{-1})} w(z) ,
\end{equation}
we have the new bi-orthogonal system
\begin{multline}
  \prod^{K^*}_{j=1}(1-\alpha^*_j z^{-1})\,\Phi_n(z) \\
 = C\cdot
  \det \begin{pmatrix}
         \phi_{n-L^*}(z)	& \ldots & \phi_{n}(z)	& \ldots & z^{-K^*}\phi_{n+K^*}(z)	\cr
         \phi_{n-L^*}(\alpha^*_1)	& \ldots & \phi_{n}(\alpha^*_1)	& \ldots & (\alpha^*_1)^{-K^*}\phi_{n+K^*}(\alpha^*_1)	\cr
         \vdots	&  & \vdots		&  & \vdots		\cr
         \phi_{n-L^*}(\alpha^*_{K^*})	& \ldots & \phi_{n}(\alpha^*_{K^*})	& \ldots & (\alpha^*_{K^*})^{-K^*}\phi_{n+K^*}(\alpha^*_{K^*})	\cr
         (\beta^*_1)^{-n}\xi_{n-L^*}(\beta^*_1)	& \ldots & (\beta^*_1)^{-n}\xi_{n}(\beta^*_1) & \ldots & (\beta^*_1)^{-n-K^*}\xi_{n+K^*}(\beta^*_1) \cr
         \vdots		&  & \vdots	&  & \vdots		\cr
         (\beta^*_{L^*})^{-n}\xi_{n-L^*}(\beta^*_{L^*})	& \ldots & (\beta^*_{L^*})^{-n}\xi_{n}(\beta^*_{L^*}) & \ldots & (\beta^*_{L^*})^{-n-K^*}\xi_{n+K^*}(\beta^*_{L^*}) \cr
       \end{pmatrix} ,
\label{Rmod*:a}
\end{multline}
\begin{multline}
  \prod^{K^*}_{j=1}(1-\alpha^*_j z^{-1})\,\Phi^*_n(z) \\
 = C^*\cdot
  \det \begin{pmatrix}
         \phi^*_{n-L^*}(z)	& \ldots & \phi^*_{n}(z)	& \ldots & z^{-K^*}\phi^*_{n+K^*}(z)		\cr
         \phi^*_{n-L^*}(\alpha^*_1)	& \ldots & \phi^*_{n}(\alpha^*_1)	& \ldots & (\alpha^*_1)^{-K^*}\phi^*_{n+K^*}(\alpha^*_1)	\cr
         \vdots			&	 & \vdots 		& 	 & \vdots			\cr
         \phi^*_{n-L^*}(\alpha^*_{K^*})	& \ldots & \phi^*_{n}(\alpha^*_{K^*})	& \ldots & (\alpha^*_{K^*})^{-K^*}\phi^*_{n+K^*}(\alpha^*_{K^*})	\cr
         (\beta^*_1)^{-n}\xi^*_{n-L^*}(\beta^*_1)& \ldots & (\beta^*_1)^{-n}\xi^*_{n}(\beta^*_1)	& \ldots & (\beta^*_1)^{-n-K^*}\xi^*_{n+K^*}(\beta^*_1)\cr
         \vdots			& 	 & \vdots 		& 	 & \vdots		\cr
         (\beta^*_{L^*})^{-n}\xi^*_{n-L^*}(\beta^*_{L^*})& \ldots & (\beta^*_{L^*})^{-n}\xi^*_{n}(\beta^*_{L^*}) & \ldots & (\beta^*_{L^*})^{-n-K^*}\xi^*_{n+K^*}(\beta^*_{L^*})\cr
       \end{pmatrix} ,
\label{Rmod*:b}
\end{multline}
\begin{multline}
  \prod^{L^*}_{j=1}(1-\beta^*_j z^{-1})\,\Xi_n(z) \\
 = C\cdot
  \det \begin{pmatrix}
         \xi_{n-L^*}(z)		& \ldots & \xi_{n}(z)	& \ldots & z^{-K^*}\xi_{n+K^*}(z)	\cr
         \phi_{n-L^*}(\alpha^*_1)	& \ldots & \phi_{n}(\alpha^*_1)	& \ldots & (\alpha^*_1)^{-K^*}\phi_{n+K^*}(\alpha^*_1)	\cr
         \vdots	&  & \vdots		&  & \vdots		\cr
         \phi_{n-L^*}(\alpha^*_{K^*})	& \ldots & \phi_{n}(\alpha^*_{K^*})	& \ldots & (\alpha^*_{K^*})^{-K^*}\phi_{n+K^*}(\alpha^*_{K^*})	\cr
         (\beta^*_1)^{-n}\xi_{n-L^*}(\beta^*_1)	& \ldots & (\beta^*_1)^{-n}\xi_{n}(\beta^*_1)	& \ldots & (\beta^*_1)^{-n-K^*}\xi_{n+K^*}(\beta^*_1) \cr
         \vdots		&  & \vdots	&  & \vdots		\cr
         (\beta^*_{L^*})^{-n}\xi_{n-L^*}(\beta^*_{L^*})	& \ldots & (\beta^*_{L^*})^{-n}\xi_{n}(\beta^*_{L^*})	& \ldots & (\beta^*_{L^*})^{-n-K^*}\xi_{n+K^*}(\beta^*_{L^*}) \cr
       \end{pmatrix} ,
\label{Rmod*:c}
\end{multline}
\begin{multline}
  \prod^{L^*}_{j=1}(1-\beta^*_j z^{-1}) \,\Xi^*_n(z) \\
 = C^*\cdot
  \det \begin{pmatrix}
         \xi^*_{n-L^*}(z)	& \ldots & \xi^*_{n}(z)	& \ldots & z^{-K^*}\xi^*_{n+K^*}(z)	\cr
         \phi^*_{n-L^*}(\alpha^*_1)	 & \ldots & \phi^*_{n}(\alpha^*_1)	& \ldots & (\alpha^*_1)^{-K^*}\phi^*_{n+K^*}(\alpha^*_1)	\cr
         \vdots	&  & \vdots 		&  & \vdots		\cr
         \phi^*_{n-L^*}(\alpha^*_{K^*})	 & \ldots & \phi^*_{n}(\alpha^*_{K^*})	& \ldots & (\alpha^*_{K^*})^{-K^*}\phi^*_{n+K^*}(\alpha^*_{K^*})	\cr
         (\beta^*_1)^{-n}\xi^*_{n-L^*}(\beta^*_1) & \ldots & (\beta^*_1)^{-n}\xi^*_{n}(\beta^*_1) & \ldots & (\beta^*_1)^{-n-K^*}\xi^*_{n+K^*}(\beta^*_1)	\cr
         \vdots	&  & \vdots 		&  & \vdots		\cr
         (\beta^*_{L^*})^{-n}\xi^*_{n-L^*}(\beta^*_{L^*}) & \ldots & (\beta^*_{L^*})^{-n}\xi^*_{n}(\beta^*_{L^*})	& \ldots & (\beta^*_{L^*})^{-n-K^*}\xi^*_{n+K^*}(\beta^*_{L^*})	\cr
       \end{pmatrix} .
\label{Rmod*:d}
\end{multline}
\end{corollary}
\begin{proof}
These results may be derived by utilising the methods of Propositions \ref{Kmod}
and \ref{Lmod}.
\end{proof}

Corollaries \ref{Rmod} and \ref{Rmod*} can be combined into a single result for
the most general rational modification of the weight (\ref{C-Uxfm:b}) and this is 
stated in (\ref{KK*LL*shift}) without proof.

\begin{remark}
We now show how the results of Ismail and Ruedemann \cite{IR_1992} 
and of Ruedemann \cite{Ru_1994} can be recovered from a specialisation 
of the result for $ K,K^* > 0 $, which follows easily by combining (\ref{Kmod:a}) and (\ref{Kmod*:a}),
\begin{multline}
  \prod^{K}_{j=1}(z-\alpha_j)\prod^{K^*}_{j=1}(z-\alpha^*_j)\,\Phi_n(z) \\
 = C\cdot
   \det \begin{pmatrix}
          z^{K^*}\phi_{n}(z)	& \ldots & \phi_{n+K^*}(z)	& \ldots & \phi_{n+K^*+K}(z) \cr
          (\alpha_1)^{K^*}\phi_{n}(\alpha_1)	& \ldots & \phi_{n+K^*}(\alpha_1)	& \ldots & \phi_{n+K^*+K}(\alpha_1)	\cr
          \vdots		& \vdots & \vdots & \vdots & \vdots	\cr
          (\alpha_{K})^{K^*}\phi_{n}(\alpha_{K})	& \ldots & \phi_{n+K^*}(\alpha_{K}) & \ldots & \phi_{n+K^*+K}(\alpha_K)	\cr
          (\alpha^*_1)^{K^*}\phi_{n}(\alpha^*_1)	& \ldots & \phi_{n+K^*}(\alpha^*_1)	& \ldots & \phi_{n+K^*+K}(\alpha^*_1)	\cr
          \vdots		& \vdots & \vdots & \vdots & \vdots	\cr
          (\alpha^*_{K^*})^{K^*}\phi_{n}(\alpha^*_{K^*})	& \ldots & \phi_{n+K^*}(\alpha^*_{K^*}) & \ldots & \phi_{n+K^*+K}(\alpha^*_{K^*})	\cr
        \end{pmatrix}, \quad n\geq 0.
\label{KK*mod:a}
\end{multline}
The other result of \cite{IR_1992} is a likewise a special case of $ L,L^* > 0 $
\begin{multline}
  \Phi_n(z) \\
 = C\cdot
   \det \begin{pmatrix}
           z^L\phi_{n-L-L^*}(z)	& \ldots & \phi_{n-L^*}(z)	& \ldots & \phi_{n}(z)	\cr
          \beta_1^{-n+L}\xi_{n-L-L^*}(\beta_1)	& \ldots & (\beta_1)^{-n}\xi_{n-L^*}(\beta_1)	& \ldots & (\beta_1)^{-n}\xi_{n}(\beta_1)	\cr
          \vdots		& \vdots & \vdots		& \vdots & \vdots		\cr
          \beta_L^{-n+L}\xi_{n-L-L^*}(\beta_L)	& \ldots & \beta_{L}^{-n}\xi_{n-L^*}(\beta_{L})	& \ldots & \beta_{L}^{-n}\xi_{n}(\beta_{L})	\cr
          (\beta^*_1)^{-n+L}\xi_{n-L-L^*}(\beta^*_1)	& \ldots & (\beta^*_1)^{-n}\xi_{n-L^*}(\beta^*_1)	& \ldots & (\beta^*_1)^{-n}\xi_{n}(\beta^*_1)	\cr
          \vdots		& \vdots & \vdots		& \vdots & \vdots		\cr
          (\beta^*_{L^*})^{-n+L}\xi_{n-L-L^*}(\beta^*_{L^*})	& \ldots & (\beta^*_{L^*})^{-n}\xi_{n-L^*}(\beta^*_{L^*})	& \ldots & (\beta^*_{L^*})^{-n}\xi_{n}(\beta^*_{L^*})	\cr
        \end{pmatrix}, \quad n\geq L+L^*.
\label{LL*mod:a}
\end{multline}
In the context of orthogonal polynomials Theorem 1 of \cite{IR_1992} states a formula
which is extended to the bi-orthogonal polynomial setting by Theorem 2 of \cite{Ru_1994}.
This is the $ K=K^*=m $ case in (\ref{KK*mod:a}) and we now proceed to show that our 
result is equivalent to theirs. Our result can be expressed, keeping only the first row, as
\begin{equation*}
    \det ( z^m\phi_n(z), \ldots, z\phi_{n+m-1}(z), \, \phi_{n+m}(z), \, \phi_{n+m+1}(z), \ldots, \phi_{n+2m}(z) ) ,
\end{equation*}
and we use elementary column operations and the forward and backward recurrence 
relations (\ref{Yrecur:a}) and (\ref{Yrecur:b}). The initial step is to recast this as
\begin{equation*}
    \det ( z^{m-1}\phi^*_{n+1}, \ldots, \phi^*_{n+m}, \, \phi_{n+m}, \, \phi_{n+m+1}, \ldots, \phi_{n+2m} ) ,
\end{equation*}
and then conduct a sequence of $ m $ consecutive sweeps. The first sweep uses 
$ \phi^*_{n+m} $ and those on and to the right of $ \phi_{n+m+1} $ in a sequence 
of successive steps. At each of these steps a recurrence relation is used to transform
the first element by incrementing its index and then another is used to decrement
the index of the second and thereby introducing a factor of $ z $.
(we use $ \cdot $ to denote an intervening element not involved in the column operations)
\begin{gather*}
    \det ( z^{m-1}\phi^*_{n+1}, \ldots, \phi^*_{n+m}, \cdot, \phi_{n+m+1}, \ldots, \phi_{n+2m} ) ,
  \\ \downarrow \\
    \det ( z^{m-1}\phi^*_{n+1}, \ldots, \phi^*_{n+m+1}, \cdot, \phi_{n+m+1}, \ldots, \phi_{n+2m} ) ,
  \\ \downarrow \\
    \det ( z^{m-1}\phi^*_{n+1}, \ldots, \phi^*_{n+m+1}, \cdot, z\phi_{n+m}, \ldots, \phi_{n+2m} ) ,
  \\ \downarrow \\
    \det ( z^{m-1}\phi^*_{n+1}, \ldots, \phi^*_{n+m+2}, \cdot, z\phi_{n+m}, \phi_{n+m+2}, \ldots, \phi_{n+2m} ) ,
  \\ \downarrow \\
    \det ( z^{m-1}\phi^*_{n+1}, \ldots, \phi^*_{n+m+2}, \cdot, z\phi_{n+m}, z\phi_{n+m+1}, \ldots, \phi_{n+2m} ) ,
  \\ \vdots \\
    \det ( z^{m-1}\phi^*_{n+1}, \ldots, \phi^*_{n+2m-1}, \cdot, z\phi_{n+m}, \ldots, z\phi_{n+2m-1} ) .
\end{gather*}
Then we successively use a recurrence involving the first element and the second in 
reverse direction in order to decrement its index, restoring the index to $ n+m $.
\begin{gather*}
    \det ( z^{m-1}\phi^*_{n+1}, \ldots, \phi^*_{n+2m-1}, \phi_{n+m}, z\phi_{n+m}, \ldots, z\phi_{n+2m-1} ) ,
  \\ \vdots \\
    \det ( z^{m-1}\phi^*_{n+1}, \ldots, \phi^*_{n+m}, \phi_{n+m}, z\phi_{n+m}, \ldots, z\phi_{n+2m-1} ) ,
\end{gather*} 
The second sweep uses $ z\phi^*_{n+m-1} $ and those of $ z\phi_{n+m+1} $ and the ones 
to the right in a similar sequence of steps 
 \begin{gather*}
    \det ( z^{m-1}\phi^*_{n+1}, \ldots, z\phi^*_{n+m-1}, \cdot, \cdot, z\phi_{n+m}, z\phi_{n+m+1}, \ldots, z\phi_{n+2m-1} ) ,
  \\ \downarrow \\
    \det ( z^{m-1}\phi^*_{n+1}, \ldots, z\phi^*_{n+m}, \cdot, \cdot, z\phi_{n+m}, z\phi_{n+m+1}, \ldots, z\phi_{n+2m-1} ) ,
  \\ \downarrow \\
    \det ( z^{m-1}\phi^*_{n+1}, \ldots, z\phi^*_{n+m}, \cdot, \cdot, z\phi_{n+m}, z^2\phi_{n+m}, \ldots, z\phi_{n+2m-1} ) ,
  \\ \vdots \\
    \det ( z^{m-1}\phi^*_{n+1}, \ldots, z\phi^*_{n+2m-2}, \cdot, \cdot, z\phi_{n+m}, z^2\phi_{n+m}, z^2\phi_{n+m+1}, \ldots, z^2\phi_{n+2m-2} ) ,
  \\ \vdots \\
    \det ( z^{m-1}\phi^*_{n+1}, \ldots, z\phi^*_{n+m}, \cdot, \cdot, z\phi_{n+m}, z^2\phi_{n+m}, \ldots, z^2\phi_{n+2m-2} ) ,
\end{gather*}
We repeat this procedure $m$ times and the final sweep yields
 \begin{gather*}
    \det ( z^{m-1}\phi^*_{n+m}, \ldots, z\phi^*_{n+m}, \phi^*_{n+m}, \phi_{n+m}, z\phi_{n+m}, z^2\phi_{n+m}, \ldots, z^m\phi_{n+m} ) ,
\end{gather*}
which is the desired expression.
\end{remark}

We now specialise the above results and consider the $ K=1,L=0 $ and $ K=0,L=1 $ 
cases in further detail as these particular transformations form the generators
for the group of Christoffel-Geronimus-Uvarov transformations. To this purpose we define the matrix generators
$ R_n^+(z;\alpha) $ and $ R_n^-(z;\beta) $ for the $ K=1,L=0 $ and $ K=0,L=1 $
cases respectively
\begin{align} 
   Y^+_n(z) = Y_n\left[\begin{array}{c} \ScSt 1 \\ \ScSt 0 \end{array};z \right]
   & = R_n^+(z;\alpha)Y_n(z) ,
   \label{Y_CUxfm:a} \\
   Y^-_n(z) = Y_n\left[\begin{array}{c} \ScSt 0 \\ \ScSt 1 \end{array};z \right]
   & = R_n^-(z;\beta)Y_n(z) .
   \label{Y_CUxfm:b}
\end{align} 
In our specialisation to the regular semi-classical weights we will draw heavily
on the results for the $ K=1,L=0 $ and $ K=0,L=1 $ cases as they will form the
generators of the Schlesinger transformations in that context. 

\begin{corollary}\label{Kshift}
For $ K=1,L=0 $ with the numerator factor of $ z-\alpha $ denote the coefficients 
of the transformed BOPS by $ \kappa^{+}_n, r^{+}_n, \bar{r}^{+}_n $ and the 
bi-orthogonal polynomials and associated functions by 
$ \Phi_{n}, \Phi^*_{n},  \Xi_{n},  \Xi^*_{n} $. 
Assuming $ \phi_n(\alpha) \neq 0 $ for $ n\in \Z_{\geq 0} $ these are given in terms of the 
base system by 
\begin{align}
  \left( \kappa^{+}_n \right)^2 & =
  -\kappa_{n+1}\kappa_{n}\frac{\phi_{n}(\alpha)}{\phi_{n+1}(\alpha)} ,
  \label{K+1:a} \\
  r^{+}_n & =
  \frac{\phi_{n}(0)\phi_{n+1}(\alpha)-\phi_{n+1}(0)\phi_{n}(\alpha)}
       {\alpha\kappa_{n+1}\phi_{n}(\alpha)},
  \label{K+1:b} \\
  \bar{r}^{+}_n & =
  \frac{\phi^{*}_{n}(\alpha)}{\phi_{n}(\alpha)} ,
  \label{K+1:c}
\end{align}
and the corresponding bi-orthogonal polynomials and associated functions are
\begin{align}
  \Phi_{n}(z) & =
  \frac{\kappa^{+}_n}{\kappa_n}\left[
  \phi_{n}(z)-\frac{\phi_{n+1}(0)}{\kappa_{n+1}\phi_{n}(\alpha)}
              \frac{\phi^*_{n}(\alpha)\phi_{n}(z)-\phi_{n}(\alpha)\phi^*_{n}(z)}{z-\alpha}
                               \right] ,
  \label{K+1:d} \\
  \Phi^*_{n}(z) & =
  \frac{\kappa^{+}_n}{\kappa_n}\left[
  \frac{\phi^*_{n}(\alpha)}{\phi_{n}(\alpha)}\phi_{n}(z)+\frac{\alpha}{\phi_{n}(\alpha)}
  \frac{\phi^*_{n}(\alpha)\phi_{n}(z)-\phi_{n}(\alpha)\phi^*_{n}(z)}{z-\alpha}
                               \right] ,
  \label{K+1:e} \\
   \Xi_{n}(z) & =
  \frac{\kappa^{+}_n}{\kappa_n}\left[
  (z-\alpha)\xi_{n}(z)-\frac{\phi_{n+1}(0)}{\kappa_{n+1}\phi_{n}(\alpha)}
              \left( \phi^*_{n}(\alpha)\xi_{n}(z)+\phi_{n}(\alpha)\xi^*_{n}(z)
              \right) \right] ,
  \label{K+1:f} \\
   \Xi^*_{n}(z) & =
  \frac{\kappa^{+}_n}{\kappa_n}\left[
  -\frac{\phi^*_{n}(\alpha)}{\phi_{n}(\alpha)}(z-\alpha)\xi_{n}(z)
  -\frac{\alpha}{\phi_{n}(\alpha)}
              \left( \phi^*_{n}(\alpha)\xi_{n}(z)+\phi_{n}(\alpha)\xi^*_{n}(z)
              \right) \right] .
  \label{K+1:g}
\end{align}
Alternatively the matrix generator $ R^{+}_{n} $ is given by
\begin{equation}
    R^{+}_n(z;\alpha) = \frac{\kappa^+_n}{\kappa_n}\frac{1}{z-\alpha}
       \begin{pmatrix}
             z-\frac{\DySt\kappa_{n}\phi_{n+1}(\alpha)}{\DySt\kappa_{n+1}\phi_{n}(\alpha)}
	& \frac{\DySt\phi_{n+1}(0)}{\DySt\kappa_{n+1}}	\cr
             \frac{\DySt\phi^*_{n}(\alpha)}{\DySt\phi_{n}(\alpha)}z
	& -\alpha					\cr
       \end{pmatrix} .
\label{K+1_R}
\end{equation}
\end{corollary}
\begin{proof}
This follows easily from Proposition \ref{Kmod} and the expansion formulae 
(\ref{phiexp:a})-(\ref{epsexp:d}).
\end{proof}

\begin{remark}
This case has been given by Zhedanov \cite{Zh_1998} in the context of BOLPS,
which he termed as a Christoffel transformation,
so we can directly compare some of the above results with his.
Accordingly we verify that Eq. (2.6) and (2.7) of \cite{Zh_1998} follow from
(\ref{K+1:a}) and (\ref{K+1:b}), while Eq. (2.1) follows from (\ref{K+1:d}).
\end{remark}

\begin{corollary}\label{Lshift}
For $ K=0,L=1 $ with the denominator factor of $ z-\beta $ denote the coefficients 
of the transformed BOPS by $ \kappa^{-}_n, r^{-}_n, \bar{r}^{-}_n $. 
Assuming $ \xi^*_{n}(\beta) \neq 0 $ for $ n\in \Z_{\geq 0} $ these are given by 
\begin{align}
  \left( \kappa^{-}_n \right)^2 & =
  -\kappa_n\kappa_{n-1}\frac{\beta\xi^*_{n-1}(\beta)}{\xi^*_{n}(\beta)} ,
  \label{L-1:a} \\
  r^{-}_n & =
  \frac{\beta\xi_{n-1}(\beta)}{\xi^{*}_{n-1}(\beta)} ,
  \label{L-1:b} \\
  \bar{r}^{-}_n & =
  \frac{\bar{\phi}_{n}(0)\beta\xi^*_{n-1}(\beta)-\bar{\phi}_{n-1}(0)\xi^*_{n}(\beta)}
       {\beta\kappa_{n}\xi^*_{n-1}(\beta)},
  \label{L-1:c}
\end{align}
and the corresponding bi-orthogonal polynomials and associated functions are
\begin{align}
  \Phi_{n}(z) & =
  -\frac{\kappa_n}{\kappa^{-}_n}\frac{\beta}{\xi^*_{n}(\beta)}
   \left[ \xi^*_{n}(\beta)\phi_{n}(z)+\xi_{n}(\beta)\phi^*_{n}(z)\right] ,
  \label{L-1:d} \\
  \Phi^*_{n}(z) & =
  -\frac{\kappa_n}{\kappa^{-}_n}\left[
  (z-\beta)\left(\frac{\bar{\phi}_n(0)}{\kappa_n}\phi_{n}(z)-\phi^*_{n}(z)\right)
  +\frac{\beta\bar{\phi}_n(0)}{\kappa_n\xi^*_{n}(\beta)}
   \left( \xi^*_{n}(\beta)\phi_{n}(z)+\xi_{n}(\beta)\phi^*_{n}(z)\right)
                               \right] ,
  \label{L-1:e} \\
   \Xi_{n}(z) & =
  -\frac{\kappa_n}{\kappa^{-}_n}\frac{\beta}{\xi^*_{n}(\beta)}
   \frac{\xi^*_{n}(\beta)\xi_{n}(z)-\xi_{n}(\beta)\xi^*_{n}(z)}{z-\beta} ,
  \label{L-1:f} \\
   \Xi^*_{n}(z) & =
  \frac{\kappa_n}{\kappa^{-}_n}\left[
   \xi^*_{n}(z)+\frac{\bar{\phi}_n(0)}{\kappa_n}\xi_n(z)
  +\frac{\bar{\phi}_n(0)}{\kappa_n}\frac{\beta}{\xi^*_{n}(\beta)}
  \frac{\xi^*_{n}(\beta)\xi_{n}(z)-\xi_{n}(\beta)\xi^*_{n}(z)}{z-\beta}
                               \right] .
  \label{L-1:g}
\end{align}
Alternatively the matrix generator $ R^{-}_n $ is given by
\begin{equation}
    R^{-}_n(z;\beta) = -\frac{\kappa_n}{\kappa^-_n}
       \begin{pmatrix}
          \beta
        & \beta\frac{\DySt \xi_{n}(\beta)}{\DySt \xi^*_{n}(\beta)}  \cr
	    \frac{\DySt\bar{\phi}_{n}(0)}{\DySt\kappa_{n}}z
   &       \beta\frac{\DySt\kappa_{n-1}\xi^*_{n-1}(\beta)}{\DySt\kappa_{n}\xi^*_{n}(\beta)}-z  \cr
       \end{pmatrix} .
\label{L-1_R}
\end{equation}
\end{corollary}
\begin{proof}
This also follows from Proposition \ref{Lmod} and the expansion formulae 
(\ref{phiexp:a})-(\ref{epsexp:d}).
\end{proof}

\begin{remark}
This case is specialisation of the Geronimus transformation given by Zhedanov \cite{Zh_1998} 
in the context of BOLPS. These transformations include a contribution from a mass point
which we do not consider here, and we effect this by the choice of the variable (in his notation)
$ \mu = \beta $, $ \phi_n = \xi_n(\beta)/\kappa_n $. Consequently we find
$ V_n = \beta \kappa_n\xi_{n-1}(\beta)/\phi_n(0)\xi^{*}_{n-1}(\beta) $ and therefore
we verify that Eq. (2.15) and (2.16) of \cite{Zh_1998} follow from
(\ref{L-1:a}) and (\ref{L-1:b}), while Eq. (2.11) follows from (\ref{L-1:d}).
\end{remark}

Compatibility of the Christoffel-Geronimus-Uvarov transformations (\ref{Y_CUxfm:a},\ref{Y_CUxfm:b})
with the other fundamental structures of the bi-orthogonal system imply a number 
of simple relations. For example compatibility with the $ n \mapsto n+1 $ recurrence
(\ref{Yrecur:a}) implies the following result.
\begin{proposition}
The $ K_n $ and $ R_n $ matrices satisfy the compatibility relation
\begin{equation}
    R^{\pm}_{n+1}(z)K_n(z) = K^{\pm}_n(z)R^{\pm}_n(z) ,
\end{equation}
for all $ z, n, \alpha, \beta $.
\end{proposition}
\begin{proof}
Considering the $ K=1,L=0 $ case first we use the explicit forms for the 
matrices $ K_n, K^{+}_n $ given by (\ref{Yrecur:a}) along with 
(\ref{K+1:a},\ref{K+1:b}), and the matrices $ R^{+}_n, R^{+}_{n+1} $ as
given in (\ref{K+1_R}). A direct calculation of the difference of the 
left-hand and right-hand sides then yields the null matrix upon repeated use of 
the recurrence relations for $ \phi_{n+2}(\alpha), \phi_{n+1}(\alpha), \phi^*_{n+1}(\alpha) $ 
implied by (\ref{Yrecur:a}).
For the $ K=0, L=1 $ case we used (\ref{Yrecur:a}) and (\ref{L-1:a},\ref{L-1:b}) 
for the matrix $ K^{-}_n $ and the formula (\ref{L-1_R}) for $ R^{-}_n, R^{-}_{n+1} $.
The appropriate difference was shown to be identically zero using the forward 
recurrence for $ \xi^*_{n+1}(\beta) $ as implied by (\ref{Yrecur:a}) 
and the backward recurrence for $ \xi^*_{n-1}(\beta) $ as implied by
(\ref{Yrecur:b}).
\end{proof}

The other structure that we have to verify compatibility with the Christoffel-Geronimus-Uvarov transformations
is the spectral derivative, namely (\ref{YzDer}). In order to verify this we have
to find expressions for the transformation of the spectral matrices $ A^{\pm}_n(z) $
under the action of the generators $ K=1, L=0 $ and $ K=0, L=1 $,
which are given in the following lemma.

\begin{lemma}
With the assumptions of Proposition \ref{Kshift}
the transformed spectral matrix $ A^+_n(z) $ under the Christoffel-Geronimus-Uvarov transformation $ K=1, L=0 $
is given by
\begin{multline}
   (z-\alpha)A^+_n(z)
  = \frac{1}{\kappa_n\phi_{n+1}(\alpha)}
       \begin{pmatrix}
          -\phi_{n+1}(0)\phi^*_{n}(\alpha)
	&  \phi_{n+1}(0)\phi_{n}(\alpha)	\cr
           \alpha\kappa_{n+1}\phi^*_{n}(\alpha)
	& -\alpha\kappa_{n+1}\phi_{n}(\alpha)	\cr
       \end{pmatrix}
   \\
   +\frac{1}{\kappa_n\kappa_{n+1}\phi_n(\alpha)\phi_{n+1}(\alpha)}
       \begin{pmatrix}
           \kappa_n\phi_{n+1}(\alpha)-\kappa_{n+1}\phi_n(\alpha)z
	& -\phi_{n+1}(0)\phi_{n}(\alpha)	\cr
          -\kappa_{n+1}\phi^*_{n}(\alpha)z
	&  \alpha\kappa_{n+1}\phi_{n}(\alpha)	\cr
       \end{pmatrix} A_n(z)
   \\ \times
       \begin{pmatrix}
          -\alpha\kappa_{n+1}\phi_n(\alpha)
	& -\phi_{n+1}(0)\phi_{n}(\alpha)	\cr
          -\kappa_{n+1}\phi^*_{n}(\alpha)z
	&  \kappa_{n+1}\phi_{n}(\alpha)z-\kappa_n\phi_{n+1}(\alpha)	\cr
       \end{pmatrix} ,
\label{K+1_Axfm}
\end{multline}
or
\begin{equation}
  A^+_n(z)
  = -\frac{\kappa^+_n}{\kappa_n}\frac{1}{z-\alpha}\left.{\rm Res}R^+_n\right|_{z=\alpha}
     +R^+_n(z)A_n(z)[R^+_n(z)]^{-1} ,
\label{}
\end{equation}
and the analogous result for $ A^-_n(z) $ under the Christoffel-Geronimus-Uvarov transformation $ K=0, L=1 $
assuming the conditions of Proposition \ref{Lshift} is given by
\begin{multline}
   (z-\beta)A^-_n(z)
   = \frac{1}{\kappa_n}
       \begin{pmatrix}
          0 & 0		\cr
          -\bar{\phi}_n(0)
	&  \kappa_n	\cr
       \end{pmatrix}
   \\
   +\frac{1}{\beta\kappa_{n-1}\kappa_n\xi^*_{n-1}(\beta)\xi^*_n(\beta)}
       \begin{pmatrix}
           \beta\kappa_n\xi^*_n(\beta)
	&  \beta\kappa_n\xi_{n}(\beta)	\cr
           \bar{\phi}_n(0)\xi^*_{n}(\beta)z
	& -\kappa_n\xi^*_{n}(\beta)z+\beta\kappa_{n-1}\xi^*_{n-1}(\beta)	\cr
       \end{pmatrix} A_n(z)
    \\ \times
       \begin{pmatrix}
           \kappa_n\xi^*_n(\beta)z-\beta\kappa_{n-1}\xi^*_{n-1}(\beta)
	&  \beta\kappa_n\xi_{n}(\beta)	\cr
           \bar{\phi}_n(0)\xi^*_{n}(\beta)z
	& -\beta\kappa_n\xi^*_{n}(\beta)	\cr
       \end{pmatrix} ,
\label{L-1_Axfm}
\end{multline} 
or
\begin{equation}
  A^-_n(z)
  = \frac{\kappa^-_n}{\kappa_n}\frac{1}{z-\beta}\left.\frac{d}{dz}R^-_n\right|_{z=\beta}
    +R^-_n(z)A_n(z)[R^-_n(z)]^{-1} .
\label{}
\end{equation}
\end{lemma}
\begin{proof}
To establish the first result (\ref{K+1_Axfm}) we differentiate the expressions 
(\ref{K+1:a},\ref{K+1:b}) with respect to $ z $ and employ the matrix equation
(\ref{YzDer}) for the derivatives of $ \phi_n, \phi^*_n $. Having rendered
the two derivatives in terms of $ \phi_n $ and $ \phi^*_n $ we invert 
(\ref{K+1:a}) and (\ref{K+1:b}) so that these two polynomials are given in
terms of $ \Phi_n $ and $ \Phi^*_n $. The second result is found by an identical
argument starting with (\ref{L-1:a}) and (\ref{L-1:b}) instead.
\end{proof}

Now compatibility of the Christoffel-Geronimus-Uvarov transformations with the spectral derivative must
imply the following result.
\begin{proposition}
Compatibility of the Christoffel-Geronimus-Uvarov transformations $ K=1, L=0 $ and $ K=0, L=1 $ with
the spectral derivative imply the condition
\begin{equation}
   \frac{d}{dz}R^{\pm}_n(z)+R^{\pm}_n(z)A_n(z) = A^{\pm}_n(z)R^{\pm}_n(z) ,
\label{Uxfm+zDer}
\end{equation}
for all $ n,z, \alpha, \beta $.
\end{proposition}
\begin{proof}
We shall verify (\ref{Uxfm+zDer}) by treating the $ \pm $ cases separately. For
the $ + $ case we employ (\ref{K+1_R}), which we also differentiate with respect to 
$ z $, and the formula for $ A^{+}_n $ (\ref{K+1_Axfm}). Then it is a matter of an
exercise to check that the difference between the left-hand and right-hand sides
are identically zero without making any assumptions on $ A_n $. The $ - $ case
is treated in a similar manner starting with (\ref{L-1_R}) and (\ref{L-1_Axfm}). 
\end{proof}

In concluding this section we shall display the matrix generators for modifications
with $ K^*=1 $ or $ L^*=1 $ for the sake of completeness even though these can be 
found from composite transformations with $ K=1, L=1 $ and $ \beta=0 $, $ \alpha=0 $
respectively. These results also exhibit an elegant symmetry with respect to the 
other generators.

\begin{corollary}\label{Kshift*}
For $ K^*=1,L^*=0 $ with the numerator factor of $ 1-\alpha^* z^{-1} $ denote the coefficients 
of the transformed BOPS by $ \kappa^{+}_n, r^{+}_n, \bar{r}^{+}_n $ and the 
bi-orthogonal polynomials and associated functions by 
$ \Phi_{n}, \Phi^*_{n},  \Xi_{n},  \Xi^*_{n} $. 
Assuming $ \phi^*_{n}(\alpha^*) \neq 0 $ for $ n\in \Z_{\geq 0} $ these are given in terms of the
base system by 
\begin{align}
  \left( \kappa^{+}_n \right)^2 & =
  \kappa_{n+1}\kappa_{n}\frac{\phi^*_{n}(\alpha^*)}{\phi^*_{n+1}(\alpha^*)} ,
  \label{K+1*:a} \\
  r^{+}_n & =
  \frac{\phi_{n}(\alpha^*)}{\phi^*_{n}(\alpha^*)} ,
  \label{K+1*:b} \\
  \bar{r}^{+}_n & =
  \frac{\bar{\phi}_{n}(0)\phi^*_{n+1}(\alpha^*)-\bar{\phi}_{n+1}(0)\alpha^*\phi^*_{n}(\alpha^*)}
       {\kappa_{n+1}\phi^*_{n}(\alpha^*)},
  \label{K+1*:c}
\end{align}
and the corresponding matrix generator is given by
\begin{equation}
    R^{+}_n(z;\alpha^*) = \frac{\kappa^+_n}{\kappa_n}\frac{1}{z-\alpha^*}
       \begin{pmatrix}
           z
        & -\frac{\DySt\alpha^*\phi_{n}(\alpha^*)}{\DySt\phi^*_{n}(\alpha^*)} \cr
	        -\frac{\DySt\bar{\phi}_{n+1}(0)\alpha^*}{\DySt\kappa_{n+1}}z
        &  \frac{\DySt\kappa_{n}\phi^*_{n+1}(\alpha^*)}{\DySt\kappa_{n+1}\phi^*_{n}(\alpha^*)}z-\alpha^* \cr
       \end{pmatrix} .
\label{K+1*_R}
\end{equation}
\end{corollary}
\begin{proof}
This follows from Proposition \ref{Kmod*} and the expansion formulae 
(\ref{phiexp:a})-(\ref{epsexp:d}).
\end{proof}

\begin{corollary}\label{Lshift*}
For $ K^*=0,L^*=1 $ with the denominator factor of $ 1-\beta^*z^{-1} $ denote the coefficients 
of the transformed BOPS by $ \kappa^{-}_n, r^{-}_n, \bar{r}^{-}_n $. 
Assuming $ \xi_{n}(\beta^*) \neq 0 $ for $ n\in \Z_{\geq 0} $ these are given by 
\begin{align}
  \left( \kappa^{-}_n \right)^2 & =
  \kappa_n\kappa_{n-1}\frac{\beta^*\xi_{n-1}(\beta^*)}{\xi_{n}(\beta^*)} ,
  \label{L-1*:a} \\
  r^{-}_n & =
  \frac{\phi_{n}(0)\xi_{n-1}(\beta^*)-\phi_{n-1}(0)\xi_{n}(\beta^*)}
       {\kappa_{n}\xi_{n-1}(\beta^*)},
  \label{L-1*:b} \\
  \bar{r}^{-}_n & =
  \frac{\xi^*_{n-1}(\beta^*)}{\beta^*\xi_{n-1}(\beta^*)} ,
  \label{L-1*:c}
\end{align}
and the corresponding matrix generator is given by
\begin{equation}
    R^{-}_n(z;\beta^*) = \frac{\kappa^-_n}{\kappa_n}
       \begin{pmatrix}
          1-\frac{\DySt\kappa_{n}\xi_{n}(\beta^*)}{\DySt\kappa_{n-1}\xi_{n-1}(\beta^*)}z^{-1}
        & \frac{\DySt\phi_n(0)\xi_{n}(\beta^*)}{\DySt\kappa_{n-1}\xi_{n-1}(\beta^*)}z^{-1} \cr
          \frac{\DySt\kappa_{n}\xi^*_{n}(\beta^*)}{\DySt\kappa_{n-1}\beta^*\xi_{n-1}(\beta^*)}
        & \frac{\DySt\kappa_{n}\xi_{n}(\beta^*)}{\DySt\kappa_{n-1}\beta^*\xi_{n-1}(\beta^*)} \cr
       \end{pmatrix} .
\label{L-1*_R}
\end{equation}
\end{corollary}
\begin{proof}
This follows from Proposition \ref{Lmod*} and the expansion formulae 
(\ref{phiexp:a})-(\ref{epsexp:d}).
\end{proof}


\section{The Semi-classical Class of Weights}
\setcounter{equation}{0}
From the viewpoint of classical orthogonal polynomials the weight (\ref{M=3_scw}) is a
particular example of the regular semi-classical class (\ref{gJwgt}), characterised
by a special structure of their logarithmic derivatives
\begin{equation}
  \frac{1}{w(z)}\frac{d}{dz}w(z) = \frac{2V(z)}{W(z)}
  = \sum^M_{j=0}\frac{\rho_j}{z-z_j}, \quad \rho_j \in \C.
\label{scwgt2}
\end{equation}
Here $ V(z) $, $ W(z) $ are polynomials with $ {\rm deg}\;V(z) < M+1 $, $ {\rm deg}\;W(z)=M+1 $.
The zeros of $ W(z) $ define $ M $ movable singularities $ \{z_j\}^{M}_{j=1} $ located at
arbitrary positions in the finite complex plane. The residues at the finite singularities, 
the parameters $ \{\rho_j\}^{M}_{j=0} $, are also arbitrary complex constants only
restricted to ensure that the trigonometric moments have meaning. 
In addition the matrix spectral differential equation in $z$ (\ref{YzDer}) has a distinct
set of fixed singularities, due to the bi-orthogonal system and its support on the unit
circle, at $ z_0=0 $ and $ z_{M+1}=\infty $ for generic values 
of the parameters $ \rho_j $. It is possible that one of the zeros of $ W(z) $
coincides with the origin and we cover this situation by denoting this zero
by $ z_0 $ and its exponent by $ \rho_0 $. 

To avoid technical complications we require the following generic conditions 
for the regular semi-classical class -
\begin{enumerate}
 \item
  $ {\rm deg}\;(W) \geq 2 $,
 \item
  $ {\rm deg}\;(V) < {\rm deg}\;(W) $,
 \item
  the $ M+1 $ zeros of $ W(z) $, $ \{z_0, z_1, \ldots ,z_M\} $ are distinct,
 \item
  the residues $ \rho_j = 2V(z_j)/W'(z_j) \notin \Z_{\geq 0} $.
\end{enumerate}
From the expansion of the denominator in terms of elementary symmetric functions
\begin{equation}
  W(z) = \prod^M_{j=0}(z-z_j(t)) = \sum^{M+1}_{l=0} (-1)^le_l[z_i] z^{M+1-l} .
\end{equation}
and (\ref{scwgt2}) one sees our system is characterised by the singularity data 
$ \{z_j,\rho_j\}^{M}_{j=0} $ of the weight, and our task is to deduce the 
consequences for the bi-orthogonal system.
If (\ref{scwgt2}) is viewed as an ordinary differential equation for $ w(z) $
given the data $ \{z_j,\rho_j\}^{M}_{j=0} $ then one particular solution is the 
generalised Jacobi weight 
\begin{equation}
  w(z) = \prod^M_{j=0}(z-z_j)^{\rho_j} , \quad \rho_j \in \C .
\label{gJwgt}
\end{equation}
We will find that a certain symmetry will be present if we consider the following
"conjugate" expression to the above weight
\begin{equation}
  w(z) = \prod^M_{j=1}(1-z_j z^{-1})^{\rho^*_j}, \quad \rho^*_j \in \C .
\label{gJwgt*}
\end{equation} 

For the regular semi-classical class the spectral coefficients $ \Theta_n(z) $,
$ \Theta^*_n(z) $, $ \Omega_n(z) $ and $ \Omega^*_n(z) $ in (\ref{YzDer}) are
polynomials in $ z $ of fixed degree $ {\rm deg}\;\Omega_n(z)={\rm deg}\;\Omega^*_n(z)=M $,
$ {\rm deg}\;\Theta_n(z)={\rm deg}\;\Theta^*_n(z)=M-1 $, independent of $ n $. This is
true provided that $ W(0)=0 $, i.e. one of the singularities is located at the origin
$ z_0 = 0 $, and is taken only in order to simplify the ensuing expressions.
Once the singularity data is given these polynomials themselves are most usefully 
characterised in two distinct ways - 
either in terms of the coefficients of the monomials in $ z $ or their evaluations
at the singular points 
$ \{\Theta_n(z_j),\Theta^*_n(z_j),\Omega_n(z_j),\Omega^*_n(z_j)\}^{M}_{j=1} $.

As a consequence of the polynomial form of the $ \Theta_n,\Theta^*_n,\Omega_n,\Omega^*_n $ 
the spectral matrix $ A_n $ appearing in (\ref{YzDer}) has
a partial fraction decomposition 
\begin{equation}
   A_n = \sum^{M}_{j=0}\frac{A_{n,j}}{z-z_j} ,
\label{An_PFexp}
\end{equation}
in terms of the residue matrices $ A_{n,j} $.
The $n$-th residue matrix $ A_{n,j} $ at the finite singularity $ z_j $ is given by
\begin{align}
   A_{n,j} & = \frac{\rho_j}{2V(z_j)}
       \begin{pmatrix}
              -\Omega_n(z_j)-V(z_j)
              +\dfrac{\kappa_{n+1}}{\kappa_n}z_j\Theta_n(z_j)
            & \dfrac{\phi_{n+1}(0)}{\kappa_n}\Theta_n(z_j)
            \cr
              -\dfrac{\bar{\phi}_{n+1}(0)}{\kappa_n}z_j\Theta^*_n(z_j)
            &  \Omega^*_n(z_j)-V(z_j)
                     -\dfrac{\kappa_{n+1}}{\kappa_n}\Theta^*_n(z_j)
            \cr
       \end{pmatrix} ,
       \nonumber \\
          & = -\frac{\rho_j}{2z_j^{n}}
       \begin{pmatrix}
              \phi^*_n(z_j)\xi_n(z_j) & -\phi_n(z_j)\xi_n(z_j) \cr
              -\phi^*_n(z_j)\xi^*_n(z_j) & \phi_n(z_j)\xi^*_n(z_j) \cr
       \end{pmatrix} ,
\label{An_resJ}
\end{align}
for $ j=1,\ldots,M $, while for $ j=0 $ the expression is
\begin{equation}
   A_{n,0} = (n-\rho_0)
       \begin{pmatrix}
              1 & -r_n \cr 0 & 0 \cr
       \end{pmatrix} ,
\label{An_res0}
\end{equation}
and for the singular point $ z_{M+1}=\infty $ 
\begin{equation}
   A_{n,\infty} = -\sum^M_{j=0}A_{n,j} =
       \begin{pmatrix}
              -n & 0 \cr -(n+\sum^M_{j=0}\rho_j)\bar{r}_n & \sum^M_{j=0}\rho_j \cr
       \end{pmatrix} .
\label{An_resInfty}
\end{equation}
Let us define the diagonal matrices of formal monodromy $ T_j $ for $ j=0,1,\ldots,M,\infty $ by 
the diagonalisation
\begin{equation}
   A_{n,j} = G_{n,j}T_j\left(G_{n,j}\right)^{-1} ,
\label{An_resDiag}
\end{equation}
and we compute these as 
\begin{align}
   T_{j} & = \begin{cases}
       \begin{pmatrix}
              0 & 0 \cr 0 & n-\rho_0 \cr
       \end{pmatrix} , & j=0 \\
       \begin{pmatrix}
              0 & 0 \cr 0 & -\rho_j \cr
       \end{pmatrix} , & j=1,\ldots,M \\
       \begin{pmatrix}
              -n & 0 \cr 0 & \sum^M_{j=0}\rho_j \cr
       \end{pmatrix} , & j=\infty
             \end{cases} ,
        \\
   G_{n,j} & = \begin{cases}
       \begin{pmatrix}
              \phi_n(0) & d_0/\kappa_n \cr \kappa_n & 0 \cr
       \end{pmatrix} , & j=0 \\
       \begin{pmatrix}
              \phi_n(z_j) & d_j\xi_n(z_j) \cr \phi^*_n(z_j) & -d_j\xi^*_n(z_j) \cr
       \end{pmatrix} , & j=1,\ldots,M \\
       \begin{pmatrix}
              \kappa_n & 0 \cr \bar{\phi}_n(0) & d_{\infty}/\kappa_n \cr
       \end{pmatrix} . & j=\infty
             \end{cases} ,
\label{Tn_Gn}
\end{align}
where $ d_0,d_1,\ldots,d_M,d_{\infty} $ are fixed by the weight and its support. 
The reader should note that our conventions are slightly different from those of 
the Kyoto School in that $ \det A_{n,j} $ vanish for $ j=0,\ldots, M $ and 
$ {\rm Tr} A_{n,j} $ are non-zero in the generic case for $ j=0,\ldots, M,\infty $.

There are sum identities implied by the above relations (\ref{An_resJ}) and 
(\ref{An_resInfty})
\begin{align}
  \sum^M_{k=1}\rho_k \frac{\DySt \phi_n(z_k)\xi_n(z_k)}{2z^n_k} = (n-\rho_0)r_n ,
  \label{An_sum:a} \\
  \sum^M_{k=1}\rho_k \frac{\DySt \phi^*_n(z_k)\xi_n(z_k)}{2z^n_k} = -\rho_0 ,
  \label{An_sum:b} \\
  \sum^M_{k=1}\rho_k \frac{\DySt \phi_n(z_k)\xi^*_n(z_k)}{2z^n_k} = \sum^M_{k=0}\rho_k ,
  \label{An_sum:c} \\
  \sum^M_{k=1}\rho_k \frac{\DySt \phi^*_n(z_k)\xi^*_n(z_k)}{2z^n_k} = \left( n+\sum^M_{k=0}\rho_k\right)\bar{r}_n .
  \label{An_sum:d}
\end{align}

What is clear from our experience is that the evaluations of the polynomials
$ \Theta_n,\Theta^*_n,\Omega_n,\Omega^*_n $ at the singular points
and their subsequent use as auxiliary variables greatly simplifies the theory. 
We shall therefore seek to cast our theory in terms of these quantities.
Firstly we note that the
bilinear products of a polynomial and an associated function evaluated at the finite
singular point $ z_j $ are related to the spectral coefficient evaluations at the
same point by the following equations, valid for all $ j=1,\ldots,M $
\begin{align}
   \phi_n(z_j)\xi_n(z_j)
  & =  2\frac{\phi_{n+1}(0)}{\kappa_n}z^n_j\frac{\Theta_n(z_j)}{2V(z_j)} ,
  \label{BilRes:a} \\
   \phi^*_n(z_j)\xi^*_n(z_j)
  & = -2\frac{\bar{\phi}_{n+1}(0)}{\kappa_n}z^{n+1}_j\frac{\Theta^*_n(z_j)}{2V(z_j)} ,
  \label{BilRes:b}
\end{align}
\begin{align}
   \phi_{n+1}(z_j)\xi_n(z_j)
  & =  2\frac{\phi_{n+1}(0)}{\kappa_n}z^n_j\frac{\Omega_n(z_j)+V(z_j)}{2V(z_j)} ,
  \label{BilRes:c} \\
   \phi_{n}(z_j)\xi_{n+1}(z_j)
  & =  2\frac{\phi_{n+1}(0)}{\kappa_n}z^n_j\frac{\Omega_n(z_j)-V(z_j)}{2V(z_j)} ,
  \label{BilRes:d} \\
   \phi^*_{n}(z_j)\xi^*_{n+1}(z_j)
  & = -2\frac{\bar{\phi}_{n+1}(0)}{\kappa_n}z^{n+1}_j
        \frac{\Omega^*_n(z_j)+V(z_j)}{2V(z_j)} ,
  \label{BilRes:e} \\
   \phi^*_{n+1}(z_j)\xi^*_{n}(z_j)
  & = -2\frac{\bar{\phi}_{n+1}(0)}{\kappa_n}z^{n+1}_j
        \frac{\Omega^*_n(z_j)-V(z_j)}{2V(z_j)} ,
  \label{BilRes:f} \\
   \phi_n(z_j)\xi^*_n(z_j)
  & = -\frac{z^n_j}{V(z_j)}\left[ \Omega_n(z_j)-V(z_j)
        -\frac{\kappa_{n+1}}{\kappa_n}z_j\Theta_n(z_j) \right] ,
  \label{BilRes:g} \\
  & = -\frac{z^n_j}{V(z_j)}\left[ \Omega^*_n(z_j)-V(z_j)
        -\frac{\kappa_{n+1}}{\kappa_n}\Theta^*_n(z_j) \right] ,
  \label{BilRes:h} \\
   \phi^*_n(z_j)\xi_n(z_j)
  & =  \frac{z^n_j}{V(z_j)}\left[ \Omega_n(z_j)+V(z_j)
        -\frac{\kappa_{n+1}}{\kappa_n}z_j\Theta_n(z_j) \right] ,
  \label{BilRes:i} \\
  & =  \frac{z^n_j}{V(z_j)}\left[ \Omega^*_n(z_j)+V(z_j)
        -\frac{\kappa_{n+1}}{\kappa_n}\Theta^*_n(z_j) \right] .
  \label{BilRes:j}
\end{align}

In the semi-classical class the singular points are free parameters in the complex
plane and another key structure in our integrable system are the 
deformation derivatives of the system (\ref{Ydefn})
with respect to arbitrary trajectories of the singularities $ z_j(t), j=1,\ldots,M $
parameterised by the deformation variable $ t $.
The deformation derivatives were derived in \cite{FW_2004a}, and given by
(where $ d/dt \coloneqq \dot{} $)
\begin{equation}
   \frac{d}{dt}Y_{n} = B_n Y_{n}
   = \left\{ B_{\infty} - \sum^{M}_{j=1}\frac{A_{n,j}}{z-z_j}\dot{z}_j
        \right\} Y_{n} ,
\label{YtDer}
\end{equation}
with
\begin{equation}
  B_{\infty} =
  \begin{pmatrix}
    \frac{\DySt\dot{\kappa}_n}{\DySt\kappa_n} & 0 \cr
    \frac{\DySt\dot{\kappa}_n\bar{\phi}_n(0)+\kappa_n\dot{\bar{\phi}}_n(0)}{\DySt\kappa_n} & -\frac{\DySt\dot{\kappa}_n}{\DySt\kappa_n} \cr
  \end{pmatrix}
  = \dot{G}_{n,\infty}G^{-1}_{n,\infty} .
\label{YtDerInfty}
\end{equation}
This latter term arises in our system because the leading term of $ Y_n $ as $ z \to \infty $ 
is not normalised to the identity matrix as is the case in most other treatments.

Compatibility between the set of deformation derivatives and the spectral derivatives
(\ref{YzDer}) leads to the system of integrable non-linear partial differential 
equations for the residue matrices known as the Schlesinger equations, 
\begin{gather}
   \dot{A}_{n,j} = \left[ B_{\infty},A_{n,j} \right]
   + \sum_{k \neq j} \frac{\dot{z}_j-\dot{z}_k}{z_j-z_k} \left[ A_{n,k}, A_{n,j} \right],
   \qquad j=1,\ldots, M ,
   \label{SchlesingerEqn}\\
   \dot{A}_{n,\infty} = \left[ B_{\infty},A_{n,\infty} \right] .
\end{gather}
A number of deformation derivatives will be required subsequently, which were 
computed in \cite{FW_2004a},
\begin{align}
   2\frac{\dot{\kappa}_n}{\kappa_n} & =
   -\sum^M_{k=1}\rho_k\frac{\dot{z}_k}{z_k}\frac{\phi_n(z_k)\xi^*_n(z_k)}{2z_k^n} ,
   \label{kappadot} \\
   \dot{r}_n & =
    \sum^M_{k=1}\rho_k\frac{\dot{z}_k}{z_k}\frac{[\phi_n(z_k)-r_n \phi^*_n(z_k)]\xi_n(z_k)}{2z_k^n} ,
   \label{rdot} \\
   \dot{\bar{r}}_n  & =
    \sum^M_{k=1}\rho_k\frac{\dot{z}_k}{z_k}\frac{[\bar{r}_n\phi_n(z_k)-\phi^*_n(z_k)]\xi^*_n(z_k)}{2z_k^n} .
   \label{rCdot}
\end{align}
In addition we will need the deformation derivatives of the polynomials and
associated functions evaluated at a singular point
$ P_{n,j} := \phi^*_n(z_j)/\phi_n(z_j) $ and $ Q_{n,j} := \xi_n(z_j)/\xi^*_n(z_j) $
\begin{multline}
   \dot{P}_{n,j} = 
  \sum^{M}_{k=1}\rho_k\frac{\dot{z}_k}{z_k}\frac{[\phi_n(z_k)P_{n,j}-\phi^*_n(z_k)]\xi^*_n(z_k)}{2z_k^n} \\
  -\sum_{k\neq j=0}\rho_k\frac{\dot{z}_j-\dot{z}_k}{z_j-z_k}
         \frac{[\phi_n(z_k)P_{n,j}-\phi^*_n(z_k)][\xi_n(z_k)P_{n,j}+\xi^*_n(z_k)]}{2z_k^n} ,
\label{Pdot}
\end{multline}
and
\begin{multline}
   \dot{Q}_{n,j} = 
  -\sum^{M}_{k=1}\rho_k\frac{\dot{z}_k}{z_k}\frac{Q_{n,j}[\phi_n(z_k)+\phi^*_n(z_k)Q_{n,j}]\xi^*_n(z_k)}{2z_k^n} \\
  -\sum_{k\neq j=0}\rho_k\frac{\dot{z}_j-\dot{z}_k}{z_j-z_k}
         \frac{[\phi_n(z_k)+\phi^*_n(z_k)Q_{n,j}][\xi_n(z_k)-\xi^*_n(z_k)Q_{n,j}]}{2z_k^n} ,
\label{Qdot}
\end{multline}
which can easily be deduced from the above work.

A consequence of the evolution of the residue matrices under (\ref{SchlesingerEqn}) 
is the result that arbitrary deformations preserve the monodromy data, in this case 
the collection of monodromy matrices $ \{ M_j\}^{M+1}_{j=0} $, with each defined by 
\begin{equation}
   Y_n(z_j+\delta e^{2\pi i}) = Y_n(z_j+\delta)M_j .
\end{equation}
Each monodromy matrix has an upper triangular form, one of the three structures
expected for classical solutions of the Garnier systems \cite{Mz_2002},
\begin{equation}
   M_j = C^{-1}_j e^{2\pi iT_j}C_j = 
       \begin{pmatrix}
              1 & c_j(1-e^{-2\pi i\rho_j}) \cr
              0 & e^{-2\pi i\rho_j} \cr
       \end{pmatrix} ,\quad j=0,\ldots,M,
\end{equation}
where the $ c_j $, and thus the $ M_j $, are independent of $ z_k $ or $ t $, and also
of $ n $.
The connection matrices are given by
\begin{equation}
   C_j =
       \begin{pmatrix}
              1 & c_j \cr
              0 & 1 \cr
       \end{pmatrix} ,\quad j=0,\ldots,M .
\end{equation}

We note that the BOPS satisfies the necessary condition for being 
a classical solution of a Garnier system \cite{Wi_2007b}. This means that the formal 
monodromy exponents given by 
\begin{align}
   \theta_0 & = n-\rho_0 ,
  \\
   \theta_j & = -\rho_j, \quad j=1,\ldots,M ,
  \\
   \theta_{\infty} & = n+\sum^{M}_{j=0}\rho_j ,
\end{align}
satisfy the classical condition $ \sum_{0,1,\ldots,M,\infty} \theta_j = 2n $.

There is a final structure in our integrable system that needs to be put into place
and this concerns the discrete deformations of the parameter set of the singularity
data, the set of exponents $ \{ \rho_j\}^{M}_{j=0} $. These are known as the 
Schlesinger transformations and are the subject of the next section.   

\section{The Schlesinger Transformations}
\setcounter{equation}{0}
When one of the free parameters $ \alpha,\beta,\alpha^*,\beta^* $ in Propositions 
\ref{Kshift}, \ref{Lshift}, \ref{Kshift*}, or \ref{Lshift*} co-incides with a 
singularity of a regular semi-classical weight, 
say $ z_j $, then the weight (\ref{gJwgt}) undergoes a modification
$ \rho_j \mapsto \rho_j \pm 1 $  or $ \rho^*_j \mapsto \rho^*_j \pm 1 $ respectively. 
In this situation we are in reality
treating a family of weights whose exponents differ by integers and in the 
context of isomonodromic systems these modifications are Schlesinger transformations.
These transformations and the basic recurrence entail the following types of
transformations in terms of monodromy exponents  
\vskip1cm
\begin{center}
\begin{tabular}{|c|c|}
   \hline
   $ n \mapsto n\pm 1 $ & $ \theta_0 \mapsto \theta_0\pm 1 $ , $ \theta_{\infty} \mapsto \theta_{\infty}\pm 1 $ \\
   $ \rho_j \mapsto  \rho_j\pm 1 $ & $ \theta_j \mapsto \theta_j\mp 1 $ , $ \theta_{\infty} \mapsto \theta_{\infty}\pm 1 $ \\
   $ \rho^*_j \mapsto  \rho^*_j\pm 1 $ & $ \theta_{0} \mapsto \theta_{0}\pm 1 $ , $ \theta_j \mapsto \theta_j\mp 1 $ \\
   $ n \mapsto n\pm 1 $ , $ \rho_j \mapsto  \rho_j\mp 1 $ & $ \theta_0 \mapsto \theta_0\pm 1 $ , $ \theta_j \mapsto \theta_{j}\pm 1 $ \\
   $ n \mapsto n\pm 1 $ , $ \rho^*_j \mapsto  \rho^*_j\mp 1 $ & $ \theta_j \mapsto \theta_{j}\pm 1 $ , $ \theta_{\infty} \mapsto \theta_{\infty}\pm 1 $ \\
   \hline
\end{tabular}
\end{center}
\vskip1cm

If we define the matrix transformations for the mappings
$ \rho_j \mapsto \rho_j\pm1 $, $ \rho^*_j \mapsto \rho^*_j \pm 1 $ for $ j=0,\ldots, M $,
\begin{align}
   Y_n(z;\rho_j\pm1) & = R^{j\pm}_n(z;\rho_j)Y_n(z;\rho_j) ,
   \label{YrhoShift}\\
   Y_n(z;\rho^*_j\pm1) & = R^{*j\pm}_n(z;\rho_j)Y_n(z;\rho^*_j) ,
   \label{YrhoShift*} 
\end{align}
(cf. (\ref{Y_CUxfm:a}),(\ref{Y_CUxfm:b}))
then these are given in the following proposition. For this and the following 
propositions it is only necessary to exhibit proofs for the $ R^{j\pm}_n $ cases as
the $ R^{*j\pm}_n $ transformations are composed from the $ R^{j\pm}_n $ and $ R^{0\pm}_n $
transformations.


\begin{proposition}\label{Prop_Rxfm}
The elementary Schlesinger transformations for the shift 
$ \rho_j \mapsto \rho_j\pm1 $, with $ z_j\neq 0 $, are given in terms of the base 
system with parameter exponents $ \{\rho_j\}^{M}_{j=0} $ by the formulae
\begin{align}
 R^{j+}_n(z;\rho_j) & = \frac{1}{z-z_j}\frac{\kappa^+_n}{\kappa_n}
       \begin{pmatrix}
             z-\frac{\DySt\kappa_{n}\phi_{n+1}(z_j)}{\DySt\kappa_{n+1}\phi_{n}(z_j)}
	& \frac{\DySt\phi_{n+1}(0)}{\DySt\kappa_{n+1}}	\cr
             \frac{\DySt\phi^*_{n}(z_j)}{\DySt\phi_{n}(z_j)}z
	& -z_j					\cr
       \end{pmatrix} ,
 \label{R+} \\
 R^{j-}_n(z;\rho_j) & = \frac{\kappa^-_n}{\kappa_n}
       \begin{pmatrix}
             \frac{\DySt\kappa_{n}\xi^*_{n}(z_j)}{\DySt\kappa_{n-1}\xi^*_{n-1}(z_j)}
	& \frac{\DySt\kappa_{n}\xi_{n}(z_j)}{\DySt\kappa_{n-1}\xi^*_{n-1}(z_j)}	\cr
             \frac{\DySt\bar{\phi}_{n}(0)\xi^*_{n}(z_j)}{\DySt\kappa_{n-1}\xi^*_{n-1}(z_j)}\frac{\DySt z}{\DySt z_j}
	& 1-\frac{\DySt\kappa_{n}\xi^*_{n}(z_j)}{\DySt\kappa_{n-1}\xi^*_{n-1}(z_j)}\frac{\DySt z}{\DySt z_j}					\cr
       \end{pmatrix} ,
 \label{R-}
\end{align}
assuming that $ \phi_{n}(z_j) \neq 0 $ and $ \xi^*_{n}(z_j) \neq 0 $ for all 
$ n\in \Z_{\geq 0} $ respectively.
\end{proposition}
\begin{proof}
The results (\ref{R+}) and (\ref{R-}) follow immediately from (\ref{K+1_R}) and 
(\ref{L-1_R}) respectively.
\end{proof}

\begin{remark}
One can verify
\begin{equation}
  \det R^{j+}_n(z;\rho_j) = \frac{1}{z-z_j}, \quad
  \det R^{j-}_n(z;\rho_j) = z-z_j ,
\end{equation}
using the forward recurrence for $ \phi_{n+1}(z_j) $, i.e. (\ref{Yrecur:a}), 
in the first equality and in the case of the second, the forward recurrence
for $ \xi_n(z_j) $ (\ref{Yrecur:a}) and the backward recurrence for 
$ \xi_{n-1}(z_j) $ (\ref{Yrecur:b}) plus (\ref{l}). 
\end{remark}

\begin{corollary}
The transformed leading coefficients are given now by  
\begin{align}
  \left( \kappa^{+}_n \right)^2 
 & = -\kappa_{n+1}\kappa_{n}\frac{\phi_{n}(z_j)}{\phi_{n+1}(z_j)}
   = -\kappa_{n+1}\kappa_{n}\frac{\Theta_{n}(z_j)}{\Omega_{n}(z_j)+V(z_j)} ,
  \label{kappa+} \\
  r^{+}_{n}
 & = \frac{\phi_{n}(0)\phi_{n+1}(z_j)-\phi_{n+1}(0)\phi_{n}(z_j)}{z_j\kappa_{n+1}\phi_{n}(z_j)}
   = \frac{\phi_{n}(0)}{\kappa_{n+1}}
     \frac{\Omega_{n}(z_j)+V(z_j)-\frac{\DySt \phi_{n+1}(0)}{\DySt \phi_{n}(0)}\Theta_{n}(z_j)}{z_j\Theta_{n}(z_j)}
  \label{refl+:a} \\
  \bar{r}^{+}_{n}
 & = \frac{\phi^*_{n}(z_j)}{\phi_{n}(z_j)}
   = \frac{\bar{\phi}_{n+1}(0)}{\kappa_{n}}
     \frac{\Theta^*_{n}(z_j)}  
        {\Omega^*_{n}(z_j)-V(z_j)-\frac{\DySt \kappa_{n+1}}{\DySt \kappa_{n}}\Theta^*_{n}(z_j)}
  \label{refl+:b}  
\end{align}
and
\begin{align}
  \left( \kappa^{-}_n \right)^2
 & = -\kappa_{n}\kappa_{n-1}\frac{z_j\xi^*_{n-1}(z_j)}{\xi^*_{n}(z_j)}
   = -\kappa_{n}\kappa_{n-1}z_j\frac{\Theta^*_{n-1}(z_j)}{\Omega^*_{n-1}(z_j)+V(z_j)} ,
  \label{kappa-} \\
  r^{-}_{n}
 & = \frac{z_j\xi_{n-1}(z_j)}{\xi^*_{n-1}(z_j)}
   = -\frac{\phi_n(0)}{\kappa_{n-1}}\frac{z_j\Theta_{n-1}(z_j)}
      {\Omega_{n-1}(z_j)-V(z_j)-\frac{\DySt \kappa_n}{\DySt \kappa_{n-1}}z_j\Theta_{n-1}(z_j)} ,
  \label{refl-:a} \\
  \bar{r}^{-}_{n}
 & = \frac{\bar{\phi}_{n}(0)z_j\xi^*_{n-1}(z_j)-\bar{\phi}_{n-1}(0)\xi^*_{n}(z_j)}{z_j\kappa_{n}\xi^*_{n-1}(z_j)}
  \label{refl-:b} \\
 & = -\frac{\bar{\phi}_{n-1}(0)}{\kappa_{n}}
  \frac{\Omega^*_{n-1}(z_j)+V(z_j)-\frac{\DySt\bar{\phi}_{n}(0)}{\DySt\bar{\phi}_{n-1}(0)}z_j\Theta^*_{n-1}(z_j)}{z_j\Theta^*_{n-1}(z_j)} .
  \nonumber
\end{align}
\end{corollary}
\begin{proof}
These six formulae follow from the Corollaries \ref{Kshift} and \ref{Lshift} 
and the bilinear product expressions - in particular
(\ref{kappa+}) is a consequence of (\ref{K+1:a}) and (\ref{BilRes:a},\ref{BilRes:c}), 
(\ref{refl+:a}) from (\ref{K+1:b}) and (\ref{BilRes:a},\ref{BilRes:c}),
(\ref{refl+:b}) from (\ref{K+1:c}) and (\ref{BilRes:b},\ref{BilRes:h}),
(\ref{kappa-}) is a consequence of (\ref{L-1:a}) and (\ref{BilRes:b},\ref{BilRes:e}), 
(\ref{refl-:a}) from (\ref{L-1:b}) and (\ref{BilRes:a},\ref{BilRes:g}) and
(\ref{refl-:b}) from (\ref{L-1:c}) and (\ref{BilRes:b},\ref{BilRes:e}).
\end{proof}

\begin{proposition}
Subject to the conditions in Proposition \ref{Prop_Rxfm} 
the inverses of the Schlesinger matrices (\ref{R+}) and (\ref{R-}) are given by
\begin{align}
 \left( R^{j+}_n(z;\rho_j) \right)^{-1} =
 R^{j-}_n(z;\rho_j+1) & = \frac{\phi_n(z_j)}{\kappa^+_n}
       \begin{pmatrix}
          \frac{\DySt\kappa_{n+1}}{\DySt\phi_{n+1}(z_j)}z_j
	& \frac{\DySt\phi_{n+1}(0)}{\DySt\phi_{n+1}(z_j)}	\cr
          \frac{\DySt\kappa_{n+1}\phi^*_{n}(z_j)}{\DySt\phi_{n+1}(z_j)\phi_{n}(z_j)}z
	& \frac{\DySt\kappa_{n}}{\DySt\phi_{n}(z_j)}-\frac{\DySt\kappa_{n+1}}{\DySt\phi_{n+1}(z_j)}z					\cr
       \end{pmatrix} ,
 \label{invR+} \\
 \left( R^{j-}_n(z;\rho_j) \right)^{-1} =
 R^{j+}_n(z;\rho_j-1) & = \frac{1}{z-z_j}\frac{\kappa_n}{\kappa^-_n}
       \begin{pmatrix}
          z-\frac{\DySt\kappa_{n-1}\xi^*_{n-1}(z_j)}{\DySt\kappa_{n}\xi^*_{n}(z_j)}z_j
	& \frac{\DySt \xi_{n}(z_j)}{\DySt \xi^*_{n}(z_j)}z_j	\cr
          \frac{\DySt\bar{\phi}_{n}(0)}{\DySt\kappa_{n}}z
	& -z_j				\cr
       \end{pmatrix} ,
 \label{invR-}
\end{align}
respectively.
\end{proposition}
\begin{proof}
To establish the first equality of (\ref{invR+}) starting from (\ref{R+}) and 
(\ref{R-}) we require $ R^{j-}_n(z;\rho_j+1) $ and this is found from (\ref{R-}) 
under the mapping $ \rho_j \mapsto \rho_j+1 $. This we compute to be the
second equality of (\ref{invR+}) and first then follows easily.
For the first equality of (\ref{invR-}) we now require $ R^{j+}_n(z;\rho_j-1) $ 
which is found from (\ref{R+}) under the mapping $ \rho_j \mapsto \rho_j-1 $.
Our computation of this yields the second equality of (\ref{invR-}) and the
first follows in a straight forward way.
\end{proof}

From the explicit forms of the Schlesinger transformations it is not obvious that
they are commutative, however this must be the case and the following proposition 
verifies this fact.
\begin{proposition}
Any two distinct Schlesinger transformations 
$ \rho_j\mapsto \rho_j+\epsilon $, $ \rho_k\mapsto \rho_k+\epsilon' $, $ k\neq j $
and $ \epsilon=\pm 1 $, $ \epsilon'=\pm 1 $ commute
\begin{equation}
    R^{j\epsilon}_n(\rho_j,\rho_k+\epsilon')R^{k\epsilon'}_n(\rho_j,\rho_k)
  = R^{k\epsilon'}_n(\rho_j+\epsilon,\rho_k)R^{j\epsilon}_n(\rho_j,\rho_k) .
\label{Rcommute}
\end{equation}
\end{proposition}
\begin{proof}
We take the upward shift case ($ \epsilon=\epsilon'=+1 $) first. To check (\ref{Rcommute}) we have to
form $ R^{j+}_n(\rho_j,\rho_k+1) $ and $ R^{k+}_n(\rho_j+1,\rho_k) $ using the
explicit formula (\ref{R+}) which entails replacing all coefficients and polynomials
appearing therein under the mappings $ \rho_k \mapsto \rho_k+1 $ and $ \rho_j \mapsto \rho_j+1 $
respectively. We observe that
\begin{equation} 
   \left( \kappa^{j+,k+}_n \right)^2
  = \kappa_{n}\kappa_{n+2}\frac{\phi_n(z_k)\phi_{n+1}(z_j)-\phi_n(z_j)\phi_{n+1}(z_k)}
                                 {\phi_{n+1}(z_k)\phi_{n+2}(z_j)-\phi_{n+1}(z_j)\phi_{n+2}(z_k)}
  = \left( \kappa^{k+,j+}_n \right)^2 . 
\label{kappa_jkUp}
\end{equation} 
All components of the difference between the left-hand and right-hand sides of 
(\ref{Rcommute}) vanish trivially except the $ (1,1) $ component. To show this 
also vanishes we have to use the three term recurrence relation
\begin{equation} 
   \kappa_{n+1}[\phi_{n+1}(0)\phi_{n+2}(z)-\phi_{n+2}(0)\phi_{n+1}(z)]
  = z[\kappa_{n+2}\phi_{n+1}(0)\phi_{n+1}(z)-\kappa_{n}\phi_{n+2}(0)\phi_{n}(z)] .
\end{equation} 
The other cases of mixed up and down shifted or down shifted parameters can be proved
using the above equality and employing the inverse formulae (\ref{invR+}) and
(\ref{invR-}).
\end{proof}

\begin{corollary}
The up-shifted bi-orthogonal polynomials and associated functions
have evaluations at the singular points given by the formulae
\begin{align}
  \Phi^{j+}_n(z_k) & = \begin{cases}
  \frac{\DySt \kappa^{+}_{n}}{\DySt \kappa_{n+1}\phi_{n}(z_j)}
  \frac{\DySt \phi_{n+1}(z_j)\phi_{n}(z_k)-\phi_{n}(z_j)\phi_{n+1}(z_k)}{\DySt z_j-z_k},
  & k\neq j \\
  \frac{\DySt \kappa^{+}_{n}}{\DySt \kappa_{n}}
  \left\{ \phi_{n}(z_j)
        +\frac{\DySt \phi_{n+1}(0)}{\DySt \kappa_{n+1}\phi_{n}(z_j)}T^{j+}_n(z_j)
  \right\},
  & k=j \end{cases} ,
  \label{BOPSxfm+:a} \\
  \Phi^{*,j+}_n(z_k) & = \begin{cases}
  \frac{\DySt \kappa^{+}_{n}}{\DySt \bar{\phi}_{n+1}(0)\phi_{n}(z_j)}
  \frac{\DySt \phi^*_{n+1}(z_j)\phi^*_{n}(z_k)-\phi^*_{n}(z_j)\phi^*_{n+1}(z_k)}{\DySt z_j-z_k},
  & k\neq j \\
  \frac{\DySt \kappa^{+}_{n}\phi^*_n(z_j)}{\DySt \kappa_{n}\phi_n(z_j)}
  \left\{ \phi_{n}(z_j)
        -\frac{\DySt z_j}{\DySt \phi^*_{n}(z_j)}T^{j+}_n(z_j)
  \right\},
  & k=j \end{cases} ,
  \label{BOPSxfm+:b} \\
   \Xi^{j+}_n(z_k) & = \begin{cases}
  \frac{\DySt \kappa^{+}_{n}}{\DySt \kappa_{n+1}\phi_{n}(z_j)}
  \left[ \phi_{n}(z_j)\xi_{n+1}(z_k)-\phi_{n+1}(z_j)\xi_{n}(z_k) \right],
  & k\neq j \\
  -\frac{\DySt \kappa^{+}_{n}\phi_{n+1}(0)}{\DySt \kappa_{n}\kappa_{n+1}\phi_{n}(z_j)}2z^n_j, & k=j
                      \end{cases} ,
  \label{BOPSxfm+:c} \\
   \Xi^{*,j+}_n(z_k) & = \begin{cases}
  \frac{\DySt \kappa^{+}_{n}}{\DySt \bar{\phi}_{n+1}(0)\phi_{n}(z_j)}
  \left[ \phi^*_{n}(z_j)\xi^*_{n+1}(z_k)-\phi^*_{n+1}(z_j)\xi^*_{n}(z_k) \right],
  & k\neq j \\
  -\frac{\DySt \kappa^{+}_{n}}{\DySt \kappa_{n}\phi_{n}(z_j)}2z^{n+1}_j, & k=j
                      \end{cases} ,
  \label{BOPSxfm+:d}
\end{align}
provided $ \phi_n(z_j) \neq 0 $ for all $ n\in \Z_{\geq 0} $, where
\begin{multline}
   T^{j+}_n(z_j) =
  \frac{1}{\rho_j+1} \sum_{k\neq j}\frac{\rho_k}{2z^n_k}
           \frac{\phi_{n}(z_j)\phi^*_{n}(z_k)-\phi^*_{n}(z_j)\phi_{n}(z_k)}{z_j-z_k}
  \\ \times\left[ \phi_{n}(z_j)\xi^*_{n}(z_k)+\phi^*_{n}(z_j)\xi_{n}(z_k) \right] .
\label{BOPSxfm+:e}
\end{multline}
The corresponding results for the down-shifted system is
\begin{align}
  \Phi^{j-}_n(z_k) & = \begin{cases}
  \frac{\DySt \kappa^{-}_{n}}{\DySt \kappa_{n-1}\xi^{*}_{n-1}(z_j)}
  \left[ \phi_{n}(z_k)\xi^*_{n}(z_j)+\phi^*_{n}(z_k)\xi_{n}(z_j)\right],
  & k\neq j \\
  \frac{\DySt \kappa^{-}_{n}}{\DySt \kappa_{n-1}\xi^{*}_{n-1}(z_j)}2z^n_j, & k=j
                      \end{cases} ,
  \label{BOPSxfm-:a} \\
  \Phi^{*,j-}_n(z_k) & = \begin{cases}
  \frac{\DySt \kappa^{-}_{n}}{\DySt \kappa_{n}z_j\xi^{*}_{n-1}(z_j)}
  \left[ z_j\xi^*_{n-1}(z_j)\phi^*_{n}(z_k)-z_k\xi^*_{n}(z_j)\phi^*_{n-1}(z_k)\right],
  & k\neq j \\
  \frac{\DySt \kappa^{-}_{n}\bar{\phi}_n(0)}{\DySt \kappa_{n}\kappa_{n-1}\xi^{*}_{n-1}(z_j)}
   2z^n_j, & k=j \end{cases} ,
  \label{BOPSxfm-:b} \\
   \Xi^{j-}_n(z_k) & = \begin{cases}
  \frac{\DySt \kappa^{-}_{n}}{\DySt \kappa_{n-1}\xi^{*}_{n-1}(z_j)}
  \frac{\DySt \xi_{n}(z_j)\xi^*_{n}(z_k)-\xi^*_{n}(z_j)\xi_{n}(z_k)}
       {\DySt z_j-z_k}, & k\neq j \\
  -\frac{\DySt \kappa^{-}_{n}}{\DySt \kappa_{n-1}\xi^*_{n-1}(z_j)}T^{j-}_n(z_j),
  & k=j \end{cases} ,
  \label{BOPSxfm-:c} \\
   \Xi^{*,j-}_n(z_k) & = \begin{cases}
  -\frac{\DySt \kappa^{-}_{n}}{\DySt \kappa_{n}z_j\xi^{*}_{n-1}(z_j)}
  \frac{\DySt z_j\xi^*_{n-1}(z_j)\xi^*_{n}(z_k)-z_k\xi^*_{n}(z_j)\xi^*_{n-1}(z_k)}
       {\DySt z_j-z_k},
  & k\neq j \\
  \frac{\DySt \kappa^{-}_{n}}{\DySt \kappa_{n}}
  \left\{ 
  \frac{\DySt \bar{\phi}_{n}(0)}{\DySt \kappa_{n-1}\xi^*_{n-1}(z_j)}T^{j-}_n(z_j)
        -\frac{\DySt \xi^*_{n}(z_j)}{\DySt z_j}
  \right\},
  & k=j \end{cases} ,
  \label{BOPSxfm-:d}
\end{align}
provided $ \xi^{*}_{n}(z_j) \neq 0 $ for all $ n\in \Z_{\geq 0} $, where
\begin{multline}
   T^{j-}_n(z_j) =
  \frac{1}{\rho_j-1} \sum_{k\neq j}\frac{\rho_k}{2z^n_k}
           \frac{\xi_{n}(z_j)\xi^*_{n}(z_k)-\xi^*_{n}(z_j)\xi_{n}(z_k)}{z_j-z_k}
  \\ \times\left[ \xi^*_{n}(z_j)\phi_{n}(z_k)+\xi_{n}(z_j)\phi^*_{n}(z_k) \right] .
\label{BOPSxfm-:e}
\end{multline}
\end{corollary}

Compatibility constraints need to be satisfied by the over-determined linear
bi-orthogonal system and we have to consider three such constraints. The first 
is the compatibility of Schlesinger transformations (\ref{YrhoShift}) with the
recurrence relation (\ref{Yrecur:a}) and this implies the following
proposition.
\begin{proposition}
The recurrence matrix and the Schlesinger transformation matrices of the bi-orthogonal 
system satisfy the compatibility condition
\begin{equation} 
   R^{j\pm}_{n+1}(z;\rho_j)K_n(z;\rho_j) = K_n(z;\rho_j\pm1)R^{j\pm}_{n}(z;\rho_j) .
\label{Recur+Schlesinger}
\end{equation} 
\end{proposition}
\begin{proof}
Taking the $+$ case first we construct $ K_n(z;\rho_j+1) $ using the specialisation
$ \alpha \to z_j $ of (\ref{K+1:b},\ref{K+1:c}). Forming the difference between
the left-hand and right-hand sides of (\ref{Recur+Schlesinger}) we find the 
$ (1,2) $ and $ (2,2) $ elements vanish trivially. Then employing the forward
recurrence for $ \phi^*_{n+1}(z_j) $ we find the $ (2,1) $ element also vanishes. 
Utilising the forward recurrences for $ \phi_{n+1}(z_j) $ and $ \phi_{n+2}(z_j) $
and finally the previous recurrence we see that the $ (1,1) $ element is also
identically zero. In the $-$ case we compute $ K_n(z;\rho_j-1) $ using the
specialisation $ \beta \to z_j $ in (\ref{L-1:b},\ref{L-1:c}). Forming the
difference again we find $ (1,1) $ element vanishes when the backward recurrence
(\ref{Yrecur:a}) is applied to $ \xi^*_{n-1}(z_j) $ and $ \xi^*_n(z_j) $.
The $ (2,1) $ element vanished after application of the forward recurrence for
$ \xi^*_{n+1}(z_j) $ and the backward recurrence for $ \xi^*_{n-1}(z_j) $.
The $ (1,2) $ element vanishes after the backward recurrence for 
$ \xi_n(z_j) $ is utilised. Finally the $ (2,2) $ element was shown to be
zero upon the application of (\ref{l}), the forward recurrence of 
$ \xi^*_{n+1}(z_j) $ and the backward recurrence for $ \xi^*_{n-1}(z_j) $.
\end{proof}

The second constraint is the compatibility of the Schlesinger transformations
with the spectral derivative (\ref{YzDer}). Again we need a preliminary
result giving explicit formulae for the transformed residue spectral matrices 
$ A^{j\pm}_{n,k} = A_{n,k}(\rho_j\pm1) $, which is contained in the following
lemma.
\begin{lemma}
The transformed residue of the spectral matrix under the Schlesinger transformation
$ \rho_j \mapsto \rho_j\pm 1 $ is given by
\begin{equation} 
\label{}
   A^{j\pm}_{n,k} = R^{j\pm}_n(z_k) A_{n,k} \left( R^{j\pm}_n(z_k) \right)^{-1} ,
\end{equation}
for $ k\neq j $ or explicitly by
\begin{multline}
   A^{j+}_{n,k}
   = \frac{1}{\kappa_n\kappa_{n+1}\phi_n(z_j)\phi_{n+1}(z_j)}
       \begin{pmatrix}
           \kappa_{n+1}\phi_n(z_j)+\frac{\DySt \phi_{n+1}(0)\phi^*_{n}(z_j)}{\DySt z_j-z_k}
	& -\frac{\DySt \phi_{n+1}(0)\phi_{n}(z_j)}{\DySt z_j-z_k}	\cr
          -z_k\frac{\DySt \kappa_{n+1}\phi^*_{n}(z_j)}{\DySt z_j-z_k}
	&  z_j\frac{\DySt \kappa_{n+1}\phi_{n}(z_j)}{\DySt z_j-z_k}	\cr
       \end{pmatrix} A_{n,k}
   \\ \times
       \begin{pmatrix}
           z_j\kappa_{n+1}\phi_n(z_j)
	&  \phi_{n+1}(0)\phi_{n}(z_j)	\cr
           z_k\kappa_{n+1}\phi^*_{n}(z_j)
	&  (z_j-z_k)\kappa_{n+1}\phi_{n}(z_j)+\phi_{n+1}(0)\phi^*_n(z_j)	\cr
       \end{pmatrix} ,
\label{Axfm+:a}
\end{multline}
and
\begin{multline}
   A^{j-}_{n,k}
   = \frac{1}{\kappa_{n-1}\kappa_{n}z_j\xi^*_{n-1}(z_j)\xi^*_{n}(z_j)}
       \begin{pmatrix}
           z_j\frac{\DySt \kappa_{n}\xi^*_{n}(z_j)}{\DySt z_j-z_k} 
	&  z_j\frac{\DySt \kappa_{n}\xi_{n}(z_j)}{\DySt z_j-z_k}	\cr
	   z_k\frac{\DySt \bar{\phi}_{n}(0)\xi^*_{n}(z_j)}{\DySt z_j-z_k}
	& \kappa_{n}\xi^*_n(z_j)+z_j\frac{\DySt \bar{\phi}_{n}(0)\xi_{n}(z_j)}{\DySt z_j-z_k} \cr
       \end{pmatrix} A_{n,k}
   \\ \times
       \begin{pmatrix}
           (z_j-z_k)\kappa_{n}\xi^*_{n}(z_j)+z_j\bar{\phi}_{n}(0)\xi_n(z_j)	
	& -z_j\kappa_{n}\xi_{n}(z_j)	\cr
          -z_k\bar{\phi}_{n}(0)\xi^*_n(z_j)
	&  z_j\kappa_{n}\xi^*_n(z_j)	\cr
       \end{pmatrix} .
\label{Axfm-:a}
\end{multline}
In the special case $ k=j $ we have
\begin{multline}
   A^{j+}_{n,j}
   = \frac{1}{\kappa_n\phi_{n+1}(z_j)}
       \begin{pmatrix}
          -\phi_{n+1}(0)\phi^*_{n}(z_j)
	&  \phi_{n+1}(0)\phi_{n}(z_j)	\cr
           z_j\kappa_{n+1}\phi^*_{n}(z_j)
	& -z_j\kappa_{n+1}\phi_{n}(z_j)	\cr
       \end{pmatrix}
   \\ \times
   \Bigg\{ (\rho_j+1)I
           + \frac{1}{\kappa_{n+1}\phi_n(z_j)}\sum_{k\neq j}\frac{A_{n,k}}{z_j-z_k}
       \begin{pmatrix}
           z_j\kappa_{n+1}\phi_n(z_j)
	&  \phi_{n+1}(0)\phi_{n}(z_j)	\cr
           z_j\kappa_{n+1}\phi^*_{n}(z_j)
	&  \phi_{n+1}(0)\phi^*_{n}(z_j)	\cr
       \end{pmatrix} \Bigg\} .
\label{Axfm+:b}
\end{multline}
and
\begin{multline}
   A^{j-}_{n,j}
   = \frac{1}{\kappa_n}
       \begin{pmatrix}
         0 & 0 \cr	-\bar{\phi}_{n}(0) & \kappa_{n}	\cr
       \end{pmatrix} (\rho_j-1)
   \\
     + \frac{z^2_j}{\kappa_{n-1}\kappa_{n}\xi^*_{n-1}(z_j)\xi^*_{n}(z_j)}
       \begin{pmatrix}
           \kappa_{n}\xi^*_n(z_j)
	&  \kappa_{n}\xi_{n}(z_j)	\cr
           \bar{\phi}_{n}(0)\xi^*_{n}(z_j)	
	&  \bar{\phi}_{n}(0)\xi_{n}(z_j)	\cr
       \end{pmatrix} 
     \sum_{k\neq j}\frac{A_{n,k}}{z_j-z_k}
       \begin{pmatrix}
          -\bar{\phi}_{n}(0)\xi_{n}(z_j)	
	&  \kappa_{n}\xi_{n}(z_j)	\cr
           \bar{\phi}_{n}(0)\xi^*_{n}(z_j)	
	& -\kappa_{n}\xi^*_n(z_j)	\cr
       \end{pmatrix} .
\label{Axfm-:b}
\end{multline}
Again we require that $ \phi_n(z_j) \neq 0 $ for all $ n\in \Z_{\geq 0 } $ for
the upward transformation and $ \xi^*_n(z_j) \neq 0 $ for all
$ n\in \Z_{\geq 0 } $ for the downward shift.
\end{lemma}
\begin{proof}
We offer two types of proof here. The first utilises the transformed bi-orthogonal
polynomials and associated functions given by (\ref{BOPSxfm+:a}-\ref{BOPSxfm+:e}) or 
(\ref{BOPSxfm-:a}-\ref{BOPSxfm-:e}) in the formula (\ref{An_resJ}). One then has to
separate out from each of the elements of the transformed $ A_n $ matrix the original
matrix and perform a matrix factorisation with the structure of a similarity 
transformation. This yields the final formulae (\ref{Axfm+:a},\ref{Axfm+:b},\ref{Axfm-:a},
\ref{Axfm-:b}). Alternatively one can start from the general Christoffel-Geronimus-Uvarov
transformation of $ A_n $ given by (\ref{K+1_Axfm}) or (\ref{L-1_Axfm}), insert
the partial fraction expansion (\ref{An_PFexp}) and resolve the resulting products
on the right-hand side of the formula in this way. In order to eliminate the unwanted
terms on the right-hand side we need to note the eigenvector identities
\begin{align*}
  A_{n,j} \begin{pmatrix} \phi_n(z_j) \cr \phi^*_n(z_j) \end{pmatrix}
  & = 0 \begin{pmatrix} \phi_n(z_j) \cr \phi^*_n(z_j) \end{pmatrix}, 
  \\ 
  A_{n,j} \begin{pmatrix} \xi_n(z_j) \cr -\xi^*_n(z_j) \end{pmatrix}
  & = -\rho_j \begin{pmatrix} \xi_n(z_j) \cr -\xi^*_n(z_j) \end{pmatrix}, 
  \\ 
  \begin{pmatrix} \phi^*_n(z_j) & -\phi_n(z_j) \end{pmatrix} A_{n,j}
  & =  -\rho_j \begin{pmatrix} \phi^*_n(z_j) & -\phi_n(z_j) \end{pmatrix},
  \\ 
 \begin{pmatrix} \xi^*_n(z_j) & \xi_n(z_j) \end{pmatrix} A_{n,j} 
  & = 0 \begin{pmatrix} \xi^*_n(z_j) & \xi_n(z_j) \end{pmatrix},
\end{align*}
and the form for $ A_{n,\infty} $ (\ref{An_resInfty}).
\end{proof}

Using the results gained up to this point we are now in a position to deduce the 
compatibility of the Schlesinger transformations with the spectral derivative
and deformation derivatives. We begin with the spectral derivative.
\begin{proposition}
The spectral matrix and the Schlesinger transformation matrices of the bi-orthogonal 
system satisfy the compatibility condition
\begin{equation} 
   \frac{d}{dz}R^{j\pm}_{n}(z;\rho_j)
  +R^{j\pm}_{n}(z;\rho_j)A_n(z;\rho_j) = A_n(z;\rho_j\pm1)R^{j\pm}_{n}(z;\rho_j) .
\label{Spectral+Schlesinger}
\end{equation}
\end{proposition}
\begin{proof}
Take the $ + $ case first. Using (\ref{R+}) we compute the difference between
the left-hand and right-hand sides of (\ref{Spectral+Schlesinger}) and resolve 
the expression into partial fractions in the $ z $ variable. The $ (z-z_j)^{-2} $ term 
has a coefficient proportional to 
\begin{equation}
         \begin{pmatrix}
           \phi_{n+1}(0)\phi^*_{n}(z_j)
	& -\phi_{n+1}(0)\phi_{n}(z_j)	\cr
          -z_j\kappa_{n+1}\phi^*_{n}(z_j)
	&  z_j\kappa_{n+1}\phi_{n}(z_j)	\cr
       \end{pmatrix} (I - A_{n,j})
  + A^{j+}_{n,j} 
         \begin{pmatrix}
           \phi_{n+1}(0)\phi^*_{n}(z_j)
	& -\phi_{n+1}(0)\phi_{n}(z_j)	\cr
          -z_j\kappa_{n+1}\phi^*_{n}(z_j)
	&  z_j\kappa_{n+1}\phi_{n}(z_j)	\cr
       \end{pmatrix} .
\label{1st_Id}
\end{equation}
Using (\ref{Axfm+:b}) this simplifies to
\begin{equation}
         \begin{pmatrix}
          -\phi_{n+1}(0)\phi^*_{n}(z_j)
	&  \phi_{n+1}(0)\phi_{n}(z_j)	\cr
           z_j\kappa_{n+1}\phi^*_{n}(z_j)
	& -z_j\kappa_{n+1}\phi_{n}(z_j)	\cr
       \end{pmatrix} (\rho_j + A_{n,j}) ,
\end{equation}
which, with the formula (\ref{An_resJ}), is identically zero. The nett coefficient
of the $ (z-z_j)^{-1} $ terms is proportional to 
\begin{multline}
       \begin{pmatrix}
         1 & 0 \cr
         \phi^*_{n}(z_j)/\phi_{n}(z_j) & 0 \cr
       \end{pmatrix} A_{n,j}
  - A^{j+}_{n,j} 
       \begin{pmatrix}
         1 & 0 \cr
         \phi^*_{n}(z_j)/\phi_{n}(z_j) & 0 \cr
       \end{pmatrix} \\
  + \sum_{k\neq j=0} \frac{\DySt 1}{\DySt z_j-z_k} \left[
       \begin{pmatrix}
          -r_{n+1}\phi^*_{n}(z_j)/\phi_{n}(z_j) & r_{n+1}	\cr
          z_j\phi^*_{n}(z_j)/\phi_{n}(z_j) &  -z_j	\cr
       \end{pmatrix} A_{n,k}
                   - A^{j+}_{n,k} 
       \begin{pmatrix}
          -r_{n+1}\phi^*_{n}(z_j)/\phi_{n}(z_j) & r_{n+1}	\cr
          z_j\phi^*_{n}(z_j)/\phi_{n}(z_j) &  -z_j	\cr
       \end{pmatrix} \right] .
\label{2nd_Id}
\end{multline}
Taking the summand of the above expression, substituting (\ref{Axfm+:a}) and
simplifying the result and then employing (\ref{An_resJ}) we arrive at
\begin{multline}
   \frac{\DySt \rho_k}{\DySt 2z_k^n}\frac{1}{\kappa_{n+1}\phi^2_n(z_j)}\Bigg[
  -    (z_j-z_k)\kappa_{n+1}\phi_{n}(z_k)\xi_{n}(z_k)
       \begin{pmatrix}
         \phi_{n}(z_j)\phi^*_{n}(z_j) & -\phi^2_{n}(z_j) \cr
         (\phi^*_{n}(z_j))^2 & -\phi_{n}(z_j)\phi^*_{n}(z_j) \cr
       \end{pmatrix}
  \\
  +    (\phi_{n}(z_j)\xi^*_{n}(z_k)+\phi^*_{n}(z_j)\xi_{n}(z_k))
       (\phi_{n}(z_j)\phi^*_{n}(z_k)-\phi^*_{n}(z_j)\phi_{n}(z_k))
       \begin{pmatrix}
         \phi_{n+1}(0) & 0 \cr
         -z_j\kappa_{n+1} & 0 \cr
       \end{pmatrix}
  \Bigg] .
\end{multline} 
When this is inserted back into (\ref{2nd_Id}) the first term is amenable to 
(\ref{An_sum:a}) and cancels with another term. Next (\ref{Axfm+:b})
is utilised in the simplified version of (\ref{2nd_Id}) and another term is 
annulled, leaving only two summation terms. Again appealing to (\ref{An_resJ})
we find the remaining two terms are equal and opposite in sign. Thus (\ref{2nd_Id})
is zero. 
Finally we are left with the $ (z-z_k)^{-1}, k\neq j $ terms whose
coefficient is proportional to
\begin{multline}
       \begin{pmatrix}
           \kappa_{n+1}\phi_n(z_j)+\frac{\DySt \phi_{n+1}(0)\phi^*_{n}(z_j)}{\DySt z_j-z_k}
	& -\frac{\DySt \phi_{n+1}(0)\phi_{n}(z_j)}{\DySt z_j-z_k}	\cr
          -z_k\frac{\DySt \kappa_{n+1}\phi^*_{n}(z_j)}{\DySt z_j-z_k}
	&  z_j\frac{\DySt \kappa_{n+1}\phi_{n}(z_j)}{\DySt z_j-z_k}	\cr
       \end{pmatrix} A_{n,k}
   \\ - A^{j+}_{n,k}
       \begin{pmatrix}
           \kappa_{n+1}\phi_n(z_j)+\frac{\DySt \phi_{n+1}(0)\phi^*_{n}(z_j)}{\DySt z_j-z_k}
	& -\frac{\DySt \phi_{n+1}(0)\phi_{n}(z_j)}{\DySt z_j-z_k}	\cr
          -z_k\frac{\DySt \kappa_{n+1}\phi^*_{n}(z_j)}{\DySt z_j-z_k}
	&  z_j\frac{\DySt \kappa_{n+1}\phi_{n}(z_j)}{\DySt z_j-z_k}	\cr
       \end{pmatrix} ,
\label{3rd_Id}
\end{multline}
and utilising (\ref{Axfm+:a}) we find this to vanish.

We now treat the $-$ case. Employing (\ref{R-}) and computing the difference of 
the left-hand and right-hand sides of (\ref{Spectral+Schlesinger}) we can resolve
it into partial fractions in $z$. The $ (z-z_j)^{-1} $ terms have a coefficient
proportional to 
\begin{equation}
   \begin{pmatrix}
      \kappa_{n}\xi^*_{n}(z_j) & \kappa_{n}\xi_{n}(z_j)	\cr
      \bar{\phi}_{n}(0)\xi^*_{n}(z_j) & \bar{\phi}_{n}(0)\xi_{n}(z_j) \cr
   \end{pmatrix} A_{n,j}
  -  A^{j-}_{n,j}
   \begin{pmatrix}
      \kappa_{n}\xi^*_{n}(z_j) & \kappa_{n}\xi_{n}(z_j)	\cr
      \bar{\phi}_{n}(0)\xi^*_{n}(z_j) & \bar{\phi}_{n}(0)\xi_{n}(z_j) \cr
   \end{pmatrix} .
\end{equation}
The first term is easily seen to be zero using (\ref{An_resJ}). In the second
term we employ (\ref{Axfm-:b}) and note that the matrix multiplications yield
a null result. 
The coefficient of the $ (z-z_k)^{-1} $ terms is proportional to 
\begin{multline}
   \begin{pmatrix}
      z_j\kappa_{n}\xi^*_{n}(z_j) & z_j\kappa_{n}\xi_{n}(z_j)	\cr
      z_k\bar{\phi}_{n}(0)\xi^*_{n}(z_j) & z_j\kappa_{n-1}\xi^*_{n-1}(z_j)-z_k\kappa_{n}\xi^*_{n}(z_j) \cr
   \end{pmatrix} A_{n,k}
  \\
  -  A^{j-}_{n,k}
   \begin{pmatrix}
      z_j\kappa_{n}\xi^*_{n}(z_j) & z_j\kappa_{n}\xi_{n}(z_j)	\cr
      z_k\bar{\phi}_{n}(0)\xi^*_{n}(z_j) & z_j\kappa_{n-1}\xi^*_{n-1}(z_j)-z_k\kappa_{n}\xi^*_{n}(z_j) \cr
   \end{pmatrix} .
\label{Id1}
\end{multline}
Substituting for $ A^{j-}_{n,k} $ from (\ref{Axfm-:a}) and carrying out the matrix 
multiplication we find this expression is identically zero. Lastly we examine the 
constant term of the partial fraction expansion, which simplifies to
\begin{equation}
   \begin{pmatrix}
      0 & 0 \cr
      \bar{\phi}_{n}(0) & -\kappa_{n} \cr
   \end{pmatrix} (I - A_{n,\infty})
  +  A^{j-}_{n,\infty}
   \begin{pmatrix}
      0 & 0 \cr
      \bar{\phi}_{n}(0) & -\kappa_{n} \cr
   \end{pmatrix} .
\label{Id2}
\end{equation}
This is easily shown to be zero using (\ref{An_resInfty}). This concludes the proof.
\end{proof} 

The final constraint is the compatibility of the Schlesinger transformations 
with the deformation derivative (\ref{YtDer}) which is established in the following
proposition.
\begin{proposition}
The deformation matrix and the Schlesinger transformation matrices of the bi-orthogonal 
system satisfy the compatibility condition
\begin{equation} 
   \frac{d}{dt}R^{j\pm}_{n}(z;\rho_j)
  +R^{j\pm}_{n}(z;\rho_j)B_n(z;\rho_j) = B_n(z;\rho_j\pm1)R^{j\pm}_{n}(z;\rho_j) .
\label{Deform+Schlesinger}
\end{equation}
\end{proposition}
\begin{proof}
Our strategy is to show that the difference between the left-hand and right-hand sides
of (\ref{Deform+Schlesinger}) is zero.
We begin with the $+$ case first.
Our initial task is to compute $ \frac{d}{dt}R^{j+}_n(z) $ starting with the 
expressions for $ R^{j+}_n(z) $ given by (\ref{R+}), (\ref{kappa+}) and employing
the derivatives (\ref{kappadot}), (\ref{rdot}) and (\ref{Pdot}). This yields the
following complicated expression (omitting an overall factor of $ \kappa^+_n/\kappa_n $)
\begin{multline}
  \frac{\dot{z}_j}{(z-z_j)^2}
   \begin{pmatrix}
      -r_{n+1}P_{n,j} & r_{n+1} \cr
      z_jP_{n,j} & -z_j \cr
   \end{pmatrix}
  \\
   + \frac{1}{z-z_j}\frac{1}{z_j+r_{n+1}P_{n,j}}\left\{
     \dot{z}_j
   \begin{pmatrix}
      -\frac{1}{2}[z-z_j-r_{n+1}P_{n,j}] & -\frac{1}{2}r_{n+1} \cr
      [-\frac{1}{2}z+z_j+r_{n+1}P_{n,j}]P_{n,j} & -\frac{1}{2}z_j-r_{n+1}P_{n,j} \cr
   \end{pmatrix}                                \right.
  \\
   +\sum_{k=1}\rho_k\frac{\dot{z}_k}{z_k} \left[
   \frac{\phi^*_n(z_k)\xi^*_n(z_k)}{2z^n_k}
   \begin{pmatrix}
      \frac{1}{2}r_{n+1}[z+z_j+r_{n+1}P_{n,j}] & \frac{1}{2}r^2_{n+1} \cr
      -z[z_j+\frac{1}{2}r_{n+1}P_{n,j}] & -\frac{1}{2}r_{n+1}z_j \cr
   \end{pmatrix}                         \right.
  \\
   +\frac{\phi_n(z_k)\xi^*_n(z_k)}{2z^n_k}
   [z_j+r_{n+1}P_{n,j}]
   \begin{pmatrix}
       0 & -r_{n+1} \cr
       zP_{n,j} & 0 \cr
   \end{pmatrix}   
  \\                                      \left.
   +\frac{\phi_n(z_k)\xi_n(z_k)}{2z^n_k}
   z_k   
   \begin{pmatrix}
      -\frac{1}{2}[z+z_j+r_{n+1}P_{n,j}]P_{n,j} & z_j+\frac{1}{2}r_{n+1}P_{n,j} \cr
      -\frac{1}{2}zP^2_{n,j} & \frac{1}{2}z_jP_{n,j} \cr
   \end{pmatrix}                          \right]
  \\
   +\sum_{k\neq j=0}\rho_k\frac{\dot{z}_j-\dot{z}_k}{z_j-z_k}
    \frac{[\phi_n(z_k)P_{n,j}-\phi^*_n(z_k)][\xi_n(z_k)P_{n,j}+\xi^*_n(z_k)]}{2z^n_k}
  \\ \times                                       \left.
   \begin{pmatrix}
      \frac{1}{2}r_{n+1}[z+z_j+r_{n+1}P_{n,j}] & \frac{1}{2}r^2_{n+1} \cr
      -z[z_j+\frac{1}{2}r_{n+1}P_{n,j}] & -\frac{1}{2}r_{n+1}z_j \cr
   \end{pmatrix}                                  \right\} .
\end{multline}
From the definition (\ref{YtDerInfty}) we can recast $ B_{n}(z;\rho_j) $ as
\begin{equation}
   \sum_{k=1}\rho_k\frac{\dot{z}_k}{z_k} \left[
   \frac{\phi^*_n(z_k)\xi^*_n(z_k)}{2z^n_k}
   \begin{pmatrix}
       0 & 0 \cr
      -1 & 0 \cr
   \end{pmatrix}
  +\frac{\phi_n(z_k)\xi^*_n(z_k)}{2z^n_k}
   \begin{pmatrix}
      -\frac{1}{2} & 0 \cr
       0 & \frac{1}{2} \cr
   \end{pmatrix}                         \right]
  -\sum_{k=1}A_{n,k}\frac{\dot{z}_k}{z-z_k} .
\label{Bn}
\end{equation}
In regard to the formula for $ B_{n}(z;\rho_j+1) $ we require deformation
derivatives of $ \kappa^{j+}_n $ and $ \bar{r}^{j+}_n $. We have found the former
in the above calculation whereas the latter can be computed starting with (\ref{refl+:b})
along with (\ref{Pdot}). In summary we find $ B_{n}(z;\rho_j+1) $ to be
\begin{multline}  
  \frac{\dot{z}_j}{z_j+r_{n+1}P_{n,j}}
   \begin{pmatrix}
       -\frac{1}{2} & 0 \cr
       -P_{n,j} & \frac{1}{2} \cr
   \end{pmatrix}  
  \\
   +\sum_{k=1}\rho_k\frac{\dot{z}_k}{z_k} \left[
    \frac{\phi^*_n(z_k)\xi^*_n(z_k)}{2z^n_k}
    \frac{1}{z_j+r_{n+1}P_{n,j}}
   \begin{pmatrix}
      \frac{1}{2}r_{n+1} & 0 \cr
      -z_j & -\frac{1}{2}r_{n+1} \cr
   \end{pmatrix}                         \right.
  \\
   +\frac{\phi_n(z_k)\xi^*_n(z_k)}{2z^n_k}
   \begin{pmatrix}
       -\frac{1}{2} & 0 \cr
        0 & \frac{1}{2} \cr
   \end{pmatrix}
  \\                                     \left.
   +\frac{\phi_n(z_k)\xi_n(z_k)}{2z^n_k}
    \frac{z_kP_{n,j}}{z_j+r_{n+1}P_{n,j}}
   \begin{pmatrix}
       -\frac{1}{2} & 0 \cr
       -P_{n,j} & \frac{1}{2} \cr
   \end{pmatrix}                         \right]
  \\
   +\sum_{k\neq j=0}\rho_k\frac{\dot{z}_j-\dot{z}_k}{z_j-z_k}
    \frac{[\phi_n(z_k)P_{n,j}-\phi^*_n(z_k)][\xi_n(z_k)P_{n,j}+\xi^*_n(z_k)]}{2z^n_k}
  \\ \times
    \frac{1}{z_j+r_{n+1}P_{n,j}}
   \begin{pmatrix}
      \frac{1}{2}r_{n+1} & 0 \cr
      -z_j & -\frac{1}{2}r_{n+1} \cr
   \end{pmatrix}
  \\
  -\sum_{k=1}A^{j+}_{n,k}\frac{\dot{z}_k}{z-z_k} .
\end{multline}  
Having all the terms in (\ref{Deform+Schlesinger}) we now perform a partial 
fraction decomposition with respect to $ z $ and examine the resulting coefficients
in this decomposition. Collecting the $ 1/(z-z_k) $ terms with $ k\neq j $ we
find its nett coefficient is proportional to (\ref{3rd_Id}) in the preceding
proposition and thus vanishes. Considering the $ 1/(z-z_j)^2 $ terms we find the
resulting coefficient is proportional to (\ref{1st_Id}) and also vanishes 
according to the preceding proposition. Assembling all the constant terms and
carrying out the matrix multiplications we find the single terms and summands
cancel internally without recourse to identities. This leaves the $ 1/(z-z_j) $
terms which we will break down into further sub-divisions. Of these the terms
with the factor $ \dot{z}_j/[z_j+r_{n+1}P_{n,j}] $ (excluding those proportional
to the monodromy exponent $ n-\rho_0 $) are easily seen to cancel. 
In the summations over $ A_{n,k} $ and $ A^{j+}_{n,k} $ we make the substitution
$ \dot{z}_k = \dot{z}_j-(\dot{z}_j-\dot{z}_k) $. We then collect the terms 
with the factor of $ \dot{z}_j $ resulting from this and all the remaining terms 
with this factor. Using (\ref{2nd_Id}) these terms are simplified to 
(they do not vanish in totality because the sum over $ k $ excludes $ k=0 $)
\begin{multline}
  -\frac{\dot{z}_j}{z_j}\left[
  -\begin{pmatrix}
      -r_{n+1}P_{n,j} & r_{n+1} \cr
      z_jP_{n,j} & -z_j \cr
   \end{pmatrix} A_{n,0}
  +A^{j+}_{n,0}
   \begin{pmatrix}
      -r_{n+1}P_{n,j} & r_{n+1} \cr
      z_jP_{n,j} & -z_j \cr
   \end{pmatrix}        \right] 
 \\
  = \frac{\dot{z}_j}{z_j}(n-\rho_0)
   \begin{pmatrix}
      r_{n}z_jP_{n,j}-r_{n+1}P_{n,j}(1-r_{n}P_{n,j}) & -r_{n}z_j \cr
      z_jP_{n,j} & -r_{n}z_jP_{n,j} \cr
   \end{pmatrix} .
\label{4th_Id}
\end{multline}
The $ k=0 $ term is extracted from the summations containing $ \dot{z}_j-\dot{z}_k $
and simplifies to 
\begin{equation} 
  \frac{\dot{z}_j}{z_j}(n-\rho_0)(1-r_{n}P_{n,j})P_{n,j}
   \begin{pmatrix}
      r_{n+1} & 0 \cr
      -z_j & 0 \cr
   \end{pmatrix} .
\label{5th_Id}  
\end{equation} 
Adding (\ref{4th_Id}) and (\ref{5th_Id}) we get
\begin{equation} 
  \dot{z}_j(n-\rho_0)r_{n}
   \begin{pmatrix}
      P_{n,j} & -1 \cr
      P^2_{n,j} & -P_{n,j} \cr
   \end{pmatrix} .
\label{6th_Id}  
\end{equation}
Now consider the $ k\neq 0 $ terms in the summations containing $ \dot{z}_j-\dot{z}_k $.
Using the working leading to (\ref{5th_Id}) and the formulae (\ref{Axfm+:a}) and (\ref{An_resJ})
we find the summands reduce to
\begin{equation} 
  (\dot{z}_j-\dot{z}_k)(A_{n,k})_{1,2}
   \begin{pmatrix}
      -P_{n,j} & 1 \cr
      -P^2_{n,j} & P_{n,j} \cr
   \end{pmatrix} .
\label{7th_Id}   
\end{equation}
Combining the contributions of (\ref{6th_Id}) and (\ref{7th_Id}) we have 
\begin{equation} 
  \sum_{k\neq j=0}(\dot{z}_j-\dot{z}_k)(A_{n,k})_{1,2}
   \begin{pmatrix}
      -P_{n,j} & 1 \cr
      -P^2_{n,j} & P_{n,j} \cr
   \end{pmatrix} .
\label{8th_Id}   
\end{equation}
There is one group of terms in the $ 1/(z-z_j) $ coefficient to be simplified, namely the 
summation with the factor $ \dot{z}_k/z_k $. After considerable internal cancellation 
this group simplifies to
\begin{equation} 
  \sum_{k=1}\dot{z}_k(A_{n,k})_{1,2}
   \begin{pmatrix}
      -P_{n,j} & 1 \cr
      -P^2_{n,j} & P_{n,j} \cr
   \end{pmatrix} .
\label{9th_Id}   
\end{equation}
In summary the $ 1/(z-z_j) $ coefficient reduces to
\begin{equation} 
  \dot{z}_j\sum_{k=0}(A_{n,k})_{1,2}
   \begin{pmatrix}
      -P_{n,j} & 1 \cr
      -P^2_{n,j} & P_{n,j} \cr
   \end{pmatrix} 
 \propto -\dot{z}_j(A_{n,\infty})_{1,2} = 0 ,
\end{equation}
where we have used (\ref{An_resInfty}) in the last step.

Next we treat the $-$ case. We recall the definition of $ Q_{n,j}:= \xi_n(z_j)/\xi^*_n(z_j) $.
We compute $ \frac{d}{dt}R^{j\pm}_{n}(z;\rho_j) $ first beginning with (\ref{R-})
and (\ref{kappa-}) and utilising the derivatives (\ref{rCdot}) and (\ref{Qdot}). 
Up to an overall factor of $ \kappa_n/\kappa^-_n $ we find this term is
\begin{multline}
  -\frac{\dot{z}_j}{2z_j}
   \begin{pmatrix}
      z_j & z_jQ_{n,j} \cr
      -\bar{r}_{n}z & z+z_j(1+\bar{r}_{n}Q_{n,j}) \cr
   \end{pmatrix}
  \\
  -\sum_{k=1}\rho_k\frac{\dot{z}_k}{z_k} \left[
   \frac{\phi^*_n(z_k)\xi^*_n(z_k)}{2z^n_k}
   \begin{pmatrix}
      \frac{1}{2}z_jQ_{n,j} & -\frac{1}{2}z_jQ^2_{n,j} \cr
      z[-1+\frac{1}{2}\bar{r}_nQ_{n,j}] & -\frac{1}{2}Q_{n,j}[z+z_j(1+\bar{r}_{n}Q_{n,j})] \cr
   \end{pmatrix}
                                         \right.
  \\                                     \left.
  +\frac{\phi_n(z_k)\xi^*_n(z_k)}{2z^n_k}
   \begin{pmatrix}
      0 & -z_jQ_{n,j} \cr
      z\bar{r}_n & 0 \cr
   \end{pmatrix}                         \right]
  \\
  -\sum_{k\neq j=0}\rho_k\frac{\dot{z}_j-\dot{z}_k}{z_j-z_k}
   \frac{[\phi_n(z_k)+Q_{n,j}\phi^*_n(z_k)][Q_{n,j}\xi^*_n(z_k)-\xi_n(z_k)]}{2z^n_k}
   \frac{1}{1+\bar{r}_{n}Q_{n,j}}
  \\ \times
   \begin{pmatrix}
      -\frac{1}{2}z_j\bar{r}_n & z_j[1+\frac{1}{2}\bar{r}_nQ_{n,j}] \cr
      -\frac{1}{2}\bar{r}^2_n z & \frac{1}{2}\bar{r}_n[z+z_j(1+\bar{r}_{n}Q_{n,j})] \cr
   \end{pmatrix} .
\end{multline}
However in contrast to the formula for $ B_{n}(z;\rho_j) $ (\ref{Bn}) the one for 
$ B_{n}(z;\rho_j-1) $ is more complicated as we have to compute the 
deformation derivatives of the shifted quantities $ \kappa^-_n $, $ \bar{r}^-_n $.
Starting with the expressions (\ref{kappa-}) and (\ref{refl-:b}) we can employ
the derivatives computed earlier in addition to (\ref{kappadot}). In summary we
find $ B_{n}(z;\rho_j-1) $ to be given by
\begin{multline}
   \frac{\dot{z}_j}{z_j}
   \begin{pmatrix}
      \frac{1}{2} & 0 \cr
      \bar{r}_{n} & -\frac{1}{2} \cr
   \end{pmatrix}
  \\
  +\sum_{k=1}\rho_k\frac{\dot{z}_k}{z_k} \left[
   \frac{[\phi^*_n(z_k)-\bar{r}_n\phi_n(z_k)][\xi^*_n(z_k)+\bar{r}_n\xi_n(z_k)]}{2z^n_k}
   \frac{z_k}{z_j[1+\bar{r}_{n}Q_{n,j}]}
   \begin{pmatrix}
      0 & 0 \cr
      1 & 0 \cr
   \end{pmatrix}                         \right.
  \\                                     \left.
  +\frac{\phi^*_n(z_k)\xi^*_n(z_k)}{2z^n_k}
   \begin{pmatrix}
      -\frac{1}{2}Q_{n,j} & 0 \cr
      -[1+\bar{r}_nQ_{n,j}] & \frac{1}{2}Q_{n,j} \cr
   \end{pmatrix}
  +\frac{\phi_n(z_k)\xi^*_n(z_k)}{2z^n_k}
   \begin{pmatrix}
      -\frac{1}{2} & 0 \cr
      0 & \frac{1}{2} \cr
   \end{pmatrix}                         \right]
  \\
  +\sum_{k\neq j=0}\rho_k\frac{\dot{z}_j-\dot{z}_k}{z_j-z_k}
   \frac{[\phi_n(z_k)+Q_{n,j}\phi^*_n(z_k)][Q_{n,j}\xi^*_n(z_k)-\xi_n(z_k)]}{2z^n_k}
   \frac{\bar{r}_n}{1+\bar{r}_{n}Q_{n,j}}
   \begin{pmatrix}
      \frac{1}{2} & 0 \cr
      \bar{r}_n & -\frac{1}{2} \cr
   \end{pmatrix}
  \\
  -\sum_{k=1}A^{j-}_{n,k}\frac{\dot{z}_k}{z-z_k} .
\end{multline}
We now decompose (\ref{Deform+Schlesinger}) into partial fractions with respect to
$ z $. Taking the $ 1/(z-z_k) $ terms we find they are proportional to the 
expression (\ref{Id1}) which was shown earlier to vanish. Collecting the $ z $ 
terms together one finds they vanish identically in a straight forward manner. 
This leaves the constant terms. After making the obvious simplifications here one 
is left with
\begin{multline}
  \dot{z}_j
   \begin{pmatrix}
      0 & 0 \cr
      \bar{r}_n & -1 \cr
   \end{pmatrix}
 + \begin{pmatrix}
      0 & 0 \cr
      \bar{r}_n & -1 \cr
   \end{pmatrix} 
   \sum_{k=1}A_{n,k}\dot{z}_k
 - \sum_{k=1}A^{j-}_{n,k}\dot{z}_k
   \begin{pmatrix}
      0 & 0 \cr
      \bar{r}_n & -1 \cr
   \end{pmatrix}
 \\
 + \frac{1}{1+\bar{r}_{n}Q_{n,j}}
   \begin{pmatrix}
      0 & 0 \cr
      1 & Q_{n,j} \cr
   \end{pmatrix}
   \sum_{k=1}\rho_k\dot{z}_k\frac{[\xi^*_n(z_k)+\bar{r}_n\xi_n(z_k)][\phi^*_n(z_k)-\bar{r}_n\phi_n(z_k)]}{2z^n_k}
\\
 + \frac{z_j}{1+\bar{r}_{n}Q_{n,j}}
   \begin{pmatrix}
      \bar{r}_n & -1 \cr
      \bar{r}^2_n & -\bar{r}_n \cr
   \end{pmatrix}
   \sum_{k\neq j=0}\rho_k\frac{\dot{z}_j-\dot{z}_k}{z_j-z_k}
   \frac{[\phi_n(z_k)+Q_{n,j}\phi^*_n(z_k)][Q_{n,j}\xi^*_n(z_k)-\xi_n(z_k)]}{2z^n_k} .
\label{Id3}
\end{multline}
In the second and third terms we substitute $ \dot{z}_k = \dot{z}_j - (\dot{z}_j-\dot{z}_k) $
and split them into two terms each. Collecting the terms with a factor of $ \dot{z}_j $
we can perform the $ k $ summation using (\ref{An_resInfty}). This group of terms vanishes
because of (\ref{Id2}). In the remaining group we employ (\ref{Axfm-:a}) and (\ref{Axfm-:b})
along with (\ref{An_resJ}) and carry out the matrix multiplications yielding
\begin{multline}
 - \sum_{k=0}\rho_k(\dot{z}_j-\dot{z}_k)\frac{[\xi^*_n(z_k)+\bar{r}_n\xi_n(z_k)]}{2z^n_k}
   \begin{pmatrix}
      0 & 0 \cr
      -\phi^*_n(z_k) & \phi_n(z_k) \cr
   \end{pmatrix}
 \\
 - \frac{z_j}{1+\bar{r}_{n}Q_{n,j}}
   \sum_{k\neq j=0}\rho_k\frac{\dot{z}_j-\dot{z}_k}{z_j-z_k}
   \frac{[\phi_n(z_k)+Q_{n,j}\phi^*_n(z_k)][Q_{n,j}\xi^*_n(z_k)-\xi_n(z_k)]}{2z^n_k} 
   \begin{pmatrix}
      \bar{r}_n & -1 \cr
      \bar{r}^2_n & -\bar{r}_n \cr
   \end{pmatrix}
 \\
 - \frac{1}{1+\bar{r}_{n}Q_{n,j}}
   \sum_{k=0}\rho_k(\dot{z}_j-\dot{z}_k)
   \frac{[\phi_n(z_k)+Q_{n,j}\phi^*_n(z_k)][\xi^*_n(z_k)+\bar{r}_n\xi_n(z_k)]}{2z^n_k} 
   \begin{pmatrix}
      0 & 0 \cr
      \bar{r}_n & -1 \cr
   \end{pmatrix} .
\end{multline}
The second term of the above expression cancels with the last term of (\ref{Id3}).
In what remains of the constant terms we separate them into those with a factor of
$ \dot{z}_j $ and those with $ \dot{z}_k $. In the former set of terms we can 
perform the $ k $ summation using (\ref{An_sum:a}-\ref{An_sum:d}) and find complete
cancellation. In the latter set we find the summands mutually cancel. Thus all the
constant terms vanish identically and the theorem is proved.
\end{proof}

It is expected that the explicit Schlesinger transformations derived here 
(\ref{R+}), (\ref{R-}), (\ref{Axfm+:a}), (\ref{Axfm+:b}), (\ref{Axfm-:a}), (\ref{Axfm-:b})
should co-incide with those found by the Kyoto school \cite{JM_1981b} however we 
refrain from carrying out this verification as this would entail a significant 
enlargement of the tasks considered here and add substantially to the length of the 
present study.

We conclude this section and our study with a selection of results on the bi-linear 
difference equations satisfied by the $ \tau $-functions of our integrable system, i.e. the 
Toeplitz determinants. These are pure difference equations in the parameters of the 
system, $ n, \{\rho_j\}_{j=0,\ldots,M} $ or the monodromy exponents $ \{\theta_j\}_{j=0,l\dots,M,\infty} $,
for products of $ \tau $-functions,
and are known as Hirota-Miwa equations. So in this sense we return to one of the initial
motivations of our study, namely to find alternative means of characterising 
averages over the unitary group (\ref{Haar_avge}) of the regular semi-classical weights (\ref{gJwgt}). 
All of the Hirota-Miwa type equations are very easily derived from the preceding analysis
and one might consider them to be immediately apparent. 
Let us define the following notation for the Toeplitz determinants
\begin{equation}
    I_{n}[\{\rho_j\}] \coloneqq \tau(\theta_0,\ldots,\theta_j,\ldots,\theta_{\infty}) .
\end{equation}

\begin{proposition}
For $ j=1,\ldots,M $ (the $ j=0 $ cases are trivial)
the $ \tau $-functions satisfy the Hirota-Miwa equations
\begin{multline}
   \tau()\tau(\theta_0+1,\theta_{\infty}+1) \\
  = \tau(\theta_j+1,\theta_{\infty}+1)\tau(\theta_0+1,\theta_{j}-1)-z_j\tau(\theta_j-1,\theta_{\infty}+1)\tau(\theta_0+1,\theta_j+1) ,
\label{HMtau:a}
\end{multline}
and
\begin{multline}
  \tau()\tau(\theta_0-1,\theta_{\infty}+1)   \\
  = \tau(\theta_j-1,\theta_{\infty}+1)\tau(\theta_0-1,\theta_{j}+1)+z_j\tau(\theta_j+1,\theta_{\infty}+1)\tau(\theta_0-1,\theta_j-1) ,
\label{HMtau:b}
\end{multline}
where $ \tau() $ denotes the $ \tau$-function with unshifted monodromy exponents.
These constitute all of the independent equations involving one generic exponent $ \theta_j $.
\end{proposition}
\begin{proof}
These follow from identifying the bi-orthogonal polynomials and associated functions
in terms of modified Toeplitz determinants using (\ref{IntRep:a})-(\ref{IntRep:d})
and employing these formulae in the Casoratian equations (\ref{Cas:a})-(\ref{Cas:c}).
We find that
\begin{align*}
   I_{n}I_{n+1} & = I_{n}^{*j+}I_{n+1}^{*j-}-z_jI_{n}^{j+}I_{n+1}^{j-} ,
  \\
   (-)^nr_{n}I_{n}^2 & = I_{n}^{j+}I_{n}^{*j-}+z_jI_{n-1}^{j+}I_{n+1}^{*j-} ,
  \\
   (-)^n\bar{r}_{n}I_{n}^2 & = I_{n}^{*j+}I_{n}^{j-}+z_jI_{n-1}^{*j+}I_{n+1}^{j-} ,
\end{align*}
and we note that these last two equations lead to the same bi-linear equation. 
Consequently we deduce (\ref{HMtau:a}) and (\ref{HMtau:b}).
\end{proof}

\begin{proposition}
For any $ j\neq k=0,\ldots, M $
the $ \tau $-functions satisfy the Hirota-Miwa equations
\begin{multline}
  (z_j-z_k)\tau(\theta_0+1,\theta_{\infty}+1)\tau(\theta_j-1,\theta_k-1,\theta_{\infty}+2) \\
  = \tau(\theta_0+1,\theta_k-1,\theta_{\infty}+2)\tau(\theta_j-1,\theta_{\infty}+1)-\tau(\theta_0+1,\theta_j-1,\theta_{\infty}+2)\tau(\theta_k-1,\theta_{\infty}+1) ,
\label{HMtau:d}
\end{multline}
and
\begin{multline}
    \tau(\theta_0+1,\theta_{\infty}+1)\tau(\theta_j-1,\theta_{k}+1) \\
  = -z_k \tau(\theta_0+1,\theta_k+1)\tau(\theta_j-1,\theta_{\infty}+1)+\tau(\theta_k+1,\theta_{\infty}+1)\tau(\theta_0+1,\theta_{j}-1) ,
\label{HMtau:e}
\end{multline}
and
\begin{multline}
  (z_j-z_k)\tau()\tau(\theta_0+1,\theta_j-1,\theta_{k}-1,\theta_{\infty}+1) \\
  = z_j \tau(\theta_j-1,\theta_{\infty}+1)\tau(\theta_0+1,\theta_k-1)-z_k \tau(\theta_0+1,\theta_j-1)\tau(\theta_k-1,\theta_{\infty}+1) ,
\label{HMtau:f}
\end{multline}
and
\begin{multline}
  \tau()\tau(\theta_0+1,\theta_j-1,\theta_k+1,\theta_{\infty}-1) \\
  = z_j\tau(\theta_j-1,\theta_{\infty}-1)\tau(\theta_0+1,\theta_k+1)+\tau(\theta_0+1,\theta_j-1)\tau(\theta_k+1,\theta_{\infty}-1) .
\label{HMtau:g}
\end{multline}
These constitute all of the independent equations involving two generic exponents $ \theta_j, \theta_k $.
\end{proposition}
\begin{proof}
Equation (\ref{HMtau:d}) is derived by combining the expression for the double shift 
$ \rho_j \mapsto \rho_j+1, \rho_k \mapsto \rho_k+1 $ given in (\ref{kappa_jkUp})
with the single shifts $ \rho_{j,k} \mapsto \rho_{j,k}+1 $ as in (\ref{kappa+})
and this results in
\begin{equation}
   I_{n+1}I_{n}^{j+,k+}\left[ I_{n+2}^{k+}I_{n+1}^{j+}-I_{n+2}^{j+}I_{n+1}^{k+} \right]
   = I_{n+2}I_{n+1}^{j+k+}\left[ I_{n+1}^{k+}I_{n}^{j+}-I_{n+1}^{j+}I_{n}^{k+} \right]
\end{equation}
This can rewritten as a perfect difference equation and solved, yielding 
\begin{equation*}
   I_{n}^{j+}I_{n+1}^{k+} - I_{n}^{k+}I_{n+1}^{j+} = (z_j-z_k)I_{n+1}I_{n}^{j+,k+} .
\end{equation*}
A similar method applies when deriving the other equations and we confine ourselves 
to recording the respective intermediate steps, namely
\begin{equation*}
   I_{n}^{*j+}I_{n+1}^{*k-} - z_kI_{n}^{j+}I_{n+1}^{k-} = I_{n+1}I_{n}^{j+,k-} ,
\end{equation*}
and
\begin{equation*}
   z_jI_{n}^{j+}I_{n}^{*k+} - z_kI_{n}^{*j+}I_{n}^{k+} = (z_j-z_k)I_{n}I_{n}^{*j+,k+} ,
\end{equation*}
and
\begin{equation*}
   z_jI_{n-1}^{*j+}I_{n+1}^{k-} + I_{n}^{*j+}I_{n}^{k-} = I_{n}I_{n}^{*j+,k-} .
\end{equation*}
\end{proof}

\section{Acknowledgments}
This research has been supported by the Australian Research Council.
The author wishes to acknowledge the generous advice given by Peter Forrester,
the interest shown by Mourad Ismail and the detailed attention given to the work
by an anonymous referee.

\setcounter{equation}{0}
\bibliographystyle{amsplain}
\bibliography{moment,nonlinear,random_matrices}

\end{document}